\documentclass[reqno,12pt]{amsart}

\usepackage{amsmath}
\usepackage{aliascnt}
\usepackage{amsfonts}
\usepackage{bm}
\usepackage{amsthm}
\usepackage[dvipsnames, svgnames]{xcolor}
\usepackage{dsfont}
\usepackage{graphicx}
\usepackage{float}
\usepackage[export]{adjustbox}
\usepackage{fancyref}
\usepackage{algorithm2e}
\usepackage[hidelinks]{hyperref}
\usepackage[colorinlistoftodos,prependcaption,textsize=small]{todonotes}
\usepackage[title,titletoc,toc]{appendix}
\usepackage{mathrsfs}
\usepackage{microtype}
\usepackage{enumerate} 
\usepackage{listings}
\usepackage{bbm}
\usepackage{tikz-cd} 
\usetikzlibrary{shapes,arrows,positioning,chains}\usepackage[a4paper,bindingoffset=0.5cm,left=1.9cm,right=1.9cm,top=2.5cm,bottom=2cm,footskip=.8cm]{geometry}
\usepackage{mathabx} 
\usepackage{pgfplots}
\usepgfplotslibrary{dateplot, statistics}
\usepackage{subcaption}

\pgfplotsset{compat=newest}
\usepackage{tikzexternal}
\definecolor{color11}{HTML}{aedea7}
\definecolor{color12}{HTML}{74c476}
\definecolor{color13}{HTML}{37a055}
\definecolor{color14}{HTML}{c6c6c6}
\definecolor{color15}{HTML}{969696}
\definecolor{color16}{HTML}{686868}

\definecolor{color31}{HTML}{bce4b5}
\definecolor{color32}{HTML}{8ed08b}
\definecolor{color33}{HTML}{57b668}
\definecolor{color34}{HTML}{2c944c}
\definecolor{color35}{HTML}{057130}
\definecolor{color36}{HTML}{00441b}

\renewcommand\sectionautorefname{Section}
\renewcommand\figureautorefname{Figure}

\newaliascnt{lemma}{theorem}

\aliascntresetthe{lemma}
\providecommand*{\lemmaautorefname}{Lemma}

\newaliascnt{corollary}{theorem}

\aliascntresetthe{corollary}
\providecommand*{\corollaryautorefname}{Corollary}

\newaliascnt{proposition}{theorem}

\aliascntresetthe{proposition}
\providecommand*{\propositionautorefname}{Proposition}

\newaliascnt{remark}{theorem}

\aliascntresetthe{remark}
\providecommand*{\remarkautorefname}{Remark}

\newaliascnt{assumption}{theorem}

\aliascntresetthe{assumption}
\providecommand*{\assumptionautorefname}{Assumption}

\newaliascnt{definition}{theorem}

\aliascntresetthe{definition}
\providecommand*{\definitionautorefname}{Definition}

\newcommand{\1}{\mathds{1}}
\newcommand{\one}[1]{\1_{\{#1\}}}
\newcommand{\comment}[1]{\footnote{\red{#1}}}
\newcommand{\contentline}[1]{[ $\rhd$ \textit{#1} $\lhd$]}

\newcommand{\R}{\mathbbm{R}}
\newcommand{\C}{\mathbbm{C}}
\newcommand{\N}{\mathbbm{N}}
\newcommand{\PP}{\mathbbm{P}}
\newcommand{\E}{\mathbbm{E}}
\newcommand{\Q}{\mathbbm{Q}}
\newcommand{\Z}{\mathbbm{Z}}
\newcommand{\F}{\mathbbm{F}}

\newcommand{\rhoh}{\hat{\rho}}
\newcommand{\rhob}{\bar{\rho}}

\newcommand{\cI}{\mathcal{I}}
\newcommand{\cJ}{\mathcal{J}}
\newcommand{\cN}{\mathcal{N}}
\newcommand{\cE}{\mathcal{E}}
\newcommand{\cF}{\mathcal{F}}
\newcommand{\cB}{\mathcal{B}}
\newcommand{\cX}{\mathcal{X}}
\newcommand{\cL}{\mathcal{L}}
\newcommand{\cS}{\mathcal{S}}
\newcommand{\cM}{\mathcal{M}}

\newcommand*\diff{\mathop{}\!\mathrm{d}}
\newcommand{\Var}{\mathrm{Var}}
\newcommand{\Cov}{\mathrm{Cov}}

\renewcommand{\labelitemi}{$\circ$}

\newcommand{\delims}[4]{\mathopen#1#2#4\mathclose#1#3}
\newcommand{\abs}[2][]{\delims{#1}\lvert\rvert{#2}}		
\newcommand{\norm}[2][]{\delims{#1}\lVert\rVert{#2}}	
\newcommand{\floor}[2][]{\delims{#1}\lfloor\rfloor{#2}}	
\newcommand{\ceil}[2][]{\delims{#1}\lceil\rceil{#2}}	
\newcommand{\red}[1]{\textcolor{red}{#1}}		

\newcommand{\dx}{\,\mathrm{d}x}
\newcommand{\xm}{\smash[t]{x_L^-}}

\DeclareMathOperator{\Ex}{E}
\DeclareMathOperator{\D10}{D_{1 \to 0}}
\DeclareMathOperator{\D21}{D_{2 \to 1}}
\DeclareMathOperator{\D32}{D_{3 \to 2}}

\renewcommand{\baselinestretch}{1.099}
\allowdisplaybreaks

\newcommand\restr[2]{{
  \left.\kern-\nulldelimiterspace 
  #1 
  \vphantom{\big|} 
  \right|_{#2} 
  }}

\let\lim\relax
\let\inf\relax
\DeclareMathOperator*{\lim}{\vphantom{p}lim}
\DeclareMathOperator*{\inf}{\vphantom{p}inf}

\let\subsectionautorefname\sectionautorefname
\let\subsubsectionautorefname\sectionautorefname

\definecolor{light-gray}{gray}{0.65}
\definecolor{darkgreen}{HTML}{44CC44}

\title[A Gaussian segment-based traffic flow model]{A Gaussian segment-based traffic flow model for the design and control of transport networks}
\author{Michel Mandjes and Jaap Storm}
\date{\today}

\begin{document}

\begin{abstract}

In the setting of a recently developed cellular stochastic traffic flow model, it has shown that the joint per-cell vehicle densities, as a function of time, can be accurately approximated by a Gaussian process, which has the attractive feature that its means and (spatial and temporal) covariances can be efficiently evaluated. 
The present article demonstrates the rich potential of this methodology in the context of road traffic control and transportation network design. To solidly provide empirical backing for the use of a multivariate Gaussian approximation, we rely on a detailed historical dataset that contains traffic flow data. Then, in the remainder of the paper, we provide a sequence of design and control related example questions that can be analyzed using the Gaussian methodology. These cover the following topics: (i)~{evaluation of stationary performance measures}, (ii)~{route selection}, (iii)~{control of traffic flows}, and {(iv)~{performance of traffic networks with arbitrary topology}}. In discussing the setup, results, and applications of these examples, we stress the appropriateness of our {stochastic} traffic model over a deterministic counterpart.

\vspace{3mm}


\vspace{3mm}

\noindent
{\sc Affiliations.} 
Michel Mandjes ({\tt \footnotesize m.r.h.mandjes@uva.nl}) is with Korteweg-de Vries Institute for Mathematics, University of Amsterdam, Science Park 904, 1098 XH Amsterdam, the Netherlands. He is also with E{\sc urandom}, Eindhoven University of Technology, Eindhoven, the Netherlands, and Amsterdam Business School, Faculty of Economics and Business, University of Amsterdam, Amsterdam, the Netherlands. 

\noindent Jaap Storm ({\tt \footnotesize p.j.storm@vu.nl}) is with {the} Department of Mathematics, Vrije Universiteit Amsterdam,
De Boelelaan {1111}, 1081 HV Amsterdam, the Netherlands. Their research is partly funded by NWO Gravitation project N{\sc etworks}, grant number 024.002.003.

\end{abstract}

\maketitle

\section{Introduction}

So as to develop an efficiently operating transport network, one relies on models that aim at describing the interaction between travelers and the underlying infrastructure. The efficiency of the network is quantified in terms of performance metrics such as the throughput of roads (i.e., the volume of vehicles that are sent over a specific road segment per time unit) or the travel time between specific origin-destination pairs. 

It is evident that a well functioning transport network has a favorable effect on society, in terms of e.g., economic growth and sustainability. This explains why one wishes to use the models mentioned above to quantify the impact of infrastructural changes, typical questions in this respect being `How much does the throughput increase by opening an additional lane?' and `By how much do travel times go down after building a new highway?' In addition, the models are also a useful tool when  developing mechanisms that effectively control streams of vehicles, with typical questions in the spirit of `What is the impact on the congestion level if one adapts the maximum speed to 100 km/h' and `Does distributing the load between a certain origin and destination over multiple routes improve travel times?' Also from an environmental point of view, there are compelling reasons to rely on traffic flow models, for instance when one aims at striking a proper balance between efficiency (say, in terms of travel times or throughput), fuel consumption, and air pollution. In this paper, we illustrate how a specific methodology that has been developed earlier can be used to solve the various control and design problems above. We proceed by providing some more background on the model and the literature.

\vspace{2mm}

Traffic {flow} models come in many varieties, each of them focusing on specific aspects{; see, e.g., \cite{Hoogendoorn2001} for a historic overview of traffic flow models developed between 1950--2000, and \cite{WageningenKessels2015genealogy} for a more recent account. It is customary to classify these models in terms of the level of detail at which traffic flows are described, and in addition the timescale considered.} 
On the one hand there are macroscopic models, which do not distinguish individual vehicles, and work with concepts like the density, flow, and velocity of the underlying traffic stream. Microscopic models, on the other hand,  describe the dynamics of individual vehicles. Existing models are predominantly of a deterministic nature, for instance the macroscopic models that aim at describing vehicles as continuous flows using conservation laws; see {e.g.,} \cite{lw1955} for a seminal contribution.
As argued extensively in \cite{QU2017}, besides such physical laws, various microscopic variables play a crucial role when modeling traffic flows. In this respect one could think of the different perceptions, moods, and driving habits that individual car drivers may have. This led to the consensus (as discussed in detail in {e.g.,} \cite[Section 1]{QU2017}) that such microscopic variables should be represented by random variables. Consequently, in order to accurately describe streams of vehicles, one should work with stochastic traffic flow models. 

In \cite{MS2019}, a stochastic traffic flow model was developed that contributed in closing the gap between deterministic, tractable models, and stochastic, intractable models. Following a line of research that started with \cite{JL2012,JL2013}, a road network, consisting of disjoint segments (referred to as cells) is considered, with the aim of finding a probabilistic description for the joint vehicle density of the cells, as a function of time. The framework developed in \cite{MS2019} has a high level of generality. It in particular explicitly covers vehicle streams consisting of multiple vehicle types, which enables modeling the different behavior and impact of, say, passenger cars and trucks, in line with the framework developed in, e.g., \cite{BC,cb2003,li2008,ww2002}. Second, in addition to evaluating standard performance indicators such as congestion levels and vehicular flows, in \cite{MS2019} a method was developed to accurately approximate the travel-time distribution between any origin and destination in the network, for any individual vehicle class.

{The first main finding of \cite{MS2019} is that, under a natural scaling of the model parameters, the vehicle densities of $m$ different vehicle types in $d$ cells, can be accurately approximated (as a function of time, that is) by a $dm$-dimensional Gaussian process. Importantly, in particular it is shown that, by applying standard numerical software, both the intra-segment correlations and the temporal correlations can be efficiently evaluated in a straightforward manner. A second main result concerns the convergence, under the same scaling, of the joint per-cell cumulative arrival process of all vehicle types to a Gaussian limit, from which the approximations to the travel-time distributions are derived. The scaling mentioned above is to ensure there is enough aggregation, concretely meaning that the cells should be sufficiently long to guarantee that the number of vehicles being present is large enough for the central limit regime to kick in. Importantly, as has been extensively explored, central limit based approximations typically do not require an excessive scale in order to be reasonably accurate; it is anticipated that a couple of tens of vehicles per segment should suffice.

\vspace{2mm}

The primary goal of the present paper is to demonstrate the potential of the methodology developed in \cite{JL2012,JL2013,MS2019} in relation to various design and control related questions. To solidly back the use of that approach, we first perform an extensive analysis to empirically verify whether, in practice, Gaussian models can be used to accurately approximate the joint distribution of the vehicle density in the cells. This we do using detailed measurements performed in the Dutch highway network, so as to systematically assess the conditions under which multivariate Gaussian distributions can be safely assumed.

Having thoroughly justified the use of a Gaussian model, we proceed by showcasing a sequence of design and control issues that can be analyzed using the methodology developed in \cite{MS2019}. These are primarily meant as illustrations of the rich potential of this methodology, in that various other control and design applications can be dealt with. More concretely, in this paper we discuss the following examples in which our Gaussian framework is particularly useful. 
\begin{itemize}
\item[$\circ$] An integral part of the model in \cite{MS2019} is the fundamental diagram (see e.g., \cite[Section 2]{WageningenKessels2015genealogy} or \cite[Section V]{maerivoet2005traffic}) that provides us with the flow for any given value of the density. This fundamental diagram has a unimodal shape: evidently, the flow is low when the density is low, but also when the density is high (think of a traffic jam), with a peak between these two extremes. Interestingly, imposing traffic control measures such as a  speed limit effectively amounts to changing the fundamental diagram. By using the Gaussian-process approximation, we can quantify the impact of such measures. We can for instance evaluate by how much the travel time in a specific scenario goes up when reducing the maximum speed on a road segment. Alternatively, the reduction in the segment's throughput can be quantified. The resulting numbers provide policy makers with the information to judge whether these losses are outweighed by the positive effects (such as reduction of fuel consumption and/or air pollution). 
\item[$\circ$] Car navigation systems typically provide their users with a proxy of the estimated travel time to a particular destination. In practice, however, in the driver's decision often also the {\it uncertainty} of the travel time plays a role. How the mean and variance are weighed typically depends on the purpose of the trip (in particular the consequences of arriving too late) and the individual driver's preference; these weights can be captured by a utility function. For example, when facing the choice between a route with a high mean travel time but a low variance, and a route with a low mean and a high variance, the driver picks the route with the highest utility  \cite{SEN}; see also, {e.g.,} \cite{LOUI,SIVA}. As the framework in \cite{MS2019} is capable of evaluating an accurate approximation of the full travel-time distribution, in particular the mean and variance of the travel time can be evaluated, thus providing the driver with all information needed to properly select the best route.
\item[$\circ$] The Gaussian framework was primarily intended to describe the distribution of the numbers of vehicles at any point in time, but also the evaluation of stationary metrics as well. In this respect one can for instance think of the mean time, in stationary conditions, to traverse a given segment, thus shedding light on capacity-related performance measures. The employed technique is highly flexible, in that for any performance measure that is a function of the vehicle density or velocity, its long-term behavior can be evaluated. Other examples include  performance measures that quantify the deterioriation of roads (which can serve as the basis for maintenance decision making) and the carbon-dioxide emission rate (which can be used when developing environmental policies).
\item[$\circ$] When analyzing road traffic dynamics, the underlying network structure evidently plays a crucial role. In particular, typical network features like merging, diverging, backpressure and forward propagation of traffic, add to the complexity of evaluating the network performance. By using an example network, containing these network features, we demonstrate how the methodology in \cite{MS2019} is well-suited to handle such features. More precisely, we argue that for any network topology we can obtain accurate Gaussian approximations for the vehicle density process, jointly in all cells of all roads of the network.
\end{itemize}
A particularly attractive feature of the methodology is that it is capable of evaluating the joint effect of the implementation of multiple control measures simultaneously, in contrast to many existing methodologies that primarily focus on quantifying the impact of one measure in isolation, cf.\ \cite{papa}. As a second advantage, our methods are of an analytic nature, so that control-related questions can be dealt with relatively quickly. As such, for online applications and sensitivity analyses our approach is better suited than simulation-based techniques. In addition, as argued in, e.g., \cite{JL2013}, analytic methods also have a considerable advantage over simulation when computing covariance matrices of traffic state variables from traffic flow models with a fine time and space discretization.


The paper is organized as follows. In \autoref{sec:DAMV}, we present the analysis of the Dutch highway data set, assessing the validity of Gaussian approximations. In \autoref{sec: stat.perf.meas}, we summarize the traffic flow model that {was} studied in \cite{MS2019}, and present a brief account of {the} main findings. This section also covers our method to evaluate long-term performance of road traffic networks, using a methodology that incorporates the stochastic nature of traffic. \autoref{sec: route choice} illustrates the importance of working with  a probabilistic traffic model in the route selection context. We present  a route choice example where the chosen route depends on the {`disutility'}  due to the uncertainty inherent in the travel-time estimation. Then, in \autoref{sec: control}, we provide numerical examples that quantify how a road segment's performance is affected by the speed limit, the number of lanes, the arrival rate, and the fraction of vehicles in each of the vehicle classes. These experiments in particular demonstrate how the Gaussian approximation can be used for the purposes of control. Finally, in \autoref{sec: network}, we specifically consider  network structures in which traffic streams split and merge.
A short abstract covering some of the material presented in this paper has appeared as \cite{MS2020}.

\section{Data analysis and model validation}
\label{sec:DAMV}

As pointed out in the introduction, in \cite{MS2019} a multi-class stochastic traffic flow model has been developed. It splits the traffic network into segments (cells), and is capable of describing, in a stochastic manner, how traffic densities jointly propagate through the cells. A Gaussian-process approximation was established, using fluid and diffusion limits. It thus leads to an explicit approximation of the joint distribution of the vehicle densities (both in space and time) in terms of a multivariate Gaussian distribution.

The fluid and diffusion results in \cite{MS2019} have been established in a specific modeling framework, with a set of assumptions being imposed. Though these assumptions, which have their origin in conservation laws, traffic flow theory and stochastic analysis, are mild and commonly accepted, it still leaves open the question whether the Gaussian diffusion limit offers a sound approximation in practice. Two concerns play a role in this context. 
\begin{itemize}
\item[$\circ$]
In the first place, as described in the introduction, the diffusion limit is an {\it asymptotic} result. This means that in principle the approximations to the joint vehicle-density process, produced on the basis of the diffusion limit, might substantially deviate from the actual situation. In situations when the aggregation level of vehicles is sufficiently high, one expects that the central-limit theorem regime kicks in, so that the approximations are sound. One wonders, however, under what precise conditions Gaussian distributions offer a sufficiently accurate approximation. 
\item[$\circ$]
Secondly, the limiting results have been established under specific {\it assumptions}: the dynamics are supposed to be Markovian, and in addition specific functions (viz.\ transition rates, that reflect the fundamental diagram) are assumed to be Lipschitz continuous (i.e., continuous with a bounded derivative). It is important to realize, though, that even if these conditions are not fulfilled, the Gaussian framework potentially still provides an accurate approximation, albeit not backed by a formal limit result.
\end{itemize}

The main finding of this section is that we provide solid backing for the claim that it is reasonable to approximate the distributions of vehicles in a segment with a Gaussian distribution (with time-dependent parameters). We do so relying on a thorough analysis of a historical data set of vehicle flows on highways. We have performed an extensive numerical study, but in this section we restrict ourselves to presenting the main conclusions.
 
\subsection{Explanation of the data set}

We proceed by providing some background in the data set that we will use. This publicly available data is retrieved from the `Nationale Databank Wegverkeersgegevens'\footnote{\url{https://www.ndw.nu/pagina/nl/103/datalevering/120/open_data/}}. The data set contains measurements of the number of vehicles passing a measuring site, expressed in vehicles per hour, at 876 measuring locations in the Netherlands along the A2 highway, one of the main highways in the Netherlands. The measurements are recorded in one minute intervals, covering 167 days: from June 3rd 2018, 04 {\sc am}, till November 26th 2018, 11 {\sc pm}. 

In our analysis, we have restricted ourselves to time intervals between 4 {\sc am} and 11 {\sc am}. As a consequence,   the amount of data is somewhat reduced, but the time window chosen contains the quiet period at the beginning of the day as  well as the morning rush-hour. In Figure~\ref{Fig: Time series data}, for a single measuring site, the first 8 days of data as well as the first 24 days of data (both of them starting at a Monday) are shown for illustration.


\begin{figure}
\begin{picture}(450,260)(0,0)
\tikzsetnextfilename{Fig_timeseries_s}
\put(35,163){\tiny \rotatebox{-30}{(Mon)}}
\put(85.5,163){\tiny \rotatebox{-30}{(Tue)}}
\put(136,163){\tiny \rotatebox{-30}{(Wed)}}
\put(186.5,163){\tiny \rotatebox{-30}{(Thu)}}
\put(237,163){\tiny \rotatebox{-30}{(Fri)}}
\put(287.5,163){\tiny \rotatebox{-30}{(Sat)}}
\put(338,163){\tiny \rotatebox{-30}{(Sun)}}
\put(388.5,163){\tiny \rotatebox{-30}{(Mon)}}
\put(0,135){
\begin{tikzpicture}
\begin{axis}[
    xlabel={ Date },    ylabel={Flow (veh/h)},
    xmin = 2018-06-03 18:00, xmax = 2018-06-11 18:00,
    ymin=-10	, ymax=230,
    ytick={0,50,100,150,200},
    xtick={2018-06-04 00:00,2018-06-05 00:00,2018-06-06 00:00,2018-06-07 00:00,
    2018-06-08 00:00,2018-06-09 00:00,2018-06-10 00:00,2018-06-11 00:00},
    xtick pos=left,    ytick pos=left,    axis lines=left,
    x axis line style=-,    y axis line style=-,
    font = \tiny, 
    date ZERO=2018-06-03 18:00,
    width = 450pt, height=125pt,
    date coordinates in=x, 
    table/col sep=comma, 
    xticklabel=\year-\month-\day,
    xticklabel style={rotate=-30, anchor=near xticklabel}]
\addplot+[no markers, color35] table[x={time},y={flow}] {df_timeseries_s.csv};
\end{axis}
\end{tikzpicture}
}
\tikzsetnextfilename{Fig_timeseries_l}
\put(0,0){
\begin{tikzpicture}
\begin{axis}[
    xlabel={ Date },    ylabel={Flow (veh/h)},
    xmin = 2018-06-03 18:00, xmax = 2018-06-29 00:00,
    ymin=-10	, ymax=230,
    ytick={0,50,100,150,200},
    xtick={2018-06-04 00:00,2018-06-07 00:00,2018-06-10 00:00,2018-06-13 00:00,
    2018-06-16 00:00,2018-06-19 00:00,2018-06-22 00:00,2018-06-25 00:00,2018-06-28 00:00},
    xtick pos=left,    ytick pos=left,    axis lines=left,
    x axis line style=-,    y axis line style=-,
    font = \tiny, 
    date ZERO=2018-06-03 00:00,
    width = 450pt, height=125pt,
    date coordinates in=x, 
    table/col sep=comma, 
    xticklabel=\year-\month-\day,
    xticklabel style={rotate=-30, anchor=near xticklabel}]
\addplot+[no markers, color35] table[x={time},y={flow}] {df_timeseries_l.csv};
\end{axis}
\end{tikzpicture}
}
\end{picture}
\caption{\label{Fig: Time series data}Illustration of time series data: flow (over 1-minute intervals) as a function of time.}
\end{figure}
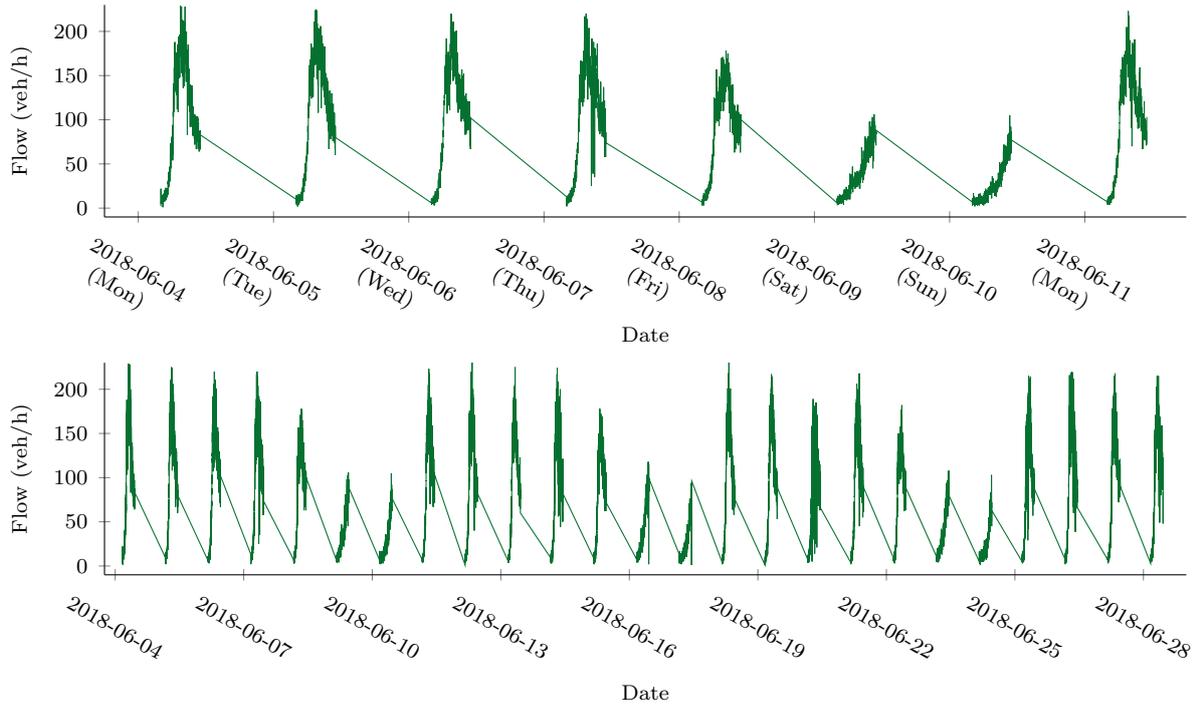

We wish to assess the validity of the claim that the number of vehicles on a road segment are approximately Gaussian. As a first step, we consider the flow over 1-minute intervals, measured at both endpoints of the segment. We collect the measurements of the 167 days, per day of the week and time~$t$ (in minutes) within the 4 {\sc am}--11 {\sc am} interval (with $t=1$ up to $7\times60=420$). 
As can be seen from Figure~\ref{Fig: Time series data}, the weekend days differ significantly from the working days, while among the working days the Fridays  slightly deviate  as well. 

If these flow measurements would stem from  a (multivariate) Gaussian distribution, then also each of the differences between the entries are (univariate) Gaussian. Focusing on adjacent (in time, that is) 1-minute intervals, this would imply that the change in density at time~$t$ is normally distributed. If we would find support for such a normal distribution for every 1-minute slot, then we have strong evidence for the desired {Gaussian distributions}. 

\subsection{Data preparation and preliminary (univariate) analysis}

As mentioned above, to validate the Gaussian approximation, we consider the empirical distribution of the flow at both endpoints of the site under consideration, for every weekday and every minute in the 4 {\sc am}--11 {\sc am} interval.

\subsubsection{Setup and methodology}

We already commented on the differences between the days of the week. Since we are interested in comparing distributions, the measurements under consideration, of which our sample consist, should ideally be drawn from identical distributions. As such, we will exclude Fridays, Saturdays and Sundays, so that we end up with data from 96 days.
We have not cleaned the data more than this, meaning that it may still contain  effects related to seasonality, holidays, etc. This means that if our analysis shows that normality is a reasonable approximation in this setting, then conceivably in reality this is even more the case.

In this section we focus on six adjacent measurement sites, on a segment without on- and off-ramps. These will be our central object during our study, that is, we show that the multivariate normal distribution is a good approximation to their empirical flow distributions; see Figure~\ref{Fig: location measurement sites} for the location of these measurement sites.\footnote{Image made with screenshots from: \url{https://dexter.ndwcloud.nu/opendata}, {which uses Mapbox and OpenStreetMap}} The names of the measurement sites are given in the Appendix~\ref{sec: appendix}, in Table~\ref{Table: Technical Info}. In light of the space available, we have decided to limit ourselves in this section to the results corresponding to three representative cases per experiment, two of which are  presented in Appendix~\ref{sec: appendix}.

\begin{figure}
\begin{tabular}{ccc}
\includegraphics[height=5cm]{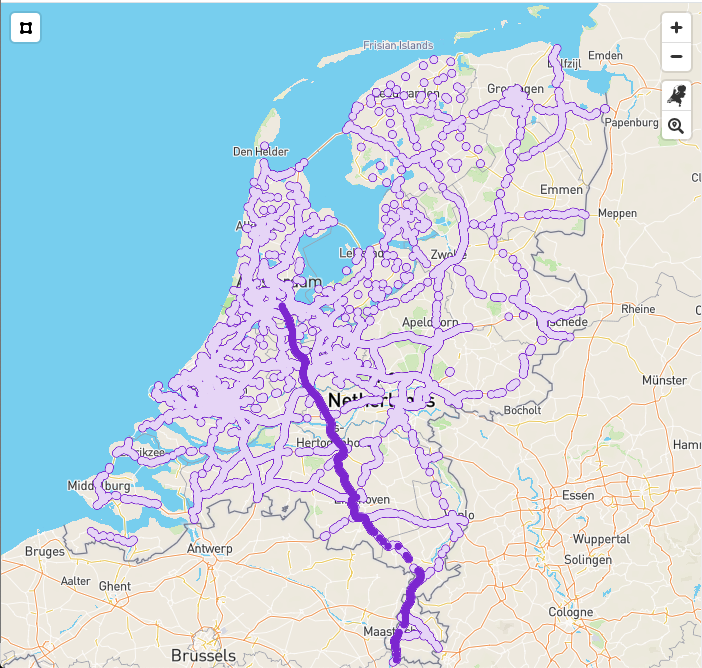} & 
\includegraphics[height=5cm]{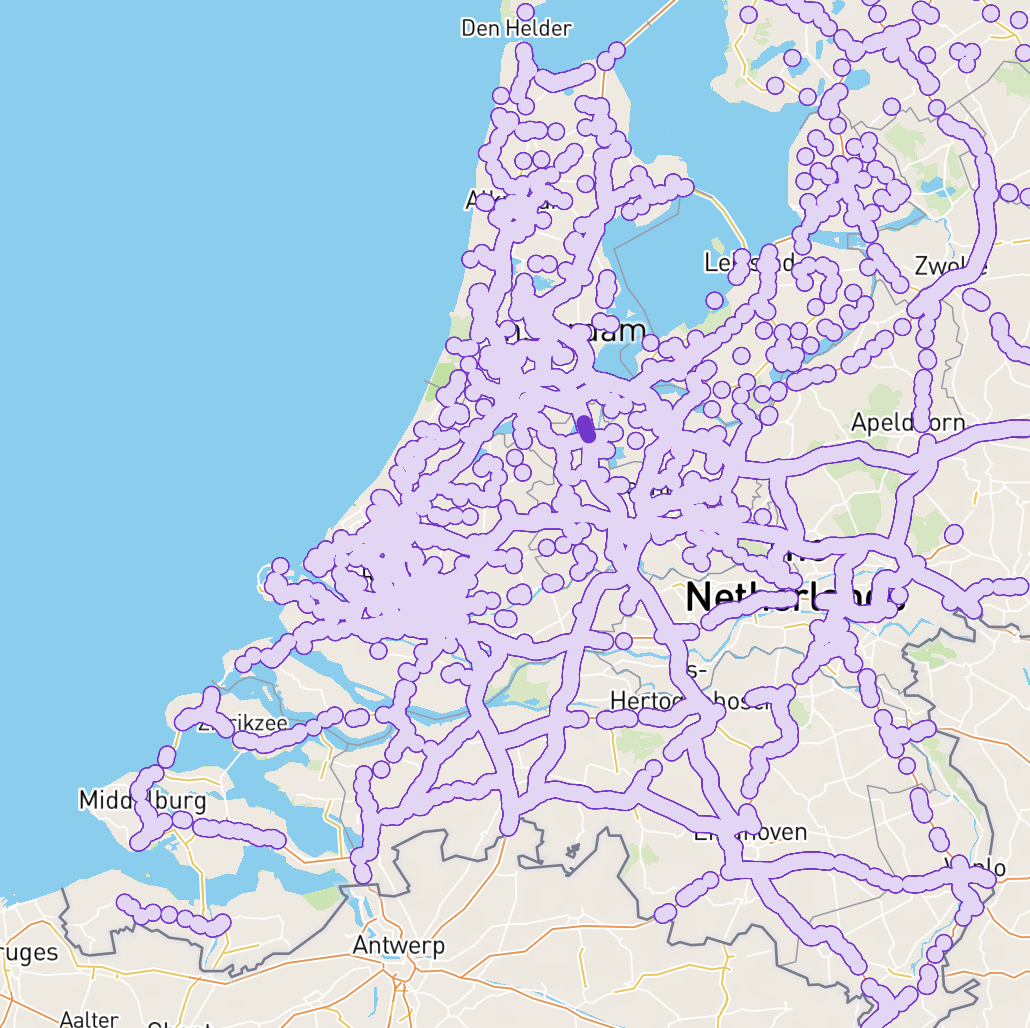} & 
\includegraphics[height=5cm]{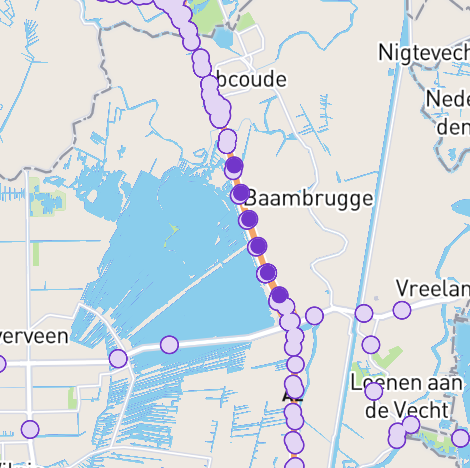}
\end{tabular}
\caption{\label{Fig: location measurement sites} Locations {(in dark purple)} of all measurement sites (left), location of {measurements under consideration (middle), and a zoomed-in image showing the analyzed measurement sites (right)}. The {pink} dots are other measurement sites in the Netherlands at which data is gathered.}
\end{figure}

\subsubsection{Results of the univariate analysis}

We now present our findings for the univariate empirical flow distributions. From Figure~\ref{Fig: Time series data}, we can see that there is a strong daily pattern, so  that when normality holds, the parameters will be time dependent. It is anticipated that the fit will become worse when we aggregate the per-minute flows into $\tau$-minute intervals, with $\tau$ an integer larger than 1. Our objective is to identify  the maximal value of $\tau$
up to which the univariate distributions are close to normal, giving a rough indication of the  maximal $\tau$ up to which we can expect a good fit for the multivariate case. 
Note that in addition in the univariate setting we can perform a graphical analysis with {\sc qq}-plots and histograms.  These graphical techniques  have their multivariate counterparts, but these are harder to visualize and interpret.

In this univariate analysis we consider the first three measurement sites. In Figure~\ref{Fig: distr flows 2} we have plotted the histograms (together with the associated best-fit normal density), and the corresponding {\sc qq}-plots, for the flow measurements at site 2, at times 05:20 {\sc am} and 07:20 {\sc am}, aggregated over $\tau$-minute intervals, with $\tau = 1,2,5,10,20$. More precisely, for a given value of~$\tau$ the sample corresponding to time 05:20 {\sc am} consists of all $\tau$ per-minute flow measurements taken in the interval between 05:20 {\sc am} and 05:($20 + \tau-1$) {\sc am}, for all days 96 days in the data set, which amounts to a sample of size  $96 \tau$ (approximately, as a consequence of missing observations).

Our first observation is that the histograms are reasonably bell-shaped, and are reasonably well approximated by the best-fit normal density, though we cannot draw strong conclusion from these figures. In addition, the {\sc qq}-plots are quite close to linear, especially for the 05:20 {\sc am} data. However, when $\tau$ becomes larger, we can see both from the histograms and (more explicitly) from the {\sc qq}-plots, that the tails of the empirical distribution become more skewed. This an be seen as a consequence of the non-stationarity, i.e., the distribution has a time-dependent mean and variance, which for larger $\tau$ results in skewed tails. In Figures~\ref{Fig: distr flows 1} and \ref{Fig: distr flows 3}, in Appendix~\ref{sec: appendix}, we have plotted the histograms and {\sc qq}-plots for sites 1 and 3, which show very similar behavior. In conclusion, the figures do not reject  normality, but for stronger backing we require additional analysis.

\begin{figure}
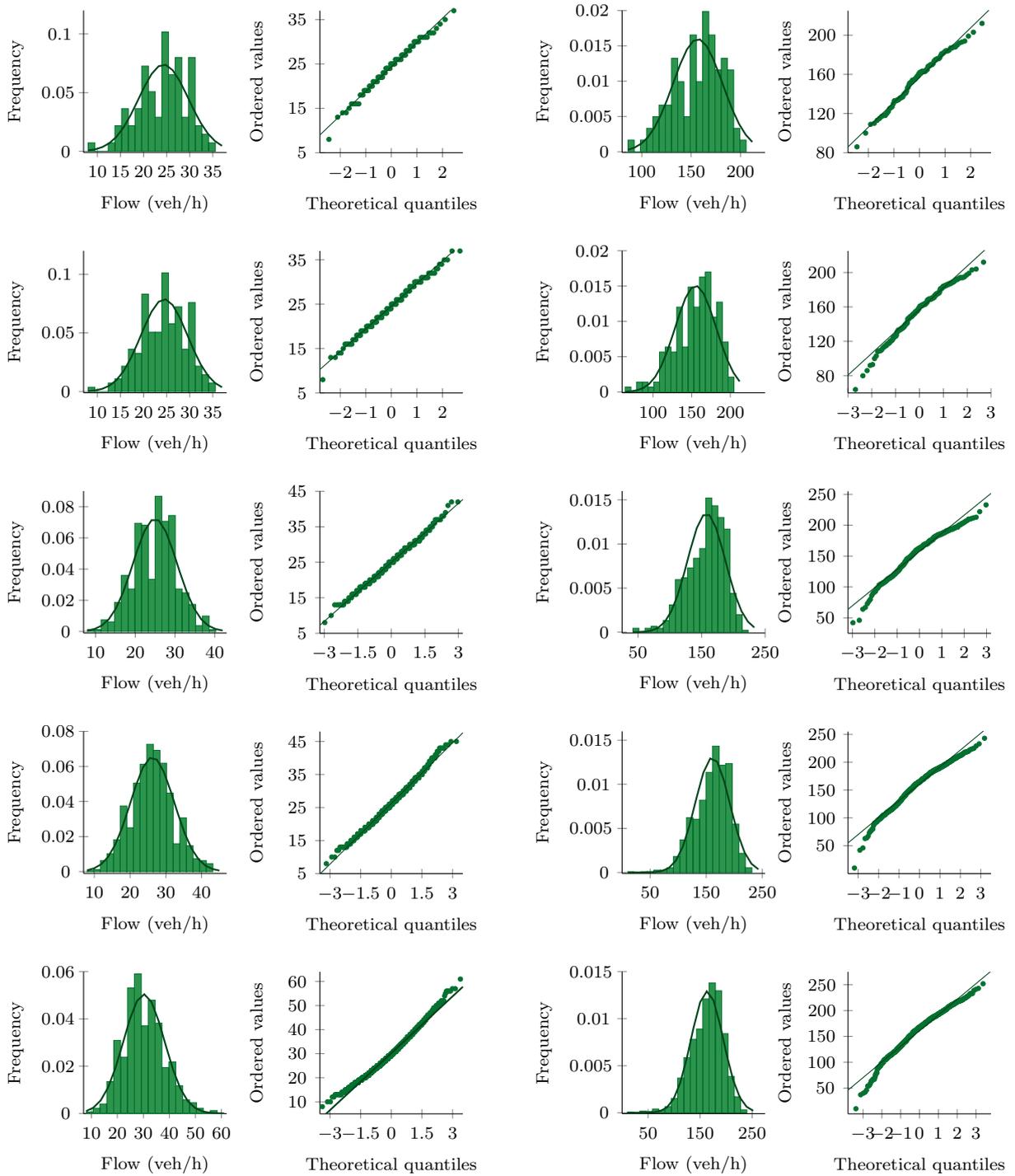


\caption{\label{Fig: distr flows 2}Distribution of flows at site 2, at times 05:20 {\sc am} (left) and 07:20 {\sc am}, measured in intervals of $\tau$ minutes, for $\tau = 1,2,5,10,20$ (read from above to below).}
\end{figure}

To find additional support, we rely on statistical goodness-of-fit tests. More precisely, for every $\tau$-minute interval between 4 {\sc am} and 11 {\sc am}., we apply a $\chi^2$ goodness-of-fit test  to the corresponding  empirical distribution, with the null-hypothesis being that the data is distributed as a normal distribution, with mean and variance as given by their respective maximum likelihood estimators. Note that these parameters are estimated for every $\tau$-minute interval, and are therefore time-varying. The $\chi^2$ test is applied to samples of size at least  $90$, with 10 bins that under the null-hypothesis contain 10\% probability mass each. Hence, the number of expected observations per bin is at least 9 for each test, so that the test is sufficiently reliable. 

We have plotted the resulting p-values cumulatively, as a function of time for $\tau = 1,2,5,10$, in Figure~\ref{Fig: p-value integrals 2}. For $\tau = 20$ the cumulative p-value was 0.25, and thus we reject the null-hypothesis for a normal distribution for almost every $\tau$-minute interval. In the figure, every line segment has a color that is associated to its slope, in the sense that the value of the p-value on the right endpoint of a segment determines the color. One can thus, from the color of the curve at that time, visually distinguish the `quality of a normal fit' at a specific moment in time. We observe that, since the p-values correspond to quantiles, that under the null-hypothesis the slope of the curve should be $0.5$. 

\begin{figure}
\begin{picture}(350 , 200)(0,0)
\put(0,100){
\tikzsetnextfilename{Fig_pvalues_univ_site=1_min=1}
\begin{tikzpicture}	
\begin{axis}[
    xlabel={time (h)},    ylabel={Cumulative p-value},
    xmin=3.5, xmax=11.5,    ymin=-5, ymax=210,
    xtick={4,6,8,10},    ytick={0,50,100,150,200},
    xticklabel style={/pgf/number format/fixed, /pgf/number format/precision=2},
    yticklabel style={/pgf/number format/fixed, /pgf/number format/precision=3,text width=1.3em,align=right},
	scaled x ticks=false,
    xtick pos=left,    ytick pos=left,    axis lines=left,
    width = 150pt,    height = 100pt,
    no markers,    x axis line style=-,    y axis line style=-,
    font = \tiny,
    legend style={at={(1,.8)},anchor=north east, draw=none},
    colorbar,
    colormap={my colormap}{
            color=(red)
            color=(yellow)
            color=(color33)
			color=(blue)
        },
    colormap access=piecewise constant,
    colorbar style={at={(1.1,1.0)},    anchor=north west,    height=55.5pt,     font = \tiny,
    width=0.2cm,    ymin=1,    ymax=4,    ytick={1,1.75,2.5,3.25,4},    yticklabels={0,0.01,0.05,0.5,1},
    }
]
\addplot+[very thick, mesh, point meta=explicit] 
table [x={time}, y={segment=1_pchi2_cum}, point meta=\thisrow{segment=1_pchi2}] {df_cumPvalues_univ_min=1.dat};
\end{axis}
\end{tikzpicture}
}
\put(190,100){
\tikzsetnextfilename{Fig_pvalues_univ_site=1_min=2}
\begin{tikzpicture}	
\begin{axis}[
    xlabel={time (h)},    ylabel={Cumulative p-value},
    xmin=3.5, xmax=11.5,    ymin=-2, ymax=105,
    xtick={4,6,8,10},    ytick={0,25,50,75,100},
    xticklabel style={/pgf/number format/fixed, /pgf/number format/precision=2},
    yticklabel style={/pgf/number format/fixed, /pgf/number format/precision=3,text width=1.3em,align=right},
	scaled x ticks=false,
    xtick pos=left,    ytick pos=left,    axis lines=left,
    width = 150pt,    height = 100pt,
    no markers,    x axis line style=-,    y axis line style=-,
    font = \tiny,
    legend style={at={(1,.8)},anchor=north east, draw=none},
    colorbar,
    colormap={my colormap}{
            color=(red)
            color=(yellow)
            color=(color33)
			color=(blue)
        },
    colormap access=piecewise constant,
    colorbar style={at={(1.1,1.0)},    anchor=north west,    height=55.5pt,     font = \tiny,
    width=0.2cm,    ymin=1,    ymax=4,    ytick={1,1.75,2.5,3.25,4},    yticklabels={0,0.01,0.05,0.5,1},
    }
]
\addplot+[very thick, mesh, point meta=explicit] 
table [x={time}, y={segment=1_pchi2_cum}, point meta=\thisrow{segment=1_pchi2}] {df_cumPvalues_univ_min=2.dat};
\end{axis}
\end{tikzpicture}
}
\put(0,0){
\tikzsetnextfilename{Fig_pvalues_univ_site=1_min=5}
\begin{tikzpicture}	
\begin{axis}[
    xlabel={time (h)},    ylabel={Cumulative p-value},
    xmin=3.5, xmax=11.5,    ymin=-1, ymax=42,
    xtick={4,6,8,10},    ytick={0,10,20,30,40,50},
    xticklabel style={/pgf/number format/fixed, /pgf/number format/precision=2},
    yticklabel style={/pgf/number format/fixed, /pgf/number format/precision=3,text width=1.3em,align=right},
	scaled x ticks=false,
    xtick pos=left,    ytick pos=left,    axis lines=left,
    width = 150pt,    height = 100pt,
    no markers,    x axis line style=-,    y axis line style=-,
    font = \tiny,
    legend style={at={(1,.8)},anchor=north east, draw=none},
    colorbar,
    colormap={my colormap}{
            color=(red)
            color=(yellow)
            color=(color33)
			color=(blue)
        },
    colormap access=piecewise constant,
    colorbar style={at={(1.1,1.0)},    anchor=north west,    height=55.5pt,     font = \tiny,
    width=0.2cm,    ymin=1,    ymax=4,    ytick={1,1.75,2.5,3.25,4},    yticklabels={0,0.01,0.05,0.5,1},
    }
]
\addplot+[very thick, mesh, point meta=explicit] 
table [x={time}, y={segment=1_pchi2_cum}, point meta=\thisrow{segment=1_pchi2}] {df_cumPvalues_univ_min=5.dat};
\end{axis}
\end{tikzpicture}
}
\put(190,0){
\tikzsetnextfilename{Fig_pvalues_univ_site=1_min=10}
\begin{tikzpicture}	
\begin{axis}[
    xlabel={time (h)},    ylabel={Cumulative p-value},
    xmin=3.5, xmax=11.5,    ymin=-0.5, ymax=21,
    xtick={4,6,8,10},    ytick={0,5,10,15,20},
    xticklabel style={/pgf/number format/fixed, /pgf/number format/precision=2},
    yticklabel style={/pgf/number format/fixed, /pgf/number format/precision=3,text width=1.3em,align=right},
	scaled x ticks=false,
    xtick pos=left,    ytick pos=left,    axis lines=left,
    width = 150pt,    height = 100pt,
    no markers,    x axis line style=-,    y axis line style=-,
    font = \tiny,
    legend style={at={(1,.8)},anchor=north east, draw=none},
    colorbar,
    colormap={my colormap}{
            color=(red)
            color=(yellow)
            color=(color33)
			color=(blue)
        },
    colormap access=piecewise constant,
    colorbar style={at={(1.1,1.0)},    anchor=north west,    height=55.5pt,     font = \tiny,
    width=0.2cm,    ymin=1,    ymax=3,    ytick={1,1.5,2,2.5,3},    yticklabels={0,0.01,0.05,0.5,1},
    }
]
\addplot+[very thick, mesh, point meta=explicit] 
table [x={time}, y={segment=1_pchi2_cum}, point meta=\thisrow{segment=1_pchi2}] {df_cumPvalues_univ_min=10.dat};
\end{axis}
\end{tikzpicture}
}
\end{picture}
\caption{\label{Fig: p-value integrals 2}
Cumulative p-values of $\chi^2$ tests at site 2, as a function of time, for $\tau$-minute intervals, with $\tau=1$ (top left), $\tau =2$ (top right), $\tau=5$ (bottom left), and $\tau = 10$ (bottom right).}
\end{figure}
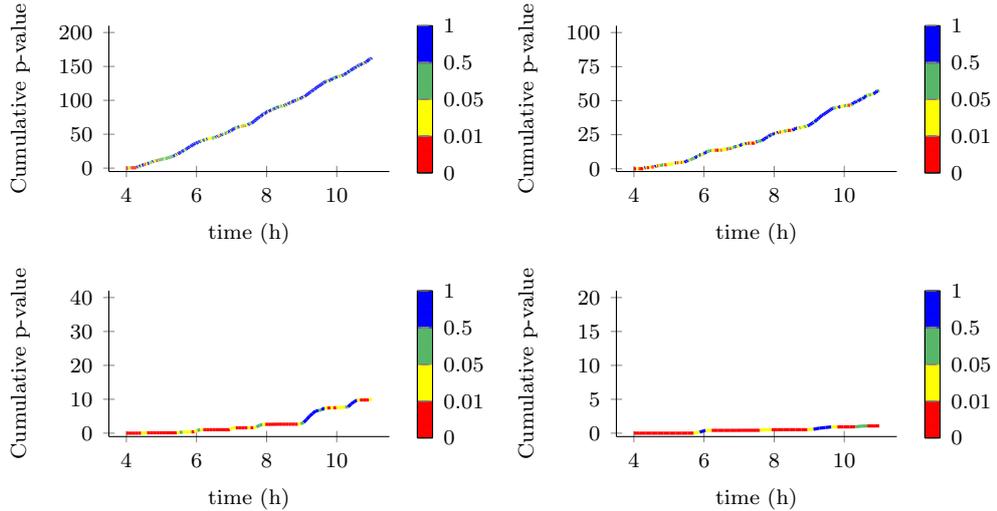
%

The main conclusion from Figure~\ref{Fig: p-value integrals 2} is that the fit becomes worse when $\tau$ increases. This is again an indication of the non-stationarity discussed above. More precisely, for $\tau \geqslant 5$, the $\chi^2$ test rejects the null-hypothesis (of a normal distribution, that is) at more than half the $\tau$-minute intervals, with the fit becoming worse for $\tau = 10$. In contrast, the fit for $\tau=2$ is already quite reasonable, with normality only rejected at a small number of intervals, mostly at around 4 {\sc am}. For $\tau = 1$, the curve is close to having the desired slope of $0.5$, which is strong indication that distribution of the $1$-minute intervals is indeed Gaussian. The analogous figures for sites 1 and 3 can be found in Appendix~\ref{sec: appendix}, in Figures~\ref{Fig: p-value integrals 1} and \ref{Fig: p-value integrals 3}, and display similar results.

Based on the above, when considering the multivariate case, we only have to consider $\tau$ equal to $1$ or $2$: at that timescale {Gaussian distributions} offer an accurate description for the one-dimensional flows.
We proceed in the next section by extending our analysis to the multivariate case.

\subsection{Analysis of the multivariate case}
We consider the joint empirical distribution at two measurement sites, so as to assess whether joint Gaussian distributions provide an accurate approximation. We first explain our methodology, after which we present our results.

\subsubsection{Setup and methodology}

The study of distributions in the multivariate case is considerably more involved than in the univariate case. To begin with, though in principle graphical analysis in the scenario of two-dimensional distributions is possible, it is quite a bit harder to interpret two-dimensional histograms and {\sc qq}-plots. Therefore, we have decided to follow another approach, and focus on hypothesis testing instead. Due to the lack of a goodness-of-fit test for multivariate normality, we utilize the result that a stochastic vector $X = (X_1,\ldots,X_p)^\top$ of length $p$, $p \in \N$, is $p$-variate normally distributed if and only if $\theta^\top X$ is univariate normally distributed, for all $\theta \in \R^p$, where a point mass is seen as a normal distribution with variance zero \cite[p.\ 383]{billingsley}. Our approach is now to apply a goodness-of-fit test to linear combinations of flow measurements, measured at identical times, at two different measurement sites. 

{In our setup, we consider two measurement sites from the segment consisting of six adjacent sites, which we index by natural numbers $i,j \in \{1,\ldots,6\}$.}
For both measurement sites, we have observations of the per-minute flow, which we interpret of realizations of the two random variables $X_i$ and $X_j$. To verify that $(X_i,X_j)$ has a bivariate normal distribution, we consider linear combinations $\alpha X_i + \beta X_j$, for a finite set of pairs $(\alpha,\beta)$. Take $\alpha,\beta\in\{-2,-1,-1/2,1/2,1,2\}$, by evident scaling properties it suffices to consider only the pairs 
\[
	\{(2,-2),(2,-1),(2,-1/2),(2,1/2),(2,1),(2,2),(-1,2),(-1/2,2),(1/2,2),(1,2)\}.
\] 
For each pair, we apply the $\chi^2$ goodness-of-fit test to $\alpha X_i + \beta X_j$, with the null-hypothesis that the distribution is normal, for each $\tau$-minute interval of observations, $\tau \in \{1,2\}$. The setup of the $\chi^2$ test is the same as in the univariate case, working with 10 bins with 10\% probability mass each under the null-hypothesis.

%
%

\subsubsection{Results of the multivariate analysis}

\begin{figure}
\begin{picture}(435 , 360)(0,0)
\put(0,270){
\tikzsetnextfilename{Fig_pvalues_multi_site=1_min=1_a=-5_b=20}
\begin{tikzpicture}	
\begin{axis}[
    xlabel={time $t$ (h)},    ylabel={$P_{(-0.5,2)}(t)$},
    xmin=3.5, xmax=11.5,    ymin=-5, ymax=210,
    xtick={4,6,8,10},    ytick={0,50,100,150,200},
    xticklabel style={/pgf/number format/fixed, /pgf/number format/precision=2},
    yticklabel style={/pgf/number format/fixed, /pgf/number format/precision=3,text width=1.3em,align=right},
	scaled x ticks=false,
    xtick pos=left,    ytick pos=left,    axis lines=left,
    width = 125pt,    height = 90pt,
    no markers,    x axis line style=-,    y axis line style=-,
    y label style={at={(axis description cs:-0.22,.5)}, anchor=south},
    font = \tiny,
    legend style={at={(1,.8)},anchor=north east, draw=none},
    colorbar,
    colormap={my colormap}{
            color=(red)
            color=(yellow)
            color=(color33)
			color=(blue)
        },
    colormap access=piecewise constant,
    colorbar style={at={(1.1,1.0)},    anchor=north west,    height=45pt,     font = \tiny,
    width=0.12cm,    ymin=1,    ymax=4,    ytick={1,1.75,2.5,3.25,4},    yticklabels={0,0.01,0.05,0.5,1},
    }
]
\addplot+[very thick, mesh, point meta=explicit] 
table [x={domain}, y={pValues_cumulative_Alpha=-5_Beta=20}, point meta=\thisrow{pValues_Alpha=-5_Beta=20}] {df_pvalues_multiv_segment=1_min=1.dat};
\end{axis}
\end{tikzpicture}
}
\put(155,270){
\tikzsetnextfilename{Fig_pvalues_multi_site=1_min=1_a=-10_b=20}
\begin{tikzpicture}	
\begin{axis}[
    xlabel={time $t$ (h)},    ylabel={$P_{(-1,2)}(t)$},
    xmin=3.5, xmax=11.5,    ymin=-5, ymax=210,
    xtick={4,6,8,10},    ytick={0,50,100,150,200},
    xticklabel style={/pgf/number format/fixed, /pgf/number format/precision=2},
    yticklabel style={/pgf/number format/fixed, /pgf/number format/precision=3,text width=1.3em,align=right},
	scaled x ticks=false,
    xtick pos=left,    ytick pos=left,    axis lines=left,
    width = 125pt,    height = 90pt,
    no markers,    x axis line style=-,    y axis line style=-,
    y label style={at={(axis description cs:-0.22,.5)}, anchor=south},
    font = \tiny,
    legend style={at={(1,.8)},anchor=north east, draw=none},
    colorbar,
    colormap={my colormap}{
            color=(red)
            color=(yellow)
            color=(color33)
			color=(blue)
        },
    colormap access=piecewise constant,
    colorbar style={at={(1.1,1.0)},    anchor=north west,    height=45pt,     font = \tiny,
    width=0.12cm,    ymin=1,    ymax=4,    ytick={1,1.75,2.5,3.25,4},    yticklabels={0,0.01,0.05,0.5,1},
    }
]
\addplot+[very thick, mesh, point meta=explicit] 
table [x={domain}, y={pValues_cumulative_Alpha=-10_Beta=20}, point meta=\thisrow{pValues_Alpha=-10_Beta=20}] {df_pvalues_multiv_segment=1_min=1.dat};
\end{axis}
\end{tikzpicture}
}
\put(310,270){
\tikzsetnextfilename{Fig_pvalues_multi_site=1_min=1_a=5_b=20}
\begin{tikzpicture}	
\begin{axis}[
    xlabel={time $t$ (h)},    ylabel={$P_{(0.5,2)}(t)$},
    xmin=3.5, xmax=11.5,    ymin=-5, ymax=210,
    xtick={4,6,8,10},    ytick={0,50,100,150,200},
    xticklabel style={/pgf/number format/fixed, /pgf/number format/precision=2},
    yticklabel style={/pgf/number format/fixed, /pgf/number format/precision=3,text width=1.3em,align=right},
	scaled x ticks=false,
    xtick pos=left,    ytick pos=left,    axis lines=left,
    width = 125pt,    height = 90pt,
    no markers,    x axis line style=-,    y axis line style=-,
    y label style={at={(axis description cs:-0.22,.5)}, anchor=south},
    font = \tiny,
    legend style={at={(1,.8)},anchor=north east, draw=none},
    colorbar,
    colormap={my colormap}{
            color=(red)
            color=(yellow)
            color=(color33)
			color=(blue)
        },
    colormap access=piecewise constant,
    colorbar style={at={(1.1,1.0)},    anchor=north west,    height=45pt,     font = \tiny,
    width=0.12cm,    ymin=1,    ymax=4,    ytick={1,1.75,2.5,3.25,4},    yticklabels={0,0.01,0.05,0.5,1},
    }
]
\addplot+[very thick, mesh, point meta=explicit] 
table [x={domain}, y={pValues_cumulative_Alpha=5_Beta=20}, point meta=\thisrow{pValues_Alpha=5_Beta=20}] {df_pvalues_multiv_segment=1_min=1.dat};
\end{axis}
\end{tikzpicture}
}
\put(0,180){
\tikzsetnextfilename{Fig_pvalues_multi_site=1_min=1_a=10_b=20}
\begin{tikzpicture}	
\begin{axis}[
    xlabel={time $t$ (h)},    ylabel={$P_{(1,2)}(t)$},
    xmin=3.5, xmax=11.5,    ymin=-5, ymax=210,
    xtick={4,6,8,10},    ytick={0,50,100,150,200},
    xticklabel style={/pgf/number format/fixed, /pgf/number format/precision=2},
    yticklabel style={/pgf/number format/fixed, /pgf/number format/precision=3,text width=1.3em,align=right},
	scaled x ticks=false,
    xtick pos=left,    ytick pos=left,    axis lines=left,
    width = 125pt,    height = 90pt,
    no markers,    x axis line style=-,    y axis line style=-,
    y label style={at={(axis description cs:-0.22,.5)}, anchor=south},
    font = \tiny,
    legend style={at={(1,.8)},anchor=north east, draw=none},
    colorbar,
    colormap={my colormap}{
            color=(red)
            color=(yellow)
            color=(color33)
			color=(blue)
        },
    colormap access=piecewise constant,
    colorbar style={at={(1.1,1.0)},    anchor=north west,    height=45pt,     font = \tiny,
    width=0.12cm,    ymin=1,    ymax=4,    ytick={1,1.75,2.5,3.25,4},    yticklabels={0,0.01,0.05,0.5,1},
    }
]
\addplot+[very thick, mesh, point meta=explicit] 
table [x={domain}, y={pValues_cumulative_Alpha=10_Beta=20}, point meta=\thisrow{pValues_Alpha=10_Beta=20}] {df_pvalues_multiv_segment=1_min=1.dat};
\end{axis}
\end{tikzpicture}
}
\put(155,180){
\tikzsetnextfilename{Fig_pvalues_multi_site=1_min=1_a=20_b=-20}
\begin{tikzpicture}	
\begin{axis}[
    xlabel={time $t$ (h)},    ylabel={$P_{(2,-2)}(t)$},
    xmin=3.5, xmax=11.5,    ymin=-5, ymax=210,
    xtick={4,6,8,10},    ytick={0,50,100,150,200},
    xticklabel style={/pgf/number format/fixed, /pgf/number format/precision=2},
    yticklabel style={/pgf/number format/fixed, /pgf/number format/precision=3,text width=1.3em,align=right},
	scaled x ticks=false,
    xtick pos=left,    ytick pos=left,    axis lines=left,
    width = 125pt,    height = 90pt,
    no markers,    x axis line style=-,    y axis line style=-,
    y label style={at={(axis description cs:-0.22,.5)}, anchor=south},
    font = \tiny,
    legend style={at={(1,.8)},anchor=north east, draw=none},
    colorbar,
    colormap={my colormap}{
            color=(red)
            color=(yellow)
            color=(color33)
			color=(blue)
        },
    colormap access=piecewise constant,
    colorbar style={at={(1.1,1.0)},    anchor=north west,    height=45pt,     font = \tiny,
    width=0.12cm,    ymin=1,    ymax=4,    ytick={1,1.75,2.5,3.25,4},    yticklabels={0,0.01,0.05,0.5,1},
    }
]
\addplot+[very thick, mesh, point meta=explicit] 
table [x={domain}, y={pValues_cumulative_Alpha=20_Beta=-20}, point meta=\thisrow{pValues_Alpha=20_Beta=-20}] {df_pvalues_multiv_segment=1_min=1.dat};
\end{axis}
\end{tikzpicture}
}
\put(310,180){
\tikzsetnextfilename{Fig_pvalues_multi_site=1_min=1_a=20_b=-10}
\begin{tikzpicture}	
\begin{axis}[
    xlabel={time $t$ (h)},    ylabel={$P_{(2,-1)}(t)$},
    xmin=3.5, xmax=11.5,    ymin=-5, ymax=210,
    xtick={4,6,8,10},    ytick={0,50,100,150,200},
    xticklabel style={/pgf/number format/fixed, /pgf/number format/precision=2},
    yticklabel style={/pgf/number format/fixed, /pgf/number format/precision=3,text width=1.3em,align=right},
	scaled x ticks=false,
    xtick pos=left,    ytick pos=left,    axis lines=left,
    width = 125pt,    height = 90pt,
    no markers,    x axis line style=-,    y axis line style=-,
    y label style={at={(axis description cs:-0.22,.5)}, anchor=south},
    font = \tiny,
    legend style={at={(1,.8)},anchor=north east, draw=none},
    colorbar,
    colormap={my colormap}{
            color=(red)
            color=(yellow)
            color=(color33)
			color=(blue)
        },
    colormap access=piecewise constant,
    colorbar style={at={(1.1,1.0)},    anchor=north west,    height=45pt,     font = \tiny,
    width=0.12cm,    ymin=1,    ymax=4,    ytick={1,1.75,2.5,3.25,4},    yticklabels={0,0.01,0.05,0.5,1},
    }
]
\addplot+[very thick, mesh, point meta=explicit] 
table [x={domain}, y={pValues_cumulative_Alpha=20_Beta=-10}, point meta=\thisrow{pValues_Alpha=20_Beta=-10}] {df_pvalues_multiv_segment=1_min=1.dat};
\end{axis}
\end{tikzpicture}
}
\put(0,90){
\tikzsetnextfilename{Fig_pvalues_multi_site=1_min=1_a=20_b=-5}
\begin{tikzpicture}	
\begin{axis}[
    xlabel={time $t$ (h)},    ylabel={$P_{(2,-0.5)}(t)$},
    xmin=3.5, xmax=11.5,    ymin=-5, ymax=210,
    xtick={4,6,8,10},    ytick={0,50,100,150,200},
    xticklabel style={/pgf/number format/fixed, /pgf/number format/precision=2},
    yticklabel style={/pgf/number format/fixed, /pgf/number format/precision=3,text width=1.3em,align=right},
	scaled x ticks=false,
    xtick pos=left,    ytick pos=left,    axis lines=left,
    width = 125pt,    height = 90pt,
    no markers,    x axis line style=-,    y axis line style=-,
    y label style={at={(axis description cs:-0.22,.5)}, anchor=south},
    font = \tiny,
    legend style={at={(1,.8)},anchor=north east, draw=none},
    colorbar,
    colormap={my colormap}{
            color=(red)
            color=(yellow)
            color=(color33)
			color=(blue)
        },
    colormap access=piecewise constant,
    colorbar style={at={(1.1,1.0)},    anchor=north west,    height=45pt,     font = \tiny,
    width=0.12cm,    ymin=1,    ymax=4,    ytick={1,1.75,2.5,3.25,4},    yticklabels={0,0.01,0.05,0.5,1},
    }
]
\addplot+[very thick, mesh, point meta=explicit] 
table [x={domain}, y={pValues_cumulative_Alpha=20_Beta=-5}, point meta=\thisrow{pValues_Alpha=20_Beta=-5}] {df_pvalues_multiv_segment=1_min=1.dat};
\end{axis}
\end{tikzpicture}
}
\put(155,90){
\tikzsetnextfilename{Fig_pvalues_multi_site=1_min=1_a=20_b=5}
\begin{tikzpicture}	
\begin{axis}[
    xlabel={time $t$ (h)},    ylabel={$P_{(2,0.5)}(t)$},
    xmin=3.5, xmax=11.5,    ymin=-5, ymax=210,
    xtick={4,6,8,10},    ytick={0,50,100,150,200},
    xticklabel style={/pgf/number format/fixed, /pgf/number format/precision=2},
    yticklabel style={/pgf/number format/fixed, /pgf/number format/precision=3,text width=1.3em,align=right},
	scaled x ticks=false,
    xtick pos=left,    ytick pos=left,    axis lines=left,
    width = 125pt,    height = 90pt,
    no markers,    x axis line style=-,    y axis line style=-,
    y label style={at={(axis description cs:-0.22,.5)}, anchor=south},
    font = \tiny,
    legend style={at={(1,.8)},anchor=north east, draw=none},
    colorbar,
    colormap={my colormap}{
            color=(red)
            color=(yellow)
            color=(color33)
			color=(blue)
        },
    colormap access=piecewise constant,
    colorbar style={at={(1.1,1.0)},    anchor=north west,    height=45pt,     font = \tiny,
    width=0.12cm,    ymin=1,    ymax=4,    ytick={1,1.75,2.5,3.25,4},    yticklabels={0,0.01,0.05,0.5,1},
    }
]
\addplot+[very thick, mesh, point meta=explicit] 
table [x={domain}, y={pValues_cumulative_Alpha=20_Beta=5}, point meta=\thisrow{pValues_Alpha=20_Beta=5}] {df_pvalues_multiv_segment=1_min=1.dat};
\end{axis}
\end{tikzpicture}
}
\put(310,90){
\tikzsetnextfilename{Fig_pvalues_multi_site=1_min=1_a=20_b=10}
\begin{tikzpicture}	
\begin{axis}[
    xlabel={time $t$ (h)},    ylabel={$P_{(2,1)}(t)$},
    xmin=3.5, xmax=11.5,    ymin=-5, ymax=210,
    xtick={4,6,8,10},    ytick={0,50,100,150,200},
    xticklabel style={/pgf/number format/fixed, /pgf/number format/precision=2},
    yticklabel style={/pgf/number format/fixed, /pgf/number format/precision=3,text width=1.3em,align=right},
	scaled x ticks=false,
    xtick pos=left,    ytick pos=left,    axis lines=left,
    width = 125pt,    height = 90pt,
    no markers,    x axis line style=-,    y axis line style=-,
    y label style={at={(axis description cs:-0.22,.5)}, anchor=south},
    font = \tiny,
    legend style={at={(1,.8)},anchor=north east, draw=none},
    colorbar,
    colormap={my colormap}{
            color=(red)
            color=(yellow)
            color=(color33)
			color=(blue)
        },
    colormap access=piecewise constant,
    colorbar style={at={(1.1,1.0)},    anchor=north west,    height=45pt,     font = \tiny,
    width=0.12cm,    ymin=1,    ymax=4,    ytick={1,1.75,2.5,3.25,4},    yticklabels={0,0.01,0.05,0.5,1},
    }
]
\addplot+[very thick, mesh, point meta=explicit] 
table [x={domain}, y={pValues_cumulative_Alpha=20_Beta=10}, point meta=\thisrow{pValues_Alpha=20_Beta=10}] {df_pvalues_multiv_segment=1_min=1.dat};
\end{axis}
\end{tikzpicture}
}
\put(155,0){
\tikzsetnextfilename{Fig_pvalues_multi_site=1_min=1_a=20_b=20}
\begin{tikzpicture}	
\begin{axis}[
    xlabel={time $t$ (h)},    ylabel={$P_{(2,2)}(t)$},
    xmin=3.5, xmax=11.5,    ymin=-5, ymax=210,
    xtick={4,6,8,10},    ytick={0,50,100,150,200},
    xticklabel style={/pgf/number format/fixed, /pgf/number format/precision=2},
    yticklabel style={/pgf/number format/fixed, /pgf/number format/precision=3,text width=1.3em,align=right},
	scaled x ticks=false,
    xtick pos=left,    ytick pos=left,    axis lines=left,
    width = 125pt,    height = 90pt,
    no markers,    x axis line style=-,    y axis line style=-,
    y label style={at={(axis description cs:-0.22,.5)}, anchor=south},
    font = \tiny,
    legend style={at={(1,.8)},anchor=north east, draw=none},
    colorbar,
    colormap={my colormap}{
            color=(red)
            color=(yellow)
            color=(color33)
			color=(blue)
        },
    colormap access=piecewise constant,
    colorbar style={at={(1.1,1.0)},    anchor=north west,    height=45pt,     font = \tiny,
    width=0.12cm,    ymin=1,    ymax=4,    ytick={1,1.75,2.5,3.25,4},    yticklabels={0,0.01,0.05,0.5,1},
    }
]
\addplot+[very thick, mesh, point meta=explicit] 
table [x={domain}, y={pValues_cumulative_Alpha=20_Beta=20}, point meta=\thisrow{pValues_Alpha=20_Beta=20}] {df_pvalues_multiv_segment=1_min=1.dat};
\end{axis}
\end{tikzpicture}
}
\end{picture}
\caption{Cumulative p-values of $\chi^2$-test ($P_{(\alpha,\beta)}(t)$), for linear combinations of measurements at site 2 and site 3, as a function of time $t$, for $\tau = 1$.}
\label{Fig: p-value integrals joints, sites 2,3, t = 1}
\end{figure}
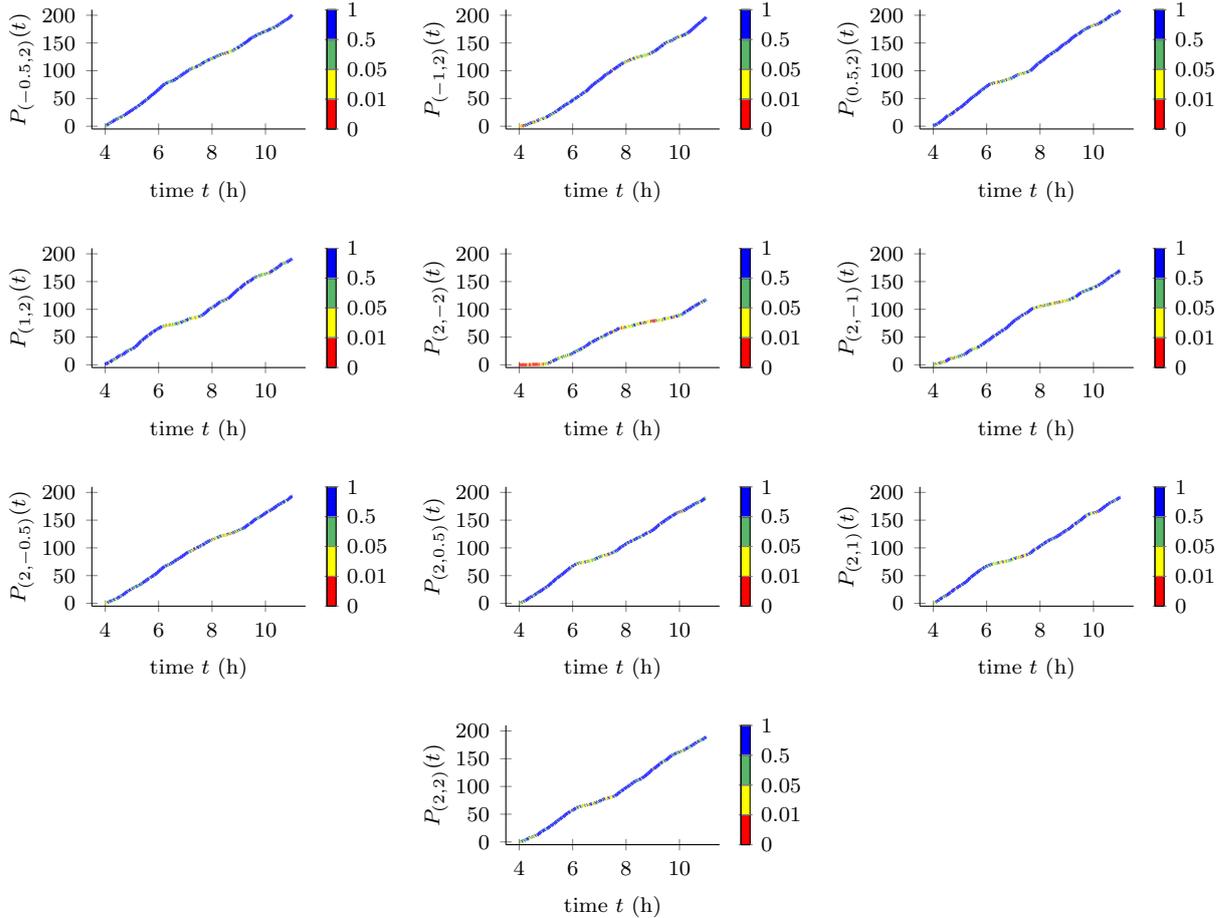

We continue by presenting the results of the $\chi^2$-test. In Figures~\ref{Fig: p-value integrals joints, sites 2,3, t = 1} and \ref{Fig: p-value integrals joints, sites 2,3, t = 2}, we have plotted the cumulative p-value as a function of time, following the procedure underlying Figure~\ref{Fig: p-value integrals 2}, with, respectively, $\tau=1$ and $\tau=2$. The vertical axes correspond to the interval $[0,210 / \tau]$, since $210 / \tau$ is the expected height of the curve under the null-hypothesis, making sure that the curve has slope $0.5$.

For $\tau = 1$, Figures~\ref{Fig: p-value integrals joints, sites 2,3, t = 1} shows that the $\chi^2$-test does not reject normality for almost every time $t$ between 4 {\sc am} and 11 {\sc am}, with the slope of the curves being between $0.4$ and $0.5$ for most figures. The only exception is $(\alpha,\beta)= (2, -2)$; here it may play a role that when subtracting two random variables of the same order of magnitude may lead to a quantity with a relatively high coefficient of variation (realizing that the mean of the difference is small). Comparing these figures to the corresponding plots for the univariate case,  in Figure~\ref{Fig: p-value integrals 2}, we see that the fit becomes slightly worse, but not significantly; in some cases the fit becomes even a bit better. For $\tau = 2$, we see the similar results, with the overall conclusion that the fit is still more than acceptable for a significant part of the times $t$ between 4 {\sc am} and 11 {\sc am}. {Similar conclusions can be drawn for the joint distribution of flows at measurement sites~1 and 2 (see Figures~\ref{Fig: p-value integrals joints, sites 1 and 2, t = 1} and \ref{Fig: p-value integrals joints, sites 1 and 2, t = 2}, in Appendix~\ref{sec: appendix}). Regarding the role of dependence, increasing the distance between the measurement sites does not clearly improve the results (see Figures~\ref{Fig: p-value integrals joints, sites 1 and 6, t = 1} and \ref{Fig: p-value integrals joints, sites 1 and 6, t = 2} for the results corresponding to site~1 and site~6), in Appendix~\ref{sec: appendix}}.

Based on these above, we conclude that bivariate Gaussian distributions are accurate approximations to the bivariate empirical flow distributions in the data set. This provides further support for the joint cumulative flow process (and hence also the joint vehicle densities process)
being accurately approximated by a Gaussian process.


\begin{figure}
\begin{picture}(435 , 360)(0,0)
\put(0,270){
\tikzsetnextfilename{Fig_pvalues_multi_site=1_min=2_a=-5_b=20}
\begin{tikzpicture}	
\begin{axis}[
    xlabel={time $t$ (h)},    ylabel={$P_{(-0.5,2)}(t)$},
    xmin=3.5, xmax=11.5,    ymin=-2.5, ymax=105,
    xtick={4,6,8,10},    ytick={0,25,50,75,100},
    xticklabel style={/pgf/number format/fixed, /pgf/number format/precision=2},
    yticklabel style={/pgf/number format/fixed, /pgf/number format/precision=3,text width=1.3em,align=right},
	scaled x ticks=false,
    xtick pos=left,    ytick pos=left,    axis lines=left,
    width = 125pt,    height = 90pt,
    no markers,    x axis line style=-,    y axis line style=-,
    y label style={at={(axis description cs:-0.22,.5)}, anchor=south},
    font = \tiny,
    legend style={at={(1,.8)},anchor=north east, draw=none},
    colorbar,
    colormap={my colormap}{
            color=(red)
            color=(yellow)
            color=(color33)
			color=(blue)
        },
    colormap access=piecewise constant,
    colorbar style={at={(1.1,1.0)},    anchor=north west,    height=45pt,     font = \tiny,
    width=0.12cm,    ymin=1,    ymax=4,    ytick={1,1.75,2.5,3.25,4},    yticklabels={0,0.01,0.05,0.5,1},
    }
]
\addplot+[very thick, mesh, point meta=explicit] 
table [x={domain}, y={pValues_cumulative_Alpha=-5_Beta=20}, point meta=\thisrow{pValues_Alpha=-5_Beta=20}] {df_pvalues_multiv_segment=1_min=2.dat};
\end{axis}
\end{tikzpicture}
}
\put(155,270){
\tikzsetnextfilename{Fig_pvalues_multi_site=1_min=2_a=-10_b=20}
\begin{tikzpicture}	
\begin{axis}[
    xlabel={time $t$ (h)},    ylabel={$P_{(-1,2)}(t)$},
    xmin=3.5, xmax=11.5,    ymin=-2.5, ymax=105,
    xtick={4,6,8,10},    ytick={0,25,50,75,100},
    xticklabel style={/pgf/number format/fixed, /pgf/number format/precision=2},
    yticklabel style={/pgf/number format/fixed, /pgf/number format/precision=3,text width=1.3em,align=right},
	scaled x ticks=false,
    xtick pos=left,    ytick pos=left,    axis lines=left,
    width = 125pt,    height = 90pt,
    no markers,    x axis line style=-,    y axis line style=-,
    y label style={at={(axis description cs:-0.22,.5)}, anchor=south},
    font = \tiny,
    legend style={at={(1,.8)},anchor=north east, draw=none},
    colorbar,
    colormap={my colormap}{
            color=(red)
            color=(yellow)
            color=(color33)
			color=(blue)
        },
    colormap access=piecewise constant,
    colorbar style={at={(1.1,1.0)},    anchor=north west,    height=45pt,     font = \tiny,
    width=0.12cm,    ymin=1,    ymax=4,    ytick={1,1.75,2.5,3.25,4},    yticklabels={0,0.01,0.05,0.5,1},
    }
]
\addplot+[very thick, mesh, point meta=explicit] 
table [x={domain}, y={pValues_cumulative_Alpha=-10_Beta=20}, point meta=\thisrow{pValues_Alpha=-10_Beta=20}] {df_pvalues_multiv_segment=1_min=2.dat};
\end{axis}
\end{tikzpicture}
}
\put(310,270){
\tikzsetnextfilename{Fig_pvalues_multi_site=1_min=2_a=5_b=20}
\begin{tikzpicture}	
\begin{axis}[
    xlabel={time $t$ (h)},    ylabel={$P_{(0.5,2)}(t)$},
    xmin=3.5, xmax=11.5,    ymin=-2.5, ymax=105,
    xtick={4,6,8,10},    ytick={0,25,50,75,100},
    xticklabel style={/pgf/number format/fixed, /pgf/number format/precision=2},
    yticklabel style={/pgf/number format/fixed, /pgf/number format/precision=3,text width=1.3em,align=right},
	scaled x ticks=false,
    xtick pos=left,    ytick pos=left,    axis lines=left,
    width = 125pt,    height = 90pt,
    no markers,    x axis line style=-,    y axis line style=-,
    y label style={at={(axis description cs:-0.22,.5)}, anchor=south},
    font = \tiny,
    legend style={at={(1,.8)},anchor=north east, draw=none},
    colorbar,
    colormap={my colormap}{
            color=(red)
            color=(yellow)
            color=(color33)
			color=(blue)
        },
    colormap access=piecewise constant,
    colorbar style={at={(1.1,1.0)},    anchor=north west,    height=45pt,     font = \tiny,
    width=0.12cm,    ymin=1,    ymax=4,    ytick={1,1.75,2.5,3.25,4},    yticklabels={0,0.01,0.05,0.5,1},
    }
]
\addplot+[very thick, mesh, point meta=explicit] 
table [x={domain}, y={pValues_cumulative_Alpha=5_Beta=20}, point meta=\thisrow{pValues_Alpha=5_Beta=20}] {df_pvalues_multiv_segment=1_min=2.dat};
\end{axis}
\end{tikzpicture}
}
\put(0,180){
\tikzsetnextfilename{Fig_pvalues_multi_site=1_min=2_a=10_b=20}
\begin{tikzpicture}	
\begin{axis}[
    xlabel={time $t$ (h)},    ylabel={$P_{(1,2)}(t)$},
    xmin=3.5, xmax=11.5,    ymin=-2.5, ymax=105,
    xtick={4,6,8,10},    ytick={0,25,50,75,100},
    xticklabel style={/pgf/number format/fixed, /pgf/number format/precision=2},
    yticklabel style={/pgf/number format/fixed, /pgf/number format/precision=3,text width=1.3em,align=right},
	scaled x ticks=false,
    xtick pos=left,    ytick pos=left,    axis lines=left,
    width = 125pt,    height = 90pt,
    no markers,    x axis line style=-,    y axis line style=-,
    y label style={at={(axis description cs:-0.22,.5)}, anchor=south},
    font = \tiny,
    legend style={at={(1,.8)},anchor=north east, draw=none},
    colorbar,
    colormap={my colormap}{
            color=(red)
            color=(yellow)
            color=(color33)
			color=(blue)
        },
    colormap access=piecewise constant,
    colorbar style={at={(1.1,1.0)},    anchor=north west,    height=45pt,     font = \tiny,
    width=0.12cm,    ymin=1,    ymax=4,    ytick={1,1.75,2.5,3.25,4},    yticklabels={0,0.01,0.05,0.5,1},
    }
]
\addplot+[very thick, mesh, point meta=explicit] 
table [x={domain}, y={pValues_cumulative_Alpha=10_Beta=20}, point meta=\thisrow{pValues_Alpha=10_Beta=20}] {df_pvalues_multiv_segment=1_min=2.dat};
\end{axis}
\end{tikzpicture}
}
\put(155,180){
\tikzsetnextfilename{Fig_pvalues_multi_site=1_min=2_a=20_b=-20}
\begin{tikzpicture}	
\begin{axis}[
    xlabel={time $t$ (h)},    ylabel={$P_{(2,-2)}(t)$},
    xmin=3.5, xmax=11.5,    ymin=-2.5, ymax=105,
    xtick={4,6,8,10},    ytick={0,25,50,75,100},
    xticklabel style={/pgf/number format/fixed, /pgf/number format/precision=2},
    yticklabel style={/pgf/number format/fixed, /pgf/number format/precision=3,text width=1.3em,align=right},
	scaled x ticks=false,
    xtick pos=left,    ytick pos=left,    axis lines=left,
    width = 125pt,    height = 90pt,
    no markers,    x axis line style=-,    y axis line style=-,
    y label style={at={(axis description cs:-0.22,.5)}, anchor=south},
    font = \tiny,
    legend style={at={(1,.8)},anchor=north east, draw=none},
    colorbar,
    colormap={my colormap}{
            color=(red)
            color=(yellow)
            color=(color33)
			color=(blue)
        },
    colormap access=piecewise constant,
    colorbar style={at={(1.1,1.0)},    anchor=north west,    height=45pt,     font = \tiny,
    width=0.12cm,    ymin=1,    ymax=4,    ytick={1,1.75,2.5,3.25,4},    yticklabels={0,0.01,0.05,0.5,1},
    }
]
\addplot+[very thick, mesh, point meta=explicit] 
table [x={domain}, y={pValues_cumulative_Alpha=20_Beta=-20}, point meta=\thisrow{pValues_Alpha=20_Beta=-20}] {df_pvalues_multiv_segment=1_min=2.dat};
\end{axis}
\end{tikzpicture}
}
\put(310,180){
\tikzsetnextfilename{Fig_pvalues_multi_site=1_min=2_a=20_b=-10}
\begin{tikzpicture}	
\begin{axis}[
    xlabel={time $t$ (h)},    ylabel={$P_{(2,-1)}(t)$},
    xmin=3.5, xmax=11.5,    ymin=-2.5, ymax=105,
    xtick={4,6,8,10},    ytick={0,25,50,75,100},
    xticklabel style={/pgf/number format/fixed, /pgf/number format/precision=2},
    yticklabel style={/pgf/number format/fixed, /pgf/number format/precision=3,text width=1.3em,align=right},
	scaled x ticks=false,
    xtick pos=left,    ytick pos=left,    axis lines=left,
    width = 125pt,    height = 90pt,
    no markers,    x axis line style=-,    y axis line style=-,
    y label style={at={(axis description cs:-0.22,.5)}, anchor=south},
    font = \tiny,
    legend style={at={(1,.8)},anchor=north east, draw=none},
    colorbar,
    colormap={my colormap}{
            color=(red)
            color=(yellow)
            color=(color33)
			color=(blue)
        },
    colormap access=piecewise constant,
    colorbar style={at={(1.1,1.0)},    anchor=north west,    height=45pt,     font = \tiny,
    width=0.12cm,    ymin=1,    ymax=4,    ytick={1,1.75,2.5,3.25,4},    yticklabels={0,0.01,0.05,0.5,1},
    }
]
\addplot+[very thick, mesh, point meta=explicit] 
table [x={domain}, y={pValues_cumulative_Alpha=20_Beta=-10}, point meta=\thisrow{pValues_Alpha=20_Beta=-10}] {df_pvalues_multiv_segment=1_min=2.dat};
\end{axis}
\end{tikzpicture}
}
\put(0,90){
\tikzsetnextfilename{Fig_pvalues_multi_site=1_min=2_a=20_b=-5}
\begin{tikzpicture}	
\begin{axis}[
    xlabel={time $t$ (h)},    ylabel={$P_{(2,-0.5)}(t)$},
    xmin=3.5, xmax=11.5,    ymin=-2.5, ymax=105,
    xtick={4,6,8,10},    ytick={0,25,50,75,100},
    xticklabel style={/pgf/number format/fixed, /pgf/number format/precision=2},
    yticklabel style={/pgf/number format/fixed, /pgf/number format/precision=3,text width=1.3em,align=right},
	scaled x ticks=false,
    xtick pos=left,    ytick pos=left,    axis lines=left,
    width = 125pt,    height = 90pt,
    no markers,    x axis line style=-,    y axis line style=-,
    y label style={at={(axis description cs:-0.22,.5)}, anchor=south},
    font = \tiny,
    legend style={at={(1,.8)},anchor=north east, draw=none},
    colorbar,
    colormap={my colormap}{
            color=(red)
            color=(yellow)
            color=(color33)
			color=(blue)
        },
    colormap access=piecewise constant,
    colorbar style={at={(1.1,1.0)},    anchor=north west,    height=45pt,     font = \tiny,
    width=0.12cm,    ymin=1,    ymax=4,    ytick={1,1.75,2.5,3.25,4},    yticklabels={0,0.01,0.05,0.5,1},
    }
]
\addplot+[very thick, mesh, point meta=explicit] 
table [x={domain}, y={pValues_cumulative_Alpha=20_Beta=-5}, point meta=\thisrow{pValues_Alpha=20_Beta=-5}] {df_pvalues_multiv_segment=1_min=2.dat};
\end{axis}
\end{tikzpicture}
}
\put(155,90){
\tikzsetnextfilename{Fig_pvalues_multi_site=1_min=2_a=20_b=5}
\begin{tikzpicture}	
\begin{axis}[
    xlabel={time $t$ (h)},    ylabel={$P_{(2,0.5)}(t)$},
    xmin=3.5, xmax=11.5,    ymin=-2.5, ymax=105,
    xtick={4,6,8,10},    ytick={0,25,50,75,100},
    xticklabel style={/pgf/number format/fixed, /pgf/number format/precision=2},
    yticklabel style={/pgf/number format/fixed, /pgf/number format/precision=3,text width=1.3em,align=right},
	scaled x ticks=false,
    xtick pos=left,    ytick pos=left,    axis lines=left,
    width = 125pt,    height = 90pt,
    no markers,    x axis line style=-,    y axis line style=-,
    y label style={at={(axis description cs:-0.22,.5)}, anchor=south},
    font = \tiny,
    legend style={at={(1,.8)},anchor=north east, draw=none},
    colorbar,
    colormap={my colormap}{
            color=(red)
            color=(yellow)
            color=(color33)
			color=(blue)
        },
    colormap access=piecewise constant,
    colorbar style={at={(1.1,1.0)},    anchor=north west,    height=45pt,     font = \tiny,
    width=0.12cm,    ymin=1,    ymax=4,    ytick={1,1.75,2.5,3.25,4},    yticklabels={0,0.01,0.05,0.5,1},
    }
]
\addplot+[very thick, mesh, point meta=explicit] 
table [x={domain}, y={pValues_cumulative_Alpha=20_Beta=5}, point meta=\thisrow{pValues_Alpha=20_Beta=5}] {df_pvalues_multiv_segment=1_min=2.dat};
\end{axis}
\end{tikzpicture}
}
\put(310,90){
\tikzsetnextfilename{Fig_pvalues_multi_site=1_min=2_a=20_b=10}
\begin{tikzpicture}	
\begin{axis}[
    xlabel={time $t$ (h)},    ylabel={$P_{(2,1)}(t)$},
    xmin=3.5, xmax=11.5,    ymin=-2.5, ymax=105,
    xtick={4,6,8,10},    ytick={0,25,50,75,100},
    xticklabel style={/pgf/number format/fixed, /pgf/number format/precision=2},
    yticklabel style={/pgf/number format/fixed, /pgf/number format/precision=3,text width=1.3em,align=right},
	scaled x ticks=false,
    xtick pos=left,    ytick pos=left,    axis lines=left,
    width = 125pt,    height = 90pt,
    no markers,    x axis line style=-,    y axis line style=-,
    y label style={at={(axis description cs:-0.22,.5)}, anchor=south},
    font = \tiny,
    legend style={at={(1,.8)},anchor=north east, draw=none},
    colorbar,
    colormap={my colormap}{
            color=(red)
            color=(yellow)
            color=(color33)
			color=(blue)
        },
    colormap access=piecewise constant,
    colorbar style={at={(1.1,1.0)},    anchor=north west,    height=45pt,     font = \tiny,
    width=0.12cm,    ymin=1,    ymax=4,    ytick={1,1.75,2.5,3.25,4},    yticklabels={0,0.01,0.05,0.5,1},
    }
]
\addplot+[very thick, mesh, point meta=explicit] 
table [x={domain}, y={pValues_cumulative_Alpha=20_Beta=10}, point meta=\thisrow{pValues_Alpha=20_Beta=10}] {df_pvalues_multiv_segment=1_min=2.dat};
\end{axis}
\end{tikzpicture}
}
\put(155,0){
\tikzsetnextfilename{Fig_pvalues_multi_site=1_min=2_a=20_b=20}
\begin{tikzpicture}	
\begin{axis}[
    xlabel={time $t$ (h)},    ylabel={$P_{(2,2)}(t)$},
    xmin=3.5, xmax=11.5,    ymin=-2.5, ymax=105,
    xtick={4,6,8,10},    ytick={0,25,50,75,100},
    xticklabel style={/pgf/number format/fixed, /pgf/number format/precision=2},
    yticklabel style={/pgf/number format/fixed, /pgf/number format/precision=3,text width=1.3em,align=right},
	scaled x ticks=false,
    xtick pos=left,    ytick pos=left,    axis lines=left,
    width = 125pt,    height = 90pt,
    no markers,    x axis line style=-,    y axis line style=-,
    y label style={at={(axis description cs:-0.22,.5)}, anchor=south},
    font = \tiny,
    legend style={at={(1,.8)},anchor=north east, draw=none},
    colorbar,
    colormap={my colormap}{
            color=(red)
            color=(yellow)
            color=(color33)
			color=(blue)
        },
    colormap access=piecewise constant,
    colorbar style={at={(1.1,1.0)},    anchor=north west,    height=45pt,     font = \tiny,
    width=0.12cm,    ymin=1,    ymax=4,    ytick={1,1.75,2.5,3.25,4},    yticklabels={0,0.01,0.05,0.5,1},
    }
]
\addplot+[very thick, mesh, point meta=explicit] 
table [x={domain}, y={pValues_cumulative_Alpha=20_Beta=20}, point meta=\thisrow{pValues_Alpha=20_Beta=20}] {df_pvalues_multiv_segment=1_min=2.dat};
\end{axis}
\end{tikzpicture}
}
\end{picture}
\caption{Cumulative p-values of $\chi^2$-test ($P_{(\alpha,\beta)}(t)$), for linear combinations of measurements at site 2 and site 3, as a function of time $t$, for $\tau = 2$.}
\label{Fig: p-value integrals joints, sites 2,3, t = 2}
\end{figure}
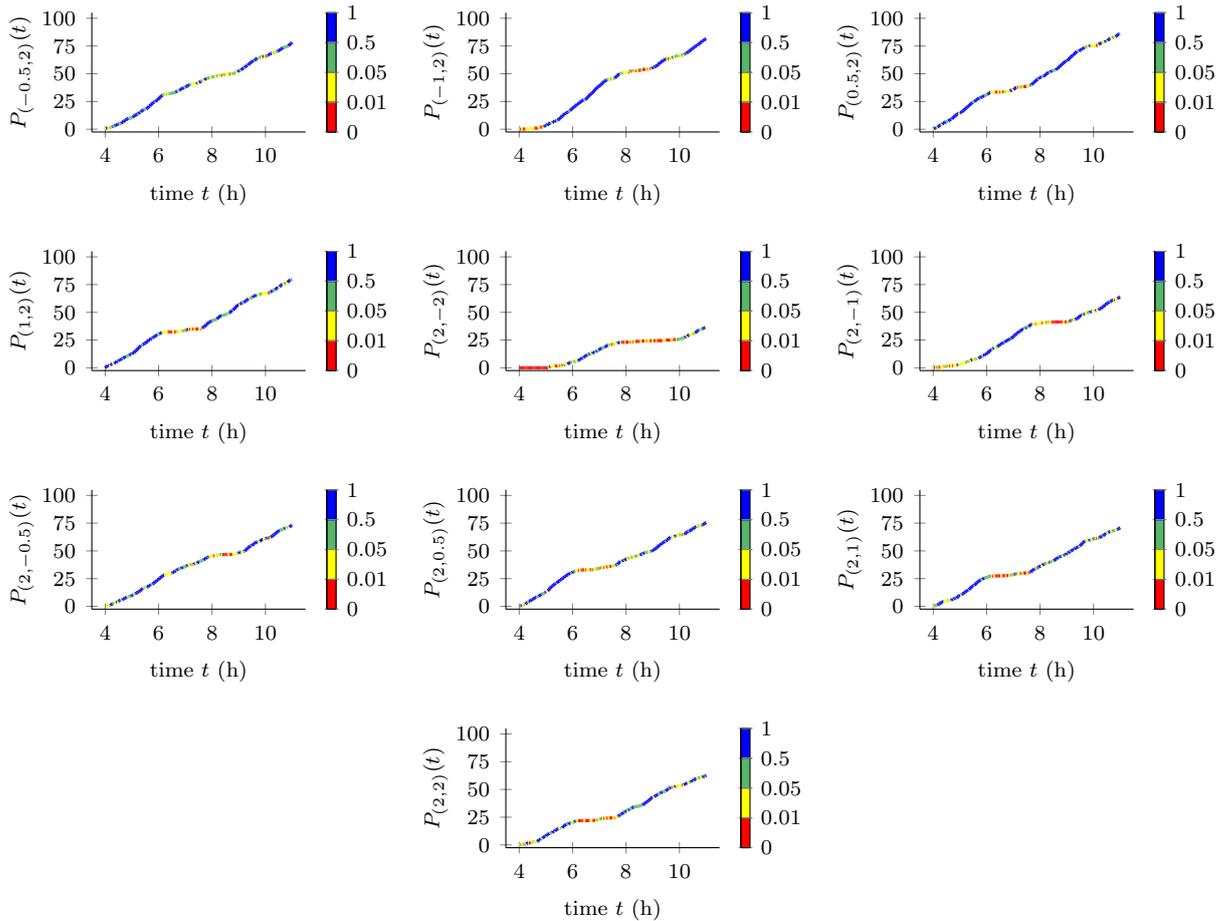


\section{Evaluation of stationary performance measures}
\label{sec: stat.perf.meas}

In this section we show the model of \cite{MS2019} can be used to evaluate long-term performance measures, that explicitly incorporate the stochastic nature of traffic. Examples of such performance metrics are the number of vehicles that traverse a road per unit of time, or the amount of carbon-dioxide emitted per unit of time. The results can be used to support decision making, as we shall illustrate in  \autoref{sec: control}. 

In this section we first briefly describe the model and main results from \cite{MS2019}. We restrict ourselves to aspects that are relevant for road traffic researchers, in that for all rigorous backing we refer to  \cite{MS2019}. 
The main idea is that we express long-term performance measures in terms of the stationary distribution of an underlying Markov chain.  By approximating this stationary distribution, using the Gaussian diffusion limit, we obtain an accurate and efficient method for computing the performance measure of interest.


\subsection{Model summary}
\label{sec:MS}

We now provide a compact model description; see \cite[Section 3]{MS2019} for more details. For ease we focus on a segment consisting of a sequence of cells, modeling a road segment without any intermediate on-ramps and off-ramps, rather than a general network. As pointed out in \cite[Section 7.1]{MS2019}, however, the setup naturally extends to networks;  see also \autoref{sec: network}.

The segment is divided in $d$ cells, with cell~$i$ having length $1/\ell_i > 0$, for $i \in \{1,\ldots,d\}$. We consider $m$ different vehicle types. We denote by $X_{ij}(t)$ the number of type-$j$ vehicles, in cell~$i$ at time $t \geqslant 0$, for $j \in \{1,\ldots,m\}$ and $i \in \{1,\ldots,d\}$. We assume $X_{ij}(t) \in \{0,1,\ldots,X^{\textrm{jam}}_{ij}\}$, where $X^{\textrm{jam}}_{ij}$ is the maximum number of type-$j$ vehicles that can be simultaneously present in cell~$i$. We denote $X(t)$ for the $md$-dimensional random vector with entries $X_{ij}(t)$. 
We define the type-$j$ vehicle density in cell~$i$ at time $t \geqslant 0$ by
\[
	\rho_{ij}(t) \colon = \frac{X_{ij}(t)}{\ell_i},
\]
attaining values in $\{0,1/\ell_i,\ldots,X^{\textrm{jam}}_{ij}/\ell_i\}$, where $\rho^{\textrm{jam}}_{ij}:=X^{\textrm{jam}}_{ij}/\ell_i$ can be interpreted as the jamming vehicle density. 
Analogously to $X(t)$, we let $\rho(t)$ be the $md$-dimensional random vector with entries $\rho_{ij}(t)$.

Vehicles traverse the consecutive cells, starting at cell~$1$ and ending at cell~$d$. In a \textit{cellular transition model} (CTM), vehicle-mass moves across cell boundaries at a rate that is given by a \textit{discrete flux-function}. This function, whose arguments are the vehicle densities in the sending and receiving cells (of all $m$ types), is derived from a macroscopic fundamental diagram (MFD) by solving an associated Riemann problem \cite[Section 3.2]{MS2019}. More formally, this means that
the discrete flux-function from cell $i$ to cell $i+1$ is given by a function of the type
\[
	\tilde q_i : \bigtimes_{j=1}^m \left( [0,\rho^{{\rm jam}}_{ij}] \times [0,\rho_{i+1,j}^{{\rm jam}}] \right) \to \bigtimes_{j=1}^m  [0,q^{{\rm max}}_{ij}],
\] 
where, as pointed out in \cite[Assumption 3.1]{MS2019}, a mild regularity condition has to be imposed on these $\tilde q_i$.

In our setup, we consider stochastic inter-cell transition times, so that the {\it mean}
dynamics correspond with a CTM. To be precise, the time it takes a type-$j$ vehicle to
move from  cell $i$ to cell $i + 1$ is an exponentially distributed random variable, with 
a mean that is in line the discrete flux-function. In addition to vehicles jumping between cells, vehicles enter the segment at cell~1 and depart from the segment at cell~$d$. These transitions are to be handled slightly different from the inter-cell transitions. Concretely, the arrival rate at cell 1 (say of type $j$) depends on the vehicle densities in cell 1 but is in addition bounded by a (given) rate $\lambda_j$. Likewise, the type-$j$ departure rate at cell $d$ is a function of the vehicle densities in cell $d$, truncated at $\nu_j$. {Since in our framework the transition times are exponentially distributed, the process under study is  a continuous-time Markov chain. }

We denote by $q_{0,j}(\rho(t))$, $q_{i,j}(\rho(t))$ and $q_{d,j}(\rho(t))$, respectively, the type-$j$ arrival rate at cell~$1$, transition rate between cell~$i$ and cell~$i+1$, and departure rate from cell~$d$, given the state of the system $\rho(t)$, where $j \in \{1,\ldots,m\}$ and $i \in \{1,\ldots,d-1\}$. We observe that these rates characterize the infinitesimal generator underlying our Markov process. Importantly, we can write the type-$j$ density in cell~$i$ by
\[
	\rho_{ij}(t) = \rho_{ij}(0) + \frac{1}{\ell_i} \int_0^t Y_{i-1,j}\left( \int_0^t q_{i-1,j}\left(\rho(s)\right) \diff s \right) - \frac{1}{\ell_i} \int_0^t Y_{i,j}\left( \int_0^t q_{i,j}\left(\rho(s)\right) \diff s \right),
\]
where $Y_{ij}(\cdot)$ are independent unit-rate Poisson processes, with $j \in \{1,\ldots,m\}$ and $i \in \{1,\ldots,d\}$. Observe that $Y_{ij}(t)$ can be seen as the cumulative number of arrivals to cell $i$ of type $j$ over the interval $[0,t].$ We denote by $Y(\cdot)$ the  process  with entries $Y_{i,j}(\cdot)$, for $i \in \{0,\ldots,d\}$, $j \in \{1,\ldots,m\}$. 

\subsection{Main results}
\label{sec:MR}

The Markov chain defined in the previous section has a huge state space. As a consequence, direct numerical evaluation of performance metrics is often not feasible. The main idea presented in this section, is to approximate the random objects under study by a suitably chosen Gaussian counterpart. The formal backing of this procedure is given by the scaling limits presented in  \cite[Section 4]{MS2019}, which we briefly summarize here.

We scale the cell lengths by a factor $n$ to ensure enough aggregation for the central limit theorem to kick in, i.e.,  $\ell_i \mapsto n\ell_i$ where $n$ will grow large. Simultaneously, we scale time by a factor $n$, i.e., $t \mapsto nt$, so that that the expected flow of density between cells per unit of time remains invariant. We denote $\rho^n(t) := \rho(nt) / n$, with the cell lengths being $n \ell_i$.

For keep notation light, let $Q(\rho(t))$ be the vector of length $(d+1)m$ with entries $q_{i-1,j}(\rho(t))$, $i \in \{1,\ldots,d+1\}$, $j \in \{1,\ldots,m\}$, ordered lexicographically, i.e., $Q_{(i-1)m+j} = q_{i-i,j}$. We define $H$ to be the $dm \times (d+1)m$ matrix with $H_{kl} := \mathbbm{1}_{\{k=l\}} - \mathbbm{1}_{\{k+m=l\}}$. Finally, we set $L$ as the $dm \times dm$-dimensional diagonal matrix,
with the $k$-th diagonal element being $1/\ell_i$ if $\lceil k/m \rceil = i$, for $k \in \{1,\ldots,dm\}$ and $i \in \{1,\ldots,d\}$.

The first result is a \textit{fluid limit}, which can be seen as a law of large numbers. It states that if  $\rho^n(0) \to \bar{\rho(0)}$, then, almost surely as $n \to \infty$, 
\begin{equation}\label{Eqn: integral equation rho_bar}
	\rho^n(t) \to \bar{\rho}(t) = \bar{\rho}(0) + \int_0^t F(\bar \rho(s)) \diff s, 
\end{equation}
with $F(\rho(t)) := L H Q(\rho(t))$; for the precise statement we refer to \cite[Thm.\ 4.1]{MS2019}.

The second result is a \textit{diffusion limit}, which can be regarded as a central limit theorem. Suppose that  $\lim_{n \to \infty} \sqrt{n}|\rho^n(0) - \rho_0| = 0$ for some $\rho_0$. Then the process
$\hat \rho^n(\cdot)$, defined through $\hat \rho^n(t) := \sqrt{n}(\rho^n(t) - \bar{\rho}(t))$, converges 
 in distribution (as $n \to \infty$) to the process $\hat{\rho}(\cdot)$ solving the stochastic differential equation
\begin{equation}\label{CLT}
	 \hat{\rho}(t) = \hat{\rho}(0) + \int_0^t \partial F(\bar \rho(s)) \hat \rho(t) \diff s + \int_0^t LH \Sigma(\bar \rho(s)) \diff B(s);
\end{equation}
see for details \cite[Thm.\ 4.3]{MS2019}.
Here $\partial F(\rho(t))$ is to be understood as the matrix of weak partial derivative of $F(\bar \rho(t))$, $\Sigma(\bar \rho(s))$ is the diagonal matrix with the square roots of $Q(\bar \rho(t))$ as entries, and $B(\cdot)$ is a length-$(d+1)m$ vector of independent standard Brownian motions.

It is known  that $\hat \rho(\cdot)$, as defined through \eqref{CLT}, is a Gaussian process, with the corresponding  {mean vector} and {covariance matrix} being defined as
\begin{equation*}
	 M(t) := \E[ \rhoh(t) ],\:\:\:\:
		 \Gamma(s,t):=\textrm{cov}\left(\rhoh(s),\rhoh(t)\right) = \E[(\rhoh(s) - M(s))(\rhoh(t)-M(t))^\top].		 
\end{equation*}
The variance of $\rhoh(t)$ is denoted by $V(t) := {\rm var}[\rhoh(t)]=\Gamma(t,t).$
As in \cite[Section 5.6, Problems 6.1, 6.2]{karatzas2012brownian}, with $\bar\Phi(s):=\Phi(s,0)$, these satisfy the explicit expressions
\begin{align}
		 M(t) &= \bar\Phi(t) \left[M(0) + \int_0^t \bar\Phi^{-1}(s) \diff s \right], \nonumber \\
		 \Gamma(s,t) &= \bar\Phi(s) \left[ V(0) + \int_0^{t \wedge s} \bar\Phi^{-1}(u) LH\,\Sigma(\rhob(u))  \left(\bar\Phi^{-1}(u) \,LH\,\Sigma(\rhob(u))\right)^\top \diff u \right] \bar\Phi^\top(t) \label{Eqn: Cov rho_t, rho_s}.
\end{align}
In addition, $M(t)$ and $V(t)$ solve the linear (matrix) differential equations
\begin{align}
	 \frac{{\rm d}M(t)}{{\rm d}t}& = \partial F(\rhob(t)) \,M(t),  \nonumber \\
		 \frac{{\rm d}V(t)}{{\rm d}t} &=  \partial F(\rhob(t)) V(t) + V(t) (\partial F(\rhob(t)))^\top + LH\,\Sigma(\rhob(t)) (LH\,\Sigma(\rhob(t)))^\top .  \label{Eqn: dVar rho}
\end{align}

For the cumulative per-cell arrival process $Y(t)$, similar results have been derived \cite[Section 5]{MS2019}, which we briefly recapitulate now. 
Note that we have 
\begin{equation}\label{Eqn: Y2rho}
	\rho(t) = L\,X(0) + LH\,Y(t),
\end{equation}
where $X(0) := L^{-1} \rho(0)$ is the initial number of vehicles per cell. Note that $Y_{i,j}(\cdot)$ is a counting process with rate $q_{i,j}(\rho(t))$. But, by \eqref{Eqn: Y2rho}, we can also write these rates as  functions $Y(t)$, say $h(Y(t))$. It effectively means that the fluid and diffusion limits that apply to $\rho(\cdot)$, also apply to $Y(\cdot)$. To be precise, if $\lim_{n \to \infty} X^n(0) = X(0)$, then for the sequence of scaled processes $\{Y^n(\cdot)\}_n$ defined through $Y^n(t) := Y(nt) / n$, almost surely as $n \to \infty$,
\[
	Y^n(\cdot)\to \bar Y(\cdot) = \int_0^t h(\bar{Y}(s)) \diff s.
\]
Suppose that $\lim_{n \to \infty} \sqrt{n}|X^n(0) - X(0)| = 0$. Then  the sequence of scaled and centered processes $\{\hat Y^n(\cdot)\}_n$ defined through $\hat Y^n(t) := \sqrt{n} (Y^n(t) - \bar{Y}(t) )$, converges in distribution (as $n\to\infty$) to $\hat Y(\cdot)$ solving the stochastic differential equation
\begin{equation}\label{eq:DL}
	\hat Y(t) = \int_0^t \partial h(\bar{Y}(s)) \hat{Y}(s) \diff s + \int_0^t \bar \Sigma(\bar{Y}(s)) \diff B(s),
\end{equation}
where $\partial h$ is to be understood as the matrix of weak partial derivatives of $h$, and $\bar \Sigma(\bar Y(s))$ is the $(d+1)m \times (d+1)m$ diagonal matrix with entries $h(\bar{Y}(s))$. 

The computation of means and (co-)variances can be done following the same procedure as the one used for the vehicle densities. Moreover, since for any $0 \leqslant t_1 < t_2 < \ldots < t_k < \infty$, for all $k \in \N$, the distribution of the random vector $(Y(t_1,\ldots,Y(t_k))^\top$ is Gaussian, we have that any linear transformation of this vector is again Gaussian. In \cite[Section 5]{MS2019}, we exploit this fact to approximate the distribution of travel times between origins and destinations in the network. For $i \in \{1,\ldots,d\}$ and $k \in \{0,\ldots, d-i-1\}$, we define the type-$j$ travel time, $T_{i,i+k,j}(t)$, as the time that it takes a vehicle of type~$j$ to depart from cell $i+k$, given that it is in cell~$i$ at time $t$. Neglecting the effect of vehicles of the same type overtaking each other (cf.\ \cite[Eqn.\ (40)]{Qian2017}), 
\begin{equation}\label{Eqn: TT tail event}
	\{T_{i,i+k,j} > x\} = \{Y_{i+k,t}(t+x) < Y_{i-1,j}(t)\}, \quad x > 0.
\end{equation}
As such, the probabilities $\PP(T_{i,i+k,j} > x_n)$, with $x_n > 0$ and $n \in \{1,\ldots,N\}$ for some $N \in \N$, can be derived from the random vectors $Y(t), Y(t + x_1), \ldots, Y(t + x_N)$, the joint distribution of which we can approximate using the Gaussian process obtained in the diffusion limit \eqref{eq:DL}. We emphasize that, due to the time-dependent parameters of the distribution of $\hat{Y}(t)$, the distribution of $T_{i,i+k,j}$ is typically not Gaussian.

\subsection{Stationary performance metrics}\label{sec:SPM}

We have constructed the continuous-time Markov process $\rho(\cdot)$. It is irreducible, as a consequence of the fact that the probability to go from any state to the empty state (all cells being empty, that is), in any given amount of time, is strictly positive.
Since the state space ${\mathcal S}$ of $\rho(\cdot)$ is finite, every state is automatically positive recurrent, and thus the Markov process has a unique invariant distribution, which is also a limiting distribution \cite[Sections 3.5 and 3.6]{Norris98}. In other words, $\rho(\cdot)$ is an ergodic Markov chain.

In this subsection we focus on the evaluation of long-term performance metrics. We illustrate our approach by an example. {Since $\rho(\cdot)$ is ergodic, every vehicle entering the road will eventually leave the segment at cell~$d$. Therefore, in the stationary regime, the number of vehicles arriving in the first cell is a measure for the throughput of the segment. In other words, picking $f(\cdot)=q_0(\cdot)$, $f(\rho(t))$ quantifies the throughput of the segment at time $t$.}
The long-term throughput is given by, for $t$ large,
\begin{equation}\label{Eqn: int f(rho(t))}
\int_0^t f(\rho(s)) \diff s.
\end{equation}
Due to $\rho(\cdot)$ being an ergodic Markov chain, the ergodic theorem \cite[Thm.\  3.8.1]{Norris98} applies, and  we thus have that, almost surely as $t\to\infty$,
\begin{equation}\label{Eqn: ergodic law}
	\frac{1}{t}\int_0^t f(\rho(s)) \diff s \to \sum_{x \in \cS} f(x) \pi(x),
\end{equation}
where $\pi(x)$ is the stationary probability of the system being in state $x$. Combining \eqref{Eqn: int f(rho(t))} and \eqref{Eqn: ergodic law}, suggests approximating the throughput in the interval $[0,t]$ by
\begin{equation}\label{Eqn: stationary performance * t}
\int_0^t f(\rho(s)) \diff s \approx	t \sum_{x \in \cS} f(x) \pi(x).
\end{equation}
One usually obtains the stationary probabilities $\pi(x)$ by solving the system's balance equations. This is, however, not a viable option:  it would take prohibitively long due to the process' large state-space. An attractive alternative is to approximate the stationary probabilities, using that $\rho(\cdot)$ is approximately a Gaussian process with mean $\rhob(t)$ and variance $\textrm{var}[\rhoh(t)]$; see the discussion in \cite[Section 6.1]{MS2019}. More concretely, we approximate   $\pi(\cdot)$ by a multivariate Gaussian distribution, having mean vector $\mu$ and covariance matrix $V$ such that 
\begin{align*}
	& F(\mu) = 0 \\
	& \dot{V}(\mu) = \partial F(\mu) V + V \partial F(\mu) + L H \Sigma(\mu) \left( L H \Sigma(\mu) \right)^\top = 0,
\end{align*}
so that, using \eqref{Eqn: integral equation rho_bar} and \eqref{Eqn: dVar rho}, we obtain the stationary distribution of the approximating Gaussian process. We compute the above stationary point $(\mu,V) \in \R^{d+1 \times d}$ using a fixed point iteration by considering the sequence $\{\mu_k\}_k$ and $\{V_k\}_k$ given by
\[
	\mu_{k+1} = \mu_k + F(\mu_k) \Delta t, \qquad V_{k+1} = V_{k} + \dot{V}(\mu_k) \Delta t,
\]
where $\mu_0$ and $V_0$ are, respectively, the length-$d$ zero vector and $d \times d$ zero matrix. In our experiments we took $\Delta t = 0.001$, and we terminated the iteration when the Euclidean distance between $(\mu_{k+1},V_{k+1})$ and {$(\mu_{k},V_{k})$} is below $10^{-9}$.

Once we have the stationary distribution of the Gaussian process, we can approximate the stationary probability $\pi(x)$, $x \in \cS$, by integrating the density of a $N_{dm}(\mu,V)$ random variable over rectangles, using a continuity correction. More precisely, with
\[
	R(x) := \{(y_{11},\ldots,y_{dm}) \colon y_{ij} \in [x_{ij} - 1/({2\ell_i}), x_{ij}+ 1/({2\ell_i})], i \in \{1,\ldots,d\}, j \in \{1,\ldots,m\}\},
\]
denoting a rectangle around $x$, we approximate $\pi(x)$ with (in self-evident notation)
\[
	 \eta(x) :=  \int_{y \in R(x)} \diff N(\mu,V)(y).
\]
Hence, based on the above reasoning our approximation becomes
\begin{equation*}
	\int_0^t f(\rho(s)) \diff s\approx t \sum_{x \in \cS} f(x) \eta(x).
\end{equation*}
Thus far we have been using the example of throughput, but we would like to stress that any other stationary performance metric can be dealt with in the precise same way, with the only restriction to be imposed that $f(\cdot)$ is bounded. Another application of this framework could relate to the quantification of  the pollution level, so as to support environmental policies. Likewise, for maintenance purposes, one could aim at assessing quantities that reflect the deterioration of the road.

\subsection{Numerical example for throughput}
We now illustrate the above procedure for computing stationary metrics by a numerical example. As in the previous section, we take $f(\cdot) = q_0(\cdot)$. This means that we aim at assessing the long-term average rate of vehicles entering the road. 

In our example, we use the single-type MFD and associated discrete flux-function that was introduced by Daganzo \cite{DA95}; cf.\ \cite[Example 3.2]{MS2019} for a short account. As $m=1$, we omit the subscript $j$ in the notation. Denote by $\eta_1(x)$ the probability of the approximating stationary Gaussian distribution corresponding to cell~$1$ being in a state in $[x-1/ (2\ell_i),x+1/(2\ell_i)]$, with $x$ in the state space $\cS_1$ of cell~$1$. In our experiments we compute the stationary throughput rate
\begin{equation}\label{Eqn: stationary throughput rate}
	\sum_{x \in \cS_1} q_0(x) \eta_1(x),
\end{equation}
for various values of $\lambda \in [0,2520]$, and for $\nu = 1200$. In addition, we take the cell lengths $\ell_i$ equal for each $i \in \{1,\ldots,d\}$, and consider various cell lengths to assess their influence on the approximation; $\ell_i \in \{11/\rho^{\textrm{max}},22/\rho^{\textrm{max}},54/\rho^{\textrm{max}},1\}$ km, with the fractions taken so that $X^{\max}_i$ is an integer, for each cell~$i$. As for the remaining parameters, we have chosen $d = 5$, $v^f = 80$ km/h, $w = 16$ km/h, $q^{\textrm{max}} = 1800$ and $\rho^{\textrm{max}} = 108$ veh.

To compare \eqref{Eqn: stationary throughput rate} with its deterministic analog, we also compute the throughput rate that is based on the mean the marginal stationary distribution of cell~1 only: \[q_0\left(\sum_{x \in \cS_1} x \eta_1(x)\right).\] To illustrate the accuracy of the approximation, we also estimate the stationary throughput rate using simulation, with the estimator $t^{-1}\int_0^t q_0(\rho(s)) \diff s$
for $t = 10$, with $\rho(\cdot)$ corresponding to the stationary regime.

We now present our numerical results. In \autoref{Fig: estimates for throughput example}, we have plotted our three estimates as a function of $\lambda$, for our choices of $\ell_i$. We observe that for small $\ell_i$  the estimate based on our Gaussian approximation (`stoch.') is closer to the simulated estimate (`sim.') than the deterministic estimate (`determ.'), indicating that an estimate purely based on the mean does not suffice in {this} model. When $\ell_i$ increases, both estimates get closer. It should be noted, however, that this is due to the fact that the variance of the diffusion limit goes to zero as $\ell_i \to \infty$, essentially being a model-property. 
The kink in the curve of the theoretical estimate  is a numerical effect caused by the stationary regime changing from free-flow to congested. This leads to a shift of the mean of the Gaussian density. In practical terms this kink 
can be easily remedied: by fitting a strictly increasing function to the blue theoretical curve, one obtains a highly accurate approximation of the red simulated curve. 

Finally, the horizontal part of the graphs shows  that the throughput rate is bounded. The truncation level associated to the Gaussian estimate could {be used as a proxy for} the capacity of the road segment.

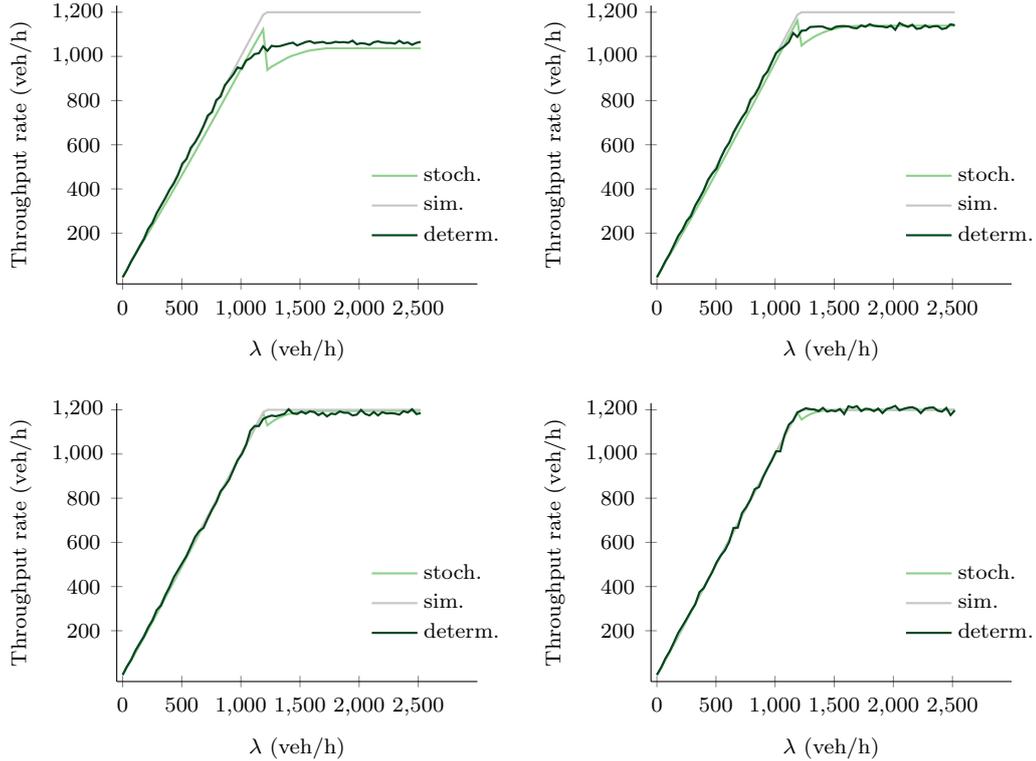
\begin{figure}
\begin{picture}(400 , 300)(0,0)
\put(0,150){
\tikzsetnextfilename{Fig_througput_l=01}
\begin{tikzpicture}	
\begin{axis}[
    xlabel={ $\lambda$ (veh/h)},    ylabel={Throughput rate (veh/h)},
    xmin=-50, xmax=3000,    ymin=-30	, ymax=1230,
    xtick={0,500,1000,1500,2000,2500},    ytick={200,400,600,800,1000,1200},
    xticklabel style={/pgf/number format/fixed, /pgf/number format/precision=2},
    yticklabel style={text width=1.5em,align=right},    
	scaled x ticks=false,
    xtick pos=left,    ytick pos=left,    axis lines=left,
    width = 180pt,    height = 150pt,
    x axis line style=-,    y axis line style=-,
    font = \tiny,
    no markers,
    legend style={at={(1.1,0.1)},anchor=south east, draw=none},
    legend cell align={left},
]
\addplot+[thick, color32] table [x={lambda}, y={th_rate}] {v2Capacity_dataFrame_nu=1200_l=01.dat};
\addplot+[thick, color14] table [x={lambda}, y={mean_rates}] {v2Capacity_dataFrame_nu=1200_l=01.dat};
\addplot+[thick, color36] table [x={lambda}, y={sim_rates}] {v2Capacity_dataFrame_nu=1200_l=01.dat};
\legend{stoch., sim., determ.}
\end{axis}
\end{tikzpicture}
}
\put(200,150){
\tikzsetnextfilename{Fig_througput_l=02}
\begin{tikzpicture}	
\begin{axis}[
    xlabel={ $\lambda$ (veh/h)},    ylabel={Throughput rate (veh/h)},
    xmin=-50, xmax=3000,    ymin=-30	, ymax=1230,
    xtick={0,500,1000,1500,2000,2500},    ytick={200,400,600,800,1000,1200},
    xticklabel style={/pgf/number format/fixed, /pgf/number format/precision=2},
    yticklabel style={text width=1.5em,align=right},    
	scaled x ticks=false,
    xtick pos=left,    ytick pos=left,    axis lines=left,
    width = 180pt,    height = 150pt,
    x axis line style=-,    y axis line style=-,
    font = \tiny,
    no markers,
    legend style={at={(1.1,0.1)},anchor=south east, draw=none},
    legend cell align={left},
]
\addplot+[thick, color32] table [x={lambda}, y={th_rate}] {v2Capacity_dataFrame_nu=1200_l=02.dat};
\addplot+[thick, color14] table [x={lambda}, y={mean_rates}] {v2Capacity_dataFrame_nu=1200_l=02.dat};
\addplot+[thick, color36] table [x={lambda}, y={sim_rates}] {v2Capacity_dataFrame_nu=1200_l=02.dat};
\legend{stoch., sim., determ.}
\end{axis}
\end{tikzpicture}
}
\put(0,0){
\tikzsetnextfilename{Fig_througput_l=05}
\begin{tikzpicture}	
\begin{axis}[
    xlabel={ $\lambda$ (veh/h)},    ylabel={Throughput rate (veh/h)},
    xmin=-50, xmax=3000,    ymin=-30	, ymax=1230,
    xtick={0,500,1000,1500,2000,2500},    ytick={200,400,600,800,1000,1200},
    xticklabel style={/pgf/number format/fixed, /pgf/number format/precision=2},
    yticklabel style={text width=1.5em,align=right},    
	scaled x ticks=false,
    xtick pos=left,    ytick pos=left,    axis lines=left,
    width = 180pt,    height = 150pt,
    x axis line style=-,    y axis line style=-,
    font = \tiny,
    no markers,
    legend style={at={(1.1,0.1)},anchor=south east, draw=none},
    legend cell align={left},
]
\addplot+[thick, color32] table [x={lambda}, y={th_rate}] {v2Capacity_dataFrame_nu=1200,0_l=0,50.dat};
\addplot+[thick, color14] table [x={lambda}, y={mean_rates}] {v2Capacity_dataFrame_nu=1200,0_l=0,50.dat};
\addplot+[thick, color36] table [x={lambda}, y={sim_rates}] {v2Capacity_dataFrame_nu=1200,0_l=0,50.dat};
\legend{stoch., sim., determ.}
\end{axis}
\end{tikzpicture}
}
\put(200,0){
\tikzsetnextfilename{Fig_througput_l=10}
\begin{tikzpicture}	
\begin{axis}[
    xlabel={ $\lambda$ (veh/h)},    ylabel={Throughput rate (veh/h)},
    xmin=-50, xmax=3000,    ymin=-30	, ymax=1230,
    xtick={0,500,1000,1500,2000,2500},    ytick={200,400,600,800,1000,1200},
    xticklabel style={/pgf/number format/fixed, /pgf/number format/precision=2},
    yticklabel style={text width=1.5em,align=right},    
	scaled x ticks=false,
    xtick pos=left,    ytick pos=left,    axis lines=left,
    width = 180pt,    height = 150pt,
    x axis line style=-,    y axis line style=-,
    font = \tiny,
    no markers,
    legend style={at={(1.1,0.1)},anchor=south east, draw=none},
    legend cell align={left},
]
\addplot+[thick, color32] table [x={lambda}, y={th_rate}] {v2Capacity_dataFrame_nu=1200,0_l=1,00.dat};
\addplot+[thick, color14] table [x={lambda}, y={mean_rates}] {v2Capacity_dataFrame_nu=1200,0_l=1,00.dat};
\addplot+[thick, color36] table [x={lambda}, y={sim_rates}] {v2Capacity_dataFrame_nu=1200,0_l=1,00.dat};
\legend{stoch., sim., determ.}
\end{axis}
\end{tikzpicture}
}
\end{picture}
\caption{Graphs of our approximation (stoch.), point estimate (determ.), and simulated estimate (sim.) for the stationary throughput rate, as a function of $\lambda$, for $\nu=1200$, and $\ell_i = 11/108$ (top left), $\ell_i = 22/108$ (top right), $\ell_i = 54/108$ (bottom left), and $\ell_i = 108/108$ (bottom right).}
\label{Fig: estimates for throughput example}
\end{figure}

\section{Route choice in stochastic environments}
\label{sec: route choice}


Nowadays, navigation software is intensively used for  route selection purposes. Often, the route that minimizes  travel time is  seen as the optimal route. Importantly, in such algorithms travel times are typically treated as deterministic quantities. As we argued, however, road traffic systems are inherently stochastic. This could mean that, for instance, when choosing between two alternatives, the route with the shorter expected travel time has the larger standard deviation. In such a situation it is up to the driver to make a choice: depending on his personal preferences (in terms of risk aversion) and the importance of the specific trip, she will choose the best alternative. A convenient framework enabling such decision uses the concept of {\it utility functions}  \cite{SEN}; see also {e.g.,} \cite{LOUI,SIVA}. Such a utility function could encompass both mean and standard deviation of the travel time, but in addition also for instance the 95\%-quantile. 

Based on the above, it makes sense to rely on a modeling framework in which 
travel times are represented by random variables. The setup of \cite{MS2019}, as summarized in Sections \ref{sec:MS}--\ref{sec:MR}, is particularly suitable for these purposes. It provides us with an accurate approximation of the travel-time distribution, covering all possible arguments of the utility function (in particular moments and quantiles).  In this section, we give an example that illustrates how different routes can be considered optimal, depending on the driver's specific preferences. 

\subsection{Route choice example}
As argued above, a driver does not only care about his expected time of arrival, but also about the corresponding uncertainty. 
In general, she wishes to minimize a utility, which is a function of various features of the travel-time distribution.
In our example, we consider the specific example of a road user with a utility function that is of the form $\mu + c \sigma$, for some constant $c > 0$, where $\mu$ and $\sigma$ are the expectation and standard deviation of the travel-time  of the chosen route. The situation is considered in which the road user can choose to drive to a destination using two different roads, both having identical characteristics, except for the free-flow (maximum) velocity. 

More concretely, let both roads consist of $d = 3$ cells, all of length $\ell_i = 1$ km. We take the Daganzo MFD again for the dynamics of vehicles, with two parameter settings. For the first set of parameters, we take $w = 16$ km/h, $\rho^{\textrm{max}} = 108$ veh/km, $q^{\textrm{max}} = 1500$ veh/h, $\lambda = 1400$ veh/h, and $\nu = 1500$ veh/h, with the free-flow velocity of route~1 being $v^f = 90$ km/h, and for route~2 $v^f = 80$ km/h. For the second setting, we take $w = 20$ km/h, $\rho^{\textrm{max}} = 108$ veh/km, $q^{\textrm{max}} = 1800$ veh/h, $\lambda = 1700$ veh/h, and $\nu = 1800$ veh/h, with $v^f = 120$ km/h for route~1, and $v^f = 110$ km/h for route~2.

For both routes, we compute the travel-time distribution by \eqref{Eqn: TT tail event}, using the Gaussian approximation of the quantities featuring in the event on the right-hand-side. The corresponding mean vector and covariance matrix can be found relying on the analog of equations \eqref{Eqn: dVar rho} for $\hat{Y}(\cdot)$. In this example, we initialize the differential equations for the mean vectors  of $\hat{Y}(\cdot)$ by the stationary mean of both routes.

For the initialization of the covariance matrix, we want to model an initial state in which the user is less certain about the state of route~1 than that of route~2. This higher uncertainty of route 1 can be, for example, a consequence of the presence of a traffic light at its entrance: the route can be congested due to the traffic light having been red for some time and turning green when the driver arrives, but it can be quiet as well due to the traffic light having been   green for a while. We model the higher uncertainty of route 1, in the differential equation describing the covariance matrix, by initializing it 
for route~1 by the diagonal matrix of the stationary mean vector divided by $b \in \{1.5,2,2.5\}$, whereas for route~2 we follow the same procedure but divide by 5.

\subsection{Numerical results}
We evaluate \eqref{Eqn: TT tail event} for $i=1$ and $k = d$, so as to obtain  these tail probabilities  at any time between $0$ and $480$ s. Using this numerical output, we compute the mean $\mu_k$ and the standard deviation $\sigma_k$ of both  travel-time distributions. For routes $k = 1,2$ in setting~1, we find (in units of seconds)

\[
	\mu_2 = 135.86, \quad \sigma_2 = 13.87,
\]
and
\[
	\begin{aligned}
		& b_1 = 2/3 : & \mu_1 = 121.56, & \;\;\; \sigma_1 = 25.43, \\
		& b_1 = 1/2 : & \mu_1 = 121.23, & \;\;\; \sigma_1 = 20.33, \\
		& b_1 = 2/5 : & \mu_1 = 121.08, & \;\;\; \sigma_1 = 17.51,
	\end{aligned}
\]
whereas for setting~2 we find,
\[
	\mu_2 = 98.93, \quad \sigma_2 = 10.56,
\]
and
\[
	\begin{aligned}
		& b_1 = 2/3 : & \mu_1 = 91.27, & \;\;\; \sigma_1 = 19.27, \\
		& b_1 = 1/2 : & \mu_1 = 91.03, & \;\;\; \sigma_1 = 15.49, \\
		& b_1 = 2/5 : & \mu_1 = 90.91, & \;\;\; \sigma_1 = 13.40.
	\end{aligned}
\]

Observe that in all cases, route~1 has the smaller mean, but the higher standard deviation. In other words: it depends on the driver's specific preference which route she will choose. 
In terms of utility functions, the route minimizing $\mu_k + c \sigma_k$ will be selected. In \autoref{Fig: route choice setting  1} and \autoref{Fig: route choice setting 2}, we have plotted both of the utility functions, for settings 1 and 2, respectively, as a function of $c$. Depending on the driver's $c$-value,  she picks route~1 or route~2.

\begin{figure}
\begin{picture}(450 , 100)(0,0)
\put(0,0){
\tikzsetnextfilename{Fig_routeChoice_setting1_b1=0667}
\begin{tikzpicture}	
\begin{axis}[
    xlabel={c},    ylabel={$\mu_k + c \sigma_k$ (s)},
    xmin=-0.2, xmax=3.2,    ymin=85, ymax=205,
    xtick={0,1,2,3},    ytick={100,120,140,160,180,200},
    xticklabel style={/pgf/number format/.cd,fixed,set thousands separator={}},
    yticklabel style={/pgf/number format/.cd,fixed,set thousands separator={}},
	scaled x ticks=false,
	scaled y ticks=false,
    xtick pos=left,    ytick pos=left,    axis lines=left,
    width = 150pt,    height = 100pt,
    no markers,    x axis line style=-,    y axis line style=-,
    font = \tiny,
    legend style={at={(1.1,0.55)},anchor=north east, draw=none},
    legend cell align={left},
]
\addplot+[thick,solid,mark=none,color32] coordinates {(0,121.56) (3,197.85)};
\addplot+[thick,solid,mark=none,color36] coordinates {(0,135.86) (3,177.47)};
\legend{$k=1$,$k=2$}
\end{axis}
\end{tikzpicture}
}
\put(150,0){
\tikzsetnextfilename{Fig_routeChoice_setting1_b1=05}
\begin{tikzpicture}	
\begin{axis}[
    xlabel={c},    ylabel={$\mu_k + c \sigma_k$ (s)},
    xmin=-0.2, xmax=3.2,    ymin=85, ymax=205,
    xtick={0,1,2,3},    ytick={100,120,140,160,180,200},
    xticklabel style={/pgf/number format/.cd,fixed,set thousands separator={}},
    yticklabel style={/pgf/number format/.cd,fixed,set thousands separator={}},
	scaled x ticks=false,
	scaled y ticks=false,
    xtick pos=left,    ytick pos=left,    axis lines=left,
    width = 150pt,    height = 100pt,
    no markers,    x axis line style=-,    y axis line style=-,
    font = \tiny,
    legend style={at={(1.1,0.55)},anchor=north east, draw=none},
    legend cell align={left},
]
\addplot+[thick,solid,mark=none,color32] coordinates {(0,121.23) (3,182.22)};
\addplot+[thick,solid,mark=none,color36] coordinates {(0,135.86) (3,177.47)};
\legend{$k=1$,$k=2$}
\end{axis}
\end{tikzpicture}
}
\put(300,0){
\tikzsetnextfilename{Fig_routeChoice_setting1_b1=04}
\begin{tikzpicture}	
\begin{axis}[
    xlabel={c},    ylabel={$\mu_k + c \sigma_k$ (s)},
    xmin=-0.2, xmax=3.2,    ymin=85, ymax=205,
    xtick={0,1,2,3},    ytick={100,120,140,160,180,200},
    xticklabel style={/pgf/number format/.cd,fixed,set thousands separator={}},
    yticklabel style={/pgf/number format/.cd,fixed,set thousands separator={}},
	scaled x ticks=false,
	scaled y ticks=false,
    xtick pos=left,    ytick pos=left,    axis lines=left,
    width = 150pt,    height = 100pt,
    no markers,    x axis line style=-,    y axis line style=-,
    font = \tiny,
    legend style={at={(1.1,0.55)},anchor=north east, draw=none},
    legend cell align={left},
]
\addplot+[thick,solid,mark=none,color32] coordinates {(0,121.08) (3,173.61)};
\addplot+[thick,solid,mark=none,color36] coordinates {(0,135.86) (3,177.47)};
\legend{$k=1$,$k=2$}
\end{axis}
\end{tikzpicture}
}
\end{picture}
\caption{Utility functions $\mu_k + c \sigma_k$ in setting~1, plotted as function of $c$, for routes $k=1,2$, with initial standard deviation factor of road~1 (from left to right) $b_1 = 2/3, 1/2, 2/5$.}
\label{Fig: route choice setting 1}
\end{figure}
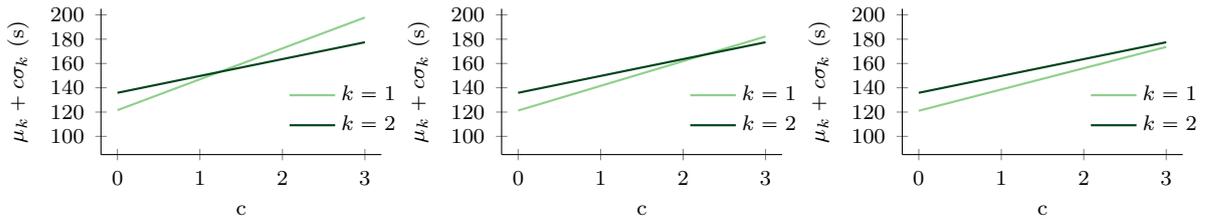

\begin{figure}
\begin{picture}(450 , 100)(0,0)
\put(0,0){
\tikzsetnextfilename{Fig_routeChoice_setting2_b1=0667}
\begin{tikzpicture}	
\begin{axis}[
    xlabel={c},    ylabel={$\mu_k + c \sigma_k$ (s)},
    xmin=-0.2, xmax=3.2,    ymin=85, ymax=195,
    xtick={0,1,2,3},    ytick={90,110,130,150,170,190},
    xticklabel style={/pgf/number format/.cd,fixed,set thousands separator={}},
    yticklabel style={/pgf/number format/.cd,fixed,set thousands separator={}},
	scaled x ticks=false,
	scaled y ticks=false,
    xtick pos=left,    ytick pos=left,    axis lines=left,
    width = 150pt,    height = 100pt,
    no markers,    x axis line style=-,    y axis line style=-,
    font = \tiny,
    legend style={at={(0.1,0.975)},anchor=north west, draw=none},
    legend cell align={left},
]
\addplot+[thick,solid,mark=none,color32] coordinates {(0,91.27) (3,149.08)};
\addplot+[thick,solid,mark=none,color36] coordinates {(0,98.93) (3,130.61)};
\legend{$k=1$,$k=2$}
\end{axis}
\end{tikzpicture}
}
\put(150,0){
\tikzsetnextfilename{Fig_routeChoice_setting2_b1=05}
\begin{tikzpicture}	
\begin{axis}[
    xlabel={c},    ylabel={$\mu_k + c \sigma_k$ (s)},
    xmin=-0.2, xmax=3.2,    ymin=85, ymax=195,
    xtick={0,1,2,3},    ytick={90,110,130,150,170,190},
    xticklabel style={/pgf/number format/.cd,fixed,set thousands separator={}},
    yticklabel style={/pgf/number format/.cd,fixed,set thousands separator={}},
	scaled x ticks=false,
	scaled y ticks=false,
    xtick pos=left,    ytick pos=left,    axis lines=left,
    width = 150pt,    height = 100pt,
    no markers,    x axis line style=-,    y axis line style=-,
    font = \tiny,
    legend style={at={(0.1,0.975)},anchor=north west, draw=none},
    legend cell align={left},
]
\addplot+[thick,solid,mark=none,color32] coordinates {(0,91.03) (3,137.5)};
\addplot+[thick,solid,mark=none,color36] coordinates {(0,98.93) (3,130.61)};
\legend{$k=1$,$k=2$}
\end{axis}
\end{tikzpicture}
}
\put(300,0){
\tikzsetnextfilename{Fig_routeChoice_setting2_b1=04}
\begin{tikzpicture}	
\begin{axis}[
    xlabel={c},    ylabel={$\mu_k + c \sigma_k$ (s)},
    xmin=-0.2, xmax=3.2,    ymin=85, ymax=195,
    xtick={0,1,2,3},    ytick={90,110,130,150,170,190},
    xticklabel style={/pgf/number format/.cd,fixed,set thousands separator={}},
    yticklabel style={/pgf/number format/.cd,fixed,set thousands separator={}},
	scaled x ticks=false,
	scaled y ticks=false,
    xtick pos=left,    ytick pos=left,    axis lines=left,
    width = 150pt,    height = 100pt,
    no markers,    x axis line style=-,    y axis line style=-,
    font = \tiny,
    legend style={at={(0.1,0.975)},anchor=north west, draw=none},
    legend cell align={left},
]
\addplot+[thick,solid,mark=none,color32] coordinates {(0,90.91) (3,131.11)};
\addplot+[thick,solid,mark=none,color36] coordinates {(0,98.93) (3,130.61)};
\legend{$k=1$,$k=2$}
\end{axis}
\end{tikzpicture}
}
\end{picture}
\caption{Utility functions $\mu_k + c \sigma_k$ in setting~2, plotted as function of $c$, for routes $k=1,2$, with initial standard deviation factor of road~1 (from left to right) $b_1 = 2/3, 1/2, 2/5$.}
\label{Fig: route choice setting 2}
\end{figure}
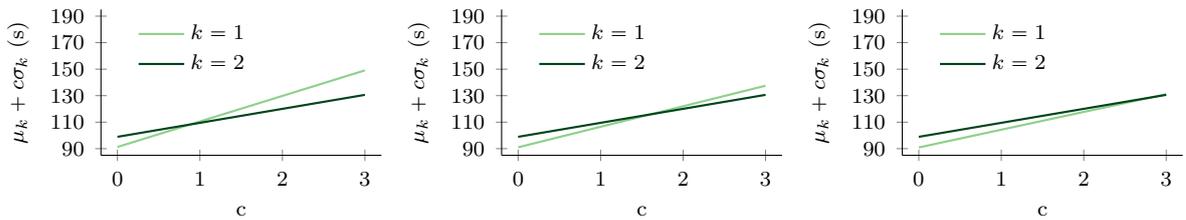

%
%

\section{Control of traffic flows}
\label{sec: control}

Arguably the ultimate goal of road traffic modeling concerns the design of traffic control policies. One often pursues  optimizing efficiency-related performance metrics such as travel-time and throughput, while taking into account constraints, e.g., relating to the policy's environmental impact. Simply put, a control policy selects values for the system's {controllable parameters}. These include, for example, the speed limit (which can be enforced by the use of matrix signs or another information system), the number of lanes (which is essentially a road network design feature), and the arrival rate of specific vehicle types (which can be enforced by applying ramp metering). Observe that some of these mechanisms can be seen as online control, whereas others are of an offline, design-related nature; see  e.g.\ \cite{papa} for a broad account of control mechanisms.

To soundly make decisions, one needs a modeling framework that is capable of quantifying the system's performance as a function of the {control parameters} mentioned above. A few typical examples are: `By how much does the throughput of a road increase when opening an additional lane?', `How does a reduction of the maximum speed impact the travel time?', and `By how much should traffic  be slowed down in case of a traffic jam?' The objective function is an integral part of the decision making. It could be based on average quantities, but it could feature e.g., standard deviations and quantiles as well, making our stochastic framework particularly useful. One could say that the problem of finding an optimal policy is parametrized by the quantities appearing in the objective function. If, for instance, we are working in a setting in which we express the system's performance through an $\alpha$-quantile (for instance of the travel-time distribution), then picking $\alpha =95\%$ leads to conservative policies in which worst-case events have a relatively large weight, whereas picking $\alpha=50\%$ puts more weight on the system's average behavior.

In this section we illustrate that the model in \cite{MS2019} is well-suited to evaluate the impact of control {parameters}. We do so by considering the travel time as a function of maximum velocity, the number of lanes, the arrival rate, and the percentage of trucks. To include the percentage of trucks as a control {parameter}, we require a multi-class MFD; note that this multi-class setting is covered by \cite{MS2019}.
We first present our experimental setup for our illustration, and then present the numerical results. We finish this section by discussing several extensions.

\subsection{Experimental setup}
We take the two-class MFD from \cite{cb2003}, which gives the joint flow of both vehicle classes as a function of their respective densities. Here vehicles from the second class are longer, and thus take more space than vehicles from the first class, and have a lower velocity. We refer to, e.g., \cite[Example 3.3]{MS2019} for more background. We refer to vehicles of class~1 as \textit{cars} and to vehicles of class~2 as \textit{trucks}. As parameters, we take those used in the numerical experiments of \cite{cb2003}, with a kilometer taken as unit length. More precisely, we take $v^f_1 = 108$~km/h, $v^f_2 = 79.2$~km/h, $v^{\rm c} = 61.2$~km/h, $L_1 = 0.0065$~km, $L_2 = 0.0165$~km, $N = 3$ lanes, and $\beta = 0.25$ as a parameter distinguishing the so-called free-flow and congestion regime. These parameters are our `base set', that is, we use the above values when the parameter is \textit{not} a control {parameter}.

In our experiment we consider a road segment of $d = 10$ cells, each of length $1$~km. Unless we specifically study the impact of the arrival rate $\lambda$,  we take $\lambda = \lambda_1 + \lambda_2 = 1200$ veh/h. In this scenario,  the road segment is moderately loaded, but clearly the experiments in this section can also be performed if the congestion level is higher. The fraction of trucks is denoted by $b \in [0,1]$, and is taken as $0.2$ (unless we specifically study its impact), so that $\lambda_1 = (1 - b) \lambda$ and $\lambda_2 = b \lambda$. Moreover, for the boundary condition at cell~10, we take $\nu_i$ equal to $2/3$ times the maximum possible flow for vehicle class~$i$, $i \in \{1,2\}$.

Using our approximations we evaluate the mean and standard deviation of the travel-time distribution for the entire segment, i.e., the time that it takes a vehicle that just entered cell~1, to depart from cell~10. In our experiment we consider the {evolution of the system, starting from its stationary mean} at time $0$ (but, evidently, any other situation could be considered as well). The travel-time distribution is computed by evaluating \eqref{Eqn: TT tail event}, for $i = 1$, $k = 9$ and for 1001 equidistant times in $[0,1400]$, measured in seconds. 

In our experiments we focus on the following control parameters:
\begin{itemize}
\item[$\circ$]
We first study the impact of  the maximum velocities, i.e., $v^f_1$ and $v^f_2$. We vary them between $[50.4,118.8]$~km/h,  corresponding to velocities between 14 and $33$~m/s. In addition, the maximum velocity of trucks is bounded at $79.2$~km/h. 
\item[$\circ$]
In the second place, we assess the impact of  the number of lanes $N$, for which we take values in $\{1,2,3,4,5\}$. The number of lanes is an offline control {parameter}, with the purpose of creating more room for vehicles, to reduce the  congestion level. 
\item[$\circ$] 
Thirdly, we evaluate the impact of the arrival rate $\lambda$, which we vary between 1000 and 4800, using $b$ to compute the arrival rates for cars and trucks. Insight into this sensitivity can be used in e.g., ramp-metering; cf.\ for instance \cite{papa}. 
\item[$\circ$]
Last, we study the impact of $b$ itself, for which we take values in $\{0.05, 0.1, 0.2\}$, and for each of these values we perform all three other experiments. These experiments shed light on the impact of the traffic mix (e.g., an increase of the fraction of vehicles that are trucks) on the performance metric of interest. 
\end{itemize}

\subsection{Results of control experiments}

We now present the results of the above control experiments. We present the results for the three control parameters (maximum velocity, number of lanes, and arrival rates) by plotting the mean and standard deviation of the travel-time distribution as a function of this control parameter, with a separate figure for each value of $b$.

\begin{figure}
\begin{picture}(450 , 200)(0,0)
\put(0,100){
\tikzsetnextfilename{Fig_control_Vs_perc20_mean}
\begin{tikzpicture}	
\begin{axis}[
    xlabel={ $v^f$ (m/s)},    ylabel={$\E [T_{1,10}(0)]$ (s)},
    xmin=13, xmax=36,    ymin=280	, ymax=720,
    xtick={15,20,25,30,35},    ytick={300,400,500,600,700},
    xticklabel style={/pgf/number format/fixed, /pgf/number format/precision=2},
    yticklabel style={text width=1.1em,align=right},    
	scaled x ticks=false,
    xtick pos=left,    ytick pos=left,    axis lines=left,
    width = 150pt,    height = 100pt,
    x axis line style=-,    y axis line style=-,
    font = \tiny,
    legend style={at={(1.1,0.975)},anchor=north east, draw=none},
    legend cell align={left},
]
\addplot+[only marks, color32,mark=*,mark options={fill=color36,scale=0.5}] table [x={Vs}, y={Mean_0_perc=20}] {pd2_Section5_Vs_df.dat};
\addplot+[only marks, color36,mark=triangle*,mark options={fill=color36,scale=0.5}] table [x={Vs}, y={Mean_1_perc=20}] {pd2_Section5_Vs_df.dat};
\legend{Cars,Trucks}
\end{axis}
\end{tikzpicture}
}
\put(150,100){
\tikzsetnextfilename{Fig_control_Vs_perc10_mean}
\begin{tikzpicture}	
\begin{axis}[
    xlabel={ $v^f$ (m/s)},    ylabel={$\E [T_{1,10}(0)]$ (s)},
    xmin=13, xmax=36,    ymin=280	, ymax=720,
    xtick={15,20,25,30,35},    ytick={300,400,500,600,700},
    xticklabel style={/pgf/number format/fixed, /pgf/number format/precision=2},
    yticklabel style={text width=1.1em,align=right},    
	scaled x ticks=false,
    xtick pos=left,    ytick pos=left,    axis lines=left,
    width = 150pt,    height = 100pt,
    x axis line style=-,    y axis line style=-,
    font = \tiny,
    legend style={at={(1.1,0.975)},anchor=north east, draw=none},
    legend cell align={left},
]
\addplot+[only marks, color32,mark=*,mark options={fill=color36,scale=0.5}] table [x={Vs}, y={Mean_0_perc=10}] {pd2_Section5_Vs_df.dat};
\addplot+[only marks, color36,mark=triangle*,mark options={fill=color36,scale=0.5}] table [x={Vs}, y={Mean_1_perc=10}] {pd2_Section5_Vs_df.dat};
\legend{Cars,Trucks}
\end{axis}
\end{tikzpicture}
}
\put(300,100){
\tikzsetnextfilename{Fig_control_Vs_perc5_mean}
\begin{tikzpicture}	
\begin{axis}[
    xlabel={ $v^f$ (m/s)},    ylabel={$\E [T_{1,10}(0)]$ (s)},
    xmin=13, xmax=36,    ymin=280	, ymax=720,
    xtick={15,20,25,30,35},    ytick={300,400,500,600,700},
    xticklabel style={/pgf/number format/fixed, /pgf/number format/precision=2},
    yticklabel style={text width=1.1em,align=right},    
	scaled x ticks=false,
    xtick pos=left,    ytick pos=left,    axis lines=left,
    width = 150pt,    height = 100pt,
    x axis line style=-,    y axis line style=-,
    font = \tiny,
    legend style={at={(1.1,0.975)},anchor=north east, draw=none},
    legend cell align={left},
]
\addplot+[only marks, color32,mark=*,mark options={fill=color36,scale=0.5}] table [x={Vs}, y={Mean_0_perc=5}] {pd2_Section5_Vs_df.dat};
\addplot+[only marks, color36,mark=triangle*,mark options={fill=color36,scale=0.5}] table [x={Vs}, y={Mean_1_perc=5}] {pd2_Section5_Vs_df.dat};
\legend{Cars,Trucks}
\end{axis}
\end{tikzpicture}
}
\put(0,0){
\tikzsetnextfilename{Fig_control_Vs_perc20_std}
\begin{tikzpicture}	
\begin{axis}[
    xlabel={ $v^f$ (m/s)},    ylabel={$\sigma_{T_{1,10}(0)}$ (s)},
    xmin=13, xmax=36,    ymin=25	, ymax=115,
    xtick={15,20,25,30,35},    ytick={30,50,70,90,110},
    xticklabel style={/pgf/number format/fixed, /pgf/number format/precision=2},
    yticklabel style={text width=1.1em,align=right},    
	scaled x ticks=false,
    xtick pos=left,    ytick pos=left,    axis lines=left,
    width = 150pt,    height = 100pt,
    x axis line style=-,    y axis line style=-,
    font = \tiny,
    legend style={at={(1.1,0.975)},anchor=north east, draw=none},
    legend cell align={left},
]
\addplot+[only marks, color32,mark=*,mark options={fill=color36,scale=0.5}] table [x={Vs}, y={Std_0_perc=20}] {pd2_Section5_Vs_df.dat};
\addplot+[only marks, color36,mark=triangle*,mark options={fill=color36,scale=0.5}] table [x={Vs}, y={Std_1_perc=20}] {pd2_Section5_Vs_df.dat};
\legend{Cars,Trucks}
\end{axis}
\end{tikzpicture}
}
\put(150,0){
\tikzsetnextfilename{Fig_control_Vs_perc10_std}
\begin{tikzpicture}	
\begin{axis}[
    xlabel={ $v^f$ (m/s)},    ylabel={$\sigma_{T_{1,10}(0)}$ (s)},
    xmin=13, xmax=36,    ymin=25	, ymax=115,
    xtick={15,20,25,30,35},    ytick={30,50,70,90,110},
    xticklabel style={/pgf/number format/fixed, /pgf/number format/precision=2},
    yticklabel style={text width=1.1em,align=right},    
	scaled x ticks=false,
    xtick pos=left,    ytick pos=left,    axis lines=left,
    width = 150pt,    height = 100pt,
    x axis line style=-,    y axis line style=-,
    font = \tiny,
    legend style={at={(1.1,0.975)},anchor=north east, draw=none},
    legend cell align={left},
]
\addplot+[only marks, color32,mark=*,mark options={fill=color36,scale=0.5}] table [x={Vs}, y={Std_0_perc=10}] {pd2_Section5_Vs_df.dat};
\addplot+[only marks, color36,mark=triangle*,mark options={fill=color36,scale=0.5}] table [x={Vs}, y={Std_1_perc=10}] {pd2_Section5_Vs_df.dat};
\legend{Cars,Trucks}
\end{axis}
\end{tikzpicture}
}
\put(300,0){
\tikzsetnextfilename{Fig_control_Vs_perc5_std}
\begin{tikzpicture}	
\begin{axis}[
    xlabel={ $v^f$ (m/s)},    ylabel={$\sigma_{T_{1,10}(0)}$ (s)},
    xmin=13, xmax=36,    ymin=25	, ymax=115,
    xtick={15,20,25,30,35},    ytick={30,50,70,90,110},
    xticklabel style={/pgf/number format/fixed, /pgf/number format/precision=2},
    yticklabel style={text width=1.1em,align=right},    
	scaled x ticks=false,
    xtick pos=left,    ytick pos=left,    axis lines=left,
    width = 150pt,    height = 100pt,
    x axis line style=-,    y axis line style=-,
    font = \tiny,
    legend style={at={(1.1,0.7)},anchor=north east, draw=none},
    legend cell align={left},
]
\addplot+[only marks, color32,mark=*,mark options={fill=color36,scale=0.5}] table [x={Vs}, y={Std_0_perc=5}] {pd2_Section5_Vs_df.dat};
\addplot+[only marks, color36,mark=triangle*,mark options={fill=color36,scale=0.5}] table [x={Vs}, y={Std_1_perc=5}] {pd2_Section5_Vs_df.dat};
\legend{Cars,Trucks}
\end{axis}
\end{tikzpicture}
}
\end{picture}
\caption{Mean ($\E [T_{1,10}(0)]$) and standard deviation ($\sigma_{T_{1,10}}(0)$) of the approximation to the travel-time distribution, plotted as function of $v^f$ (the free-flow velocity), for $b=0.2$ (left), $b=0.1$ (middle) and $b=0.05$ (right).}
\label{Fig: control met v}
\end{figure}
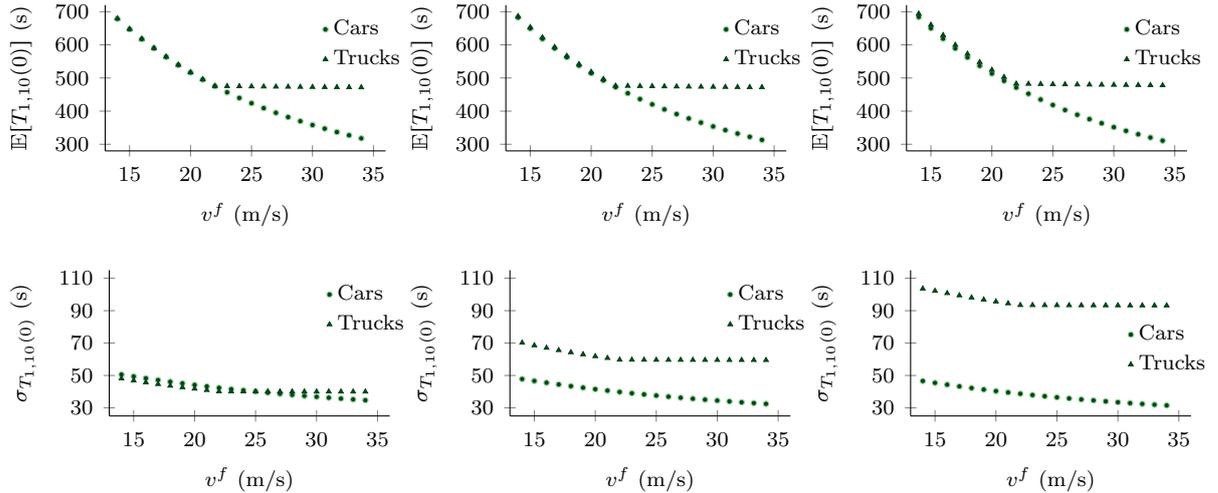

We first focus on the impact of  the maximum velocity. \autoref{Fig: control met v} shows the the mean and standard deviation of the travel time as a function of the maximum velocity. As expected, the mean  and the standard deviation are non-increasing, with the mean being hardly affected by $b$. The standard deviation of the number of cars (trucks) is decreasing (increasing) in $b$: when there are more vehicles involved there is more averaging, such that the variance goes down, and vice versa. 

Then we consider the impact of the number of lanes; see \autoref{Fig: control met N}. When $N$ is already relatively large, the effect of adding an extra lane is modest. This is mainly due to the fact that when there is enough space, vehicles drive at their maximum velocity (which is an intrinsic feature of {the} model). 
In addition, the mean is just very mildly affected by $b$, whereas the standard deviation {\it does} substantially decrease in $b$.

\begin{figure}
\begin{picture}(450 , 200)(0,0)
\put(0,100){
\tikzsetnextfilename{Fig_control_Ns_perc20_mean}
\begin{tikzpicture}	
\begin{axis}[
    xlabel={  \# lanes},    ylabel={$\E [T_{1,10}(0)]$ (s)},
    xmin=0.5, xmax=5.5,    ymin=95, ymax=210,
    xtick={1,2,3,4,5},    ytick={100,120,140,160,180,200},
    xticklabel style={/pgf/number format/fixed, /pgf/number format/precision=2},
    yticklabel style={text width=1.1em,align=right},    
	scaled x ticks=false,
    xtick pos=left,    ytick pos=left,    axis lines=left,
    width = 150pt,    height = 100pt,
    x axis line style=-,    y axis line style=-,
    font = \tiny,
    legend style={at={(1.1,0.975)},anchor=north east, draw=none},
    legend cell align={left},
]
\addplot+[only marks, color32,mark=*,mark options={fill=color36,scale=0.5}] table [x={Ns}, y={Mean_0_perc=20}] {pd_Section5_Ns_df.dat};
\addplot+[only marks, color36,mark=triangle*,mark options={fill=color36,scale=0.5}] table [x={Ns}, y={Mean_1_perc=20}] {pd_Section5_Ns_df.dat};
\legend{Cars,Trucks}
\end{axis}
\end{tikzpicture}
}
\put(150,100){
\tikzsetnextfilename{Fig_control_Ns_perc10_mean}
\begin{tikzpicture}	
\begin{axis}[
    xlabel={  \# lanes},    ylabel={$\E [T_{1,10}(0)]$ (s)},
    xmin=0.5, xmax=5.5,    ymin=95, ymax=210,
    xtick={1,2,3,4,5},    ytick={100,120,140,160,180,200},
    xticklabel style={/pgf/number format/fixed, /pgf/number format/precision=2},
    yticklabel style={text width=1.1em,align=right},    
	scaled x ticks=false,
    xtick pos=left,    ytick pos=left,    axis lines=left,
    width = 150pt,    height = 100pt,
    x axis line style=-,    y axis line style=-,
    font = \tiny,
    legend style={at={(1.1,0.975)},anchor=north east, draw=none},
    legend cell align={left},
]
\addplot+[only marks, color32,mark=*,mark options={fill=color36,scale=0.5}] table [x={Ns}, y={Mean_0_perc=10}] {pd_Section5_Ns_df.dat};
\addplot+[only marks, color36,mark=triangle*,mark options={fill=color36,scale=0.5}] table [x={Ns}, y={Mean_1_perc=10}] {pd_Section5_Ns_df.dat};
\legend{Cars,Trucks}
\end{axis}
\end{tikzpicture}
}
\put(300,100){
\tikzsetnextfilename{Fig_control_Ns_perc5_mean}
\begin{tikzpicture}	
\begin{axis}[
    xlabel={  \# lanes},    ylabel={$\E [T_{1,10}(0)]$ (s)},
    xmin=0.5, xmax=5.5,    ymin=95, ymax=210,
    xtick={1,2,3,4,5},    ytick={100,120,140,160,180,200},
    xticklabel style={/pgf/number format/fixed, /pgf/number format/precision=2},
    yticklabel style={text width=1.1em,align=right},    
	scaled x ticks=false,
    xtick pos=left,    ytick pos=left,    axis lines=left,
    width = 150pt,    height = 100pt,
    x axis line style=-,    y axis line style=-,
    font = \tiny,
    legend style={at={(1.1,0.975)},anchor=north east, draw=none},
    legend cell align={left},
]
\addplot+[only marks, color32,mark=*,mark options={fill=color36,scale=0.5}] table [x={Ns}, y={Mean_0_perc=5}] {pd_Section5_Ns_df.dat};
\addplot+[only marks, color36,mark=triangle*,mark options={fill=color36,scale=0.5}] table [x={Ns}, y={Mean_1_perc=5}] {pd_Section5_Ns_df.dat};
\legend{Cars,Trucks}
\end{axis}
\end{tikzpicture}
}
\put(0,0){
\tikzsetnextfilename{Fig_control_Ns_perc20_std}
\begin{tikzpicture}	
\begin{axis}[
    xlabel={  \# lanes},    ylabel={$\sigma_{T_{1,10}(0)}$ (s)},
    xmin=0.5, xmax=5.5,    ymin=11, ymax=75,
    xtick={1,2,3,4,5},    ytick={20,30,40,50,60,70},
    xticklabel style={/pgf/number format/fixed, /pgf/number format/precision=2},
    yticklabel style={text width=1.1em,align=right},        
	scaled x ticks=false,
    xtick pos=left,    ytick pos=left,    axis lines=left,
    width = 150pt,    height = 100pt,
    x axis line style=-,    y axis line style=-,
    font = \tiny,
    legend style={at={(1.1,0.975)},anchor=north east, draw=none},
    legend cell align={left},
]
\addplot+[only marks, color32,mark=*,mark options={fill=color36,scale=0.5}] table [x={Ns}, y={Std_0_perc=20}] {pd_Section5_Ns_df.dat};
\addplot+[only marks, color36,mark=triangle*,mark options={fill=color36,scale=0.5}] table [x={Ns}, y={Std_1_perc=20}] {pd_Section5_Ns_df.dat};
\legend{Cars,Trucks}
\end{axis}
\end{tikzpicture}
}
\put(150,0){
\tikzsetnextfilename{Fig_control_Ns_perc10_std}
\begin{tikzpicture}	
\begin{axis}[
    xlabel={ \# lanes},    ylabel={$\sigma_{T_{1,10}(0)}$ (s)},
    xmin=0.5, xmax=5.5,    ymin=11, ymax=105,
    xtick={1,2,3,4,5},    ytick={20,40,60,80,100},
    xticklabel style={/pgf/number format/fixed, /pgf/number format/precision=2},
    yticklabel style={text width=1.1em,align=right},    
	scaled x ticks=false,
    xtick pos=left,    ytick pos=left,    axis lines=left,
    width = 150pt,    height = 100pt,
    x axis line style=-,    y axis line style=-,
    font = \tiny,
    legend style={at={(1.1,0.975)},anchor=north east, draw=none},
    legend cell align={left},
]
\addplot+[only marks, color32,mark=*,mark options={fill=color36,scale=0.5}] table [x={Ns}, y={Std_0_perc=10}] {pd_Section5_Ns_df.dat};
\addplot+[only marks, color36,mark=triangle*,mark options={fill=color36,scale=0.5}] table [x={Ns}, y={Std_1_perc=10}] {pd_Section5_Ns_df.dat};
\legend{Cars,Trucks}
\end{axis}
\end{tikzpicture}
}
\put(300,0){
\tikzsetnextfilename{Fig_control_Ns_perc5_std}
\begin{tikzpicture}	
\begin{axis}[
    xlabel={  \# lanes},    ylabel={$\sigma_{T_{1,10}(0)}$ (s)},
    xmin=0.5, xmax=5.5,    ymin=11, ymax=145,
    xtick={1,2,3,4,5},    ytick={20,50,80,110,140},
    xticklabel style={/pgf/number format/fixed, /pgf/number format/precision=2},
    yticklabel style={text width=1.1em,align=right},    
	scaled x ticks=false,
    xtick pos=left,    ytick pos=left,    axis lines=left,
    width = 150pt,    height = 100pt,
    x axis line style=-,    y axis line style=-,
    font = \tiny,
    legend style={at={(1.1,0.975)},anchor=north east, draw=none},
    legend cell align={left},
]
\addplot+[only marks, color32,mark=*,mark options={fill=color36,scale=0.5}] table [x={Ns}, y={Std_0_perc=5}] {pd_Section5_Ns_df.dat};
\addplot+[only marks, color36,mark=triangle*,mark options={fill=color36,scale=0.5}] table [x={Ns}, y={Std_1_perc=5}] {pd_Section5_Ns_df.dat};
\legend{Cars,Trucks}
\end{axis}
\end{tikzpicture}
}
\end{picture}
\caption{Mean ($\E [T_{1,10}(0)]$) and standard deviation ($\sigma_{T_{1,10}}(0)$) of the approximation to the travel-time distribution, plotted as function of $N$ (the number of lanes), for $b=0.2$ (left), $b=0.1$ (middle) and $b=0.05$ (right).}
\label{Fig: control met N}
\end{figure}
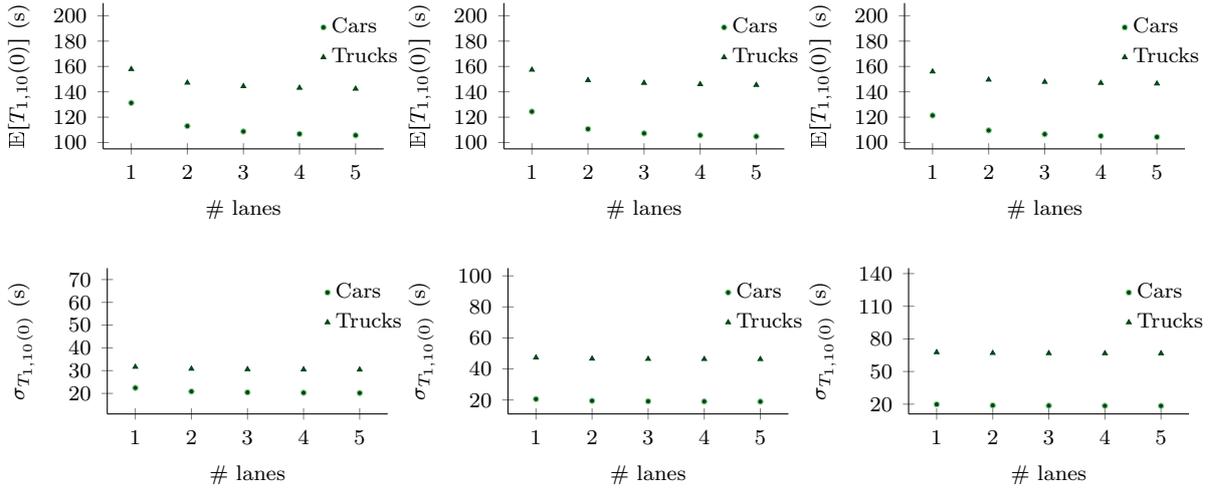

Then we focus in \autoref{Fig: control met arr_rates} on the impact of  the arrival rates. The mean travel-time increases when $\lambda$ increases, which is due to the segment becoming more congested when the arrival rates increases. In addition, the standard deviation is decreasing, which could be explained by the fact that the vehicles drive at lower velocities, and thus create less variability in the travel time. In the upper-right figure, there is a small discrepancy close to zero in the curve for the trucks. This is due to the standard deviation of the travel-time distribution being too large so that there is a significant amount of probability mass below zero (being a artefact of the Gaussian approximation). We corrected for this by only considering the probability mass on the positive part of the distribution's domain and normalizing the corresponding probabilities so that they sum to one. The discrepancy remains however, but can, for practical applications, be remedied by fitting a curve that matches the functional form of the other curves.

\subsection{Extensions}

We now provide a number of extensions of the above experimental setup.
\begin{itemize}
\item In the above experiment, we varied the maximum velocity over an entire segment. This experiment, however, can be easily adapted to a setup in which we vary  the maximum velocity of only a subset of the cells. This is particularly useful, for instance, when evaluating by how much the maximum velocity should be reduced in case of a traffic jam, and what its impact is on parts of the network where the maximum velocity is not reduced.
\item In our experiments we chose the travel-time distribution as the performance metric of interest. Clearly, different performance metrics can be considered, relating to, e.g., congestion, pollution, and emission. In addition, stationary performance measures, as discussed in \autoref{sec: stat.perf.meas}, could be evaluated.
\item There are more control parameters than the ones we considered. For instance, {the modeling framework from \cite{MS2019}} can easily handle varying the per-cell MFD (or, more precisely, the \textit{discrete flux-function}), so as to quantify the effect of e.g., geometric features of the road. 
\item Finally, with respect to control, one could say that our approach is to a large exntent future-proof. If for instance  autonomous vehicles are introduced in today's traffic structure, then these are potentially controlled by artificial intelligence, which, together with interactions with non-autonomous vehicles, inherently leads to stochastic behavior. Therefore, stochastic models will still be required, but with an adapted  shape of the MFD, and with autonomous vehicles included as a specific vehicle type. 
\end{itemize}

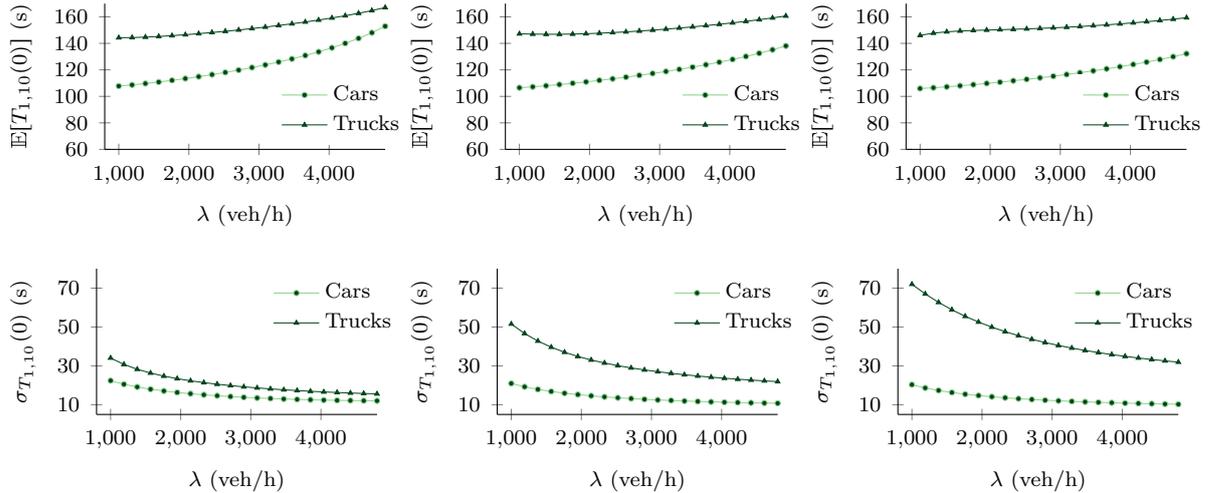
\begin{figure}
\begin{picture}(450 , 200)(0,0)
\put(0,100){
\tikzsetnextfilename{Fig_control_arr_perc20_mean}
\begin{tikzpicture}	
\begin{axis}[
    xlabel={$\lambda$ (veh/h)},    ylabel={$\E [T_{1,10}(0)]$ (s)},
    xmin=800, xmax=4800,    ymin=60, ymax=170,
    xtick={1000,2000,3000,4000},    ytick={60,80,100,120,140,160},
    xticklabel style={/pgf/number format/fixed, /pgf/number format/precision=2},
	scaled x ticks=false,
    xtick pos=left,    ytick pos=left,    axis lines=left,
    width = 150pt,    height = 100pt,
    x axis line style=-,    y axis line style=-,
    font = \tiny,
    legend cell align={left},
    legend style={at={(1.1,0.025)},anchor=south east, draw=none}
]
\addplot+[very thin, color32,mark=*,mark options={fill=color36,scale=0.5}] table [x={Arr_rates}, y={Mean_0_perc=20}] {pd_Section5_Arr_rates_df.dat};
\addplot+[very thin, color36,mark=triangle*,mark options={fill=color36,scale=0.5}] table [x={Arr_rates}, y={Mean_1_perc=20}] {pd_Section5_Arr_rates_df.dat};
\legend{Cars,Trucks}
\end{axis}
\end{tikzpicture}
}
\put(150,100){
\tikzsetnextfilename{Fig_control_arr_perc10_mean}
\begin{tikzpicture}	
\begin{axis}[
    xlabel={$\lambda$ (veh/h)},    ylabel={$\E [T_{1,10}(0)]$ (s)},
    xmin=800, xmax=4800,    ymin=60, ymax=170,
    xtick={1000,2000,3000,4000},    ytick={60,80,100,120,140,160},
    xticklabel style={/pgf/number format/fixed, /pgf/number format/precision=2},
	scaled x ticks=false,
    xtick pos=left,    ytick pos=left,    axis lines=left,
    width = 150pt,    height = 100pt,
    x axis line style=-,    y axis line style=-,
    font = \tiny,
    legend cell align={left},
    legend style={at={(1.1,0.025)},anchor=south east, draw=none}
]
\addplot+[very thin, color32,mark=*,mark options={fill=color36,scale=0.5}] table [x={Arr_rates}, y={Mean_0_perc=10}] {pd_Section5_Arr_rates_df.dat};
\addplot+[very thin, color36,mark=triangle*,mark options={fill=color36,scale=0.5}] table [x={Arr_rates}, y={Mean_1_perc=10}] {pd_Section5_Arr_rates_df.dat};
\legend{Cars,Trucks}
\end{axis}
\end{tikzpicture}
}
\put(300,100){
\tikzsetnextfilename{Fig_control_arr_perc5_mean}
\begin{tikzpicture}	
\begin{axis}[
    xlabel={$\lambda$ (veh/h)},    ylabel={$\E [T_{1,10}(0)]$ (s)},
    xmin=800, xmax=4800,    ymin=60, ymax=170,
    xtick={1000,2000,3000,4000},    ytick={60,80,100,120,140,160},
    xticklabel style={/pgf/number format/fixed, /pgf/number format/precision=2},
	scaled x ticks=false,
    xtick pos=left,    ytick pos=left,    axis lines=left,
    width = 150pt,    height = 100pt,
    x axis line style=-,    y axis line style=-,
    font = \tiny,
    legend cell align={left},
    legend style={at={(1.1,0.025)},anchor=south east, draw=none}
]
\addplot+[very thin, color32,mark=*,mark options={fill=color36,scale=0.5}] table [x={Arr_rates}, y={Mean_0_perc=5}] {pd_Section5_Arr_rates_df.dat};
\addplot+[very thin, color36,mark=triangle*,mark options={fill=color36,scale=0.5}] table [x={Arr_rates}, y={Mean_1_perc=5}] {pd_Section5_Arr_rates_df.dat};
\legend{Cars,Trucks}
\end{axis}
\end{tikzpicture}
}
\put(0,0){
\tikzsetnextfilename{Fig_control_arr_perc20_std}
\begin{tikzpicture}	
\begin{axis}[
    xlabel={$\lambda$ (veh/h)},    ylabel={$\sigma_{T_{1,10}}(0)$ (s)},
    xmin=800, xmax=4800,    ymin=5, ymax=80,
    xtick={1000,2000,3000,4000},    ytick={10,30,50,70},
    xticklabel style={/pgf/number format/fixed, /pgf/number format/precision=2},
	scaled x ticks=false,
    xtick pos=left,    ytick pos=left,    axis lines=left,
    width = 150pt,    height = 100pt,
    x axis line style=-,    y axis line style=-,
    font = \tiny,
    legend cell align={left},
    legend style={at={(1.1,0.975)},anchor=north east, draw=none}
]
\addplot+[very thin, color32,mark=*,mark options={fill=color36,scale=0.5}] table [x={Arr_rates}, y={Std_0_perc=20}] {pd_Section5_Arr_rates_df.dat};
\addplot+[very thin, color36,mark=triangle*,mark options={fill=color36,scale=0.5}] table [x={Arr_rates}, y={Std_1_perc=20}] {pd_Section5_Arr_rates_df.dat};
\legend{Cars,Trucks}
\end{axis}
\end{tikzpicture}
}
\put(150,0){
\tikzsetnextfilename{Fig_control_arr_perc10_std}
\begin{tikzpicture}	
\begin{axis}[
    xlabel={$\lambda$ (veh/h)},    ylabel={$\sigma_{T_{1,10}}(0)$ (s)},
    xmin=800, xmax=4800,    ymin=5, ymax=80,
    xtick={1000,2000,3000,4000},    ytick={10,30,50,70},
    xticklabel style={/pgf/number format/fixed, /pgf/number format/precision=2},
	scaled x ticks=false,
    xtick pos=left,    ytick pos=left,    axis lines=left,
    width = 150pt,    height = 100pt,
    x axis line style=-,    y axis line style=-,
    font = \tiny,
    legend cell align={left},
    legend style={at={(1.1,0.975)},anchor=north east, draw=none}
]
\addplot+[very thin, color32,mark=*,mark options={fill=color36,scale=0.5}] table [x={Arr_rates}, y={Std_0_perc=10}] {pd_Section5_Arr_rates_df.dat};
\addplot+[very thin, color36,mark=triangle*,mark options={fill=color36,scale=0.5}] table [x={Arr_rates}, y={Std_1_perc=10}] {pd_Section5_Arr_rates_df.dat};
\legend{Cars,Trucks}
\end{axis}
\end{tikzpicture}
}
\put(300,0){
\tikzsetnextfilename{Fig_control_arr_perc5_std}
\begin{tikzpicture}	
\begin{axis}[
    xlabel={$\lambda$ (veh/h)},    ylabel={$\sigma_{T_{1,10}}(0)$ (s)},
    xmin=800, xmax=4800,    ymin=5, ymax=80,
    xtick={1000,2000,3000,4000},    ytick={10,30,50,70},
    xticklabel style={/pgf/number format/fixed, /pgf/number format/precision=2},
	scaled x ticks=false,
    xtick pos=left,    ytick pos=left,    axis lines=left,
    width = 150pt,    height = 100pt,
    x axis line style=-,    y axis line style=-,
    font = \tiny,
    legend cell align={left},
    legend style={at={(1.1,0.975)},anchor=north east, draw=none}
]
\addplot+[very thin, color32,mark=*,mark options={fill=color36,scale=0.5}] table [x={Arr_rates}, y={Std_0_perc=5}] {pd_Section5_Arr_rates_df.dat};
\addplot+[very thin, color36,mark=triangle*,mark options={fill=color36,scale=0.5}] table [x={Arr_rates}, y={Std_1_perc=5}] {pd_Section5_Arr_rates_df.dat};
\legend{Cars,Trucks}
\end{axis}
\end{tikzpicture}
}
\end{picture}
\caption{Mean ($\E [T_{1,10}(0)]$) and standard deviation ($\sigma_{T_{1,10}}(0)$) of the approximation to the travel-time distribution, plotted as function of $\lambda$ (the arrival rate), for $b=0.2$ (left), $b=0.1$ (middle) and $b=0.05$ (right).}
\label{Fig: control met arr_rates}
\end{figure}

As a final remark, we want to emphasize that for all control measures and extensions discussed in this section, their implementation in {the} model amounts to an adaptation of (some of) the discrete flux-functions, and/or a change in the model parameters. Therefore, with the methodology that we use it is possible, and in fact quite simple, to evaluate the impact of several implemented control measures simultaneously, instead of quantifying the impact of single control measures in isolation.


\section{A network example}
\label{sec: network}

In this section  we illustrate the use of Gaussian approximations in a network setting. We do so using an example network that is detailed enough to cover all relevant features: there is merging and diverging behavior of traffic, with arrivals and departures at the boundary of the network. We work with a relatively small structure, so as to present the underlying principles as cleanly as possible.  
In this respect, we stress that the same methodology can be used to assess any network topology, under the proviso that one is able to solve the spatial Riemann problems that correspond to the used MFD, yielding the required discrete flux-function. Importantly, these Riemann problems can be restricted to the case where two cells \textit{merge} into one cell, and the case where one cell diverges into two cells, as with these two `building blocks' any  network topology can be constructed, cf.\ \cite[Figs.~2 and 3]{DaganzoNetworks}. 
For the precise description of the underlying construction, including the solution of the Riemann problems pertaining to the Daganzo MFD that we will be using in this section, we refer to~\cite{DaganzoNetworks}.

\begin{figure}
\tikzsetnextfilename{Fig_Network}
\begin{tikzpicture}[
mynode/.style={
  fill=color32,
  fill opacity=0.75,
  minimum size=1em,
  node distance=1.25cm,
  },
mynode2/.style={
  fill=color32,
  fill opacity=0.75,
  minimum size=1em,
  node distance=0.25cm,
  },
mynode3/.style={
  fill=color31,
  fill opacity=0.5,
  minimum size=0.1em,
  node distance=0.25cm,
  }, 
mynode4/.style={
  fill=color34,
  fill opacity=0.5,
  minimum size=0.1em,
  node distance=0.5cm,
  },  
mynode5/.style={
  fill=white,
  fill opacity=0.5,
  minimum size=0.1em,
  node distance=0.5cm,
  },  
every loop/.append style={-latex},  
start chain=going right  
]
\node (s-1) [diamond,mynode5] {};
\node (s1) [circle,mynode4, right=of s-1] {};
\node (s2) [circle,mynode, right=of s1] {};
\node (s3) [circle,mynode, above right=of s2] {}; 
\node (s4) [circle,mynode, below right=of s2] {}; 
\node (s5) [circle,mynode, above right=of s4] {};
\node (s6) [circle,mynode, right=of s5] {}; 
\node (s10) [diamond,mynode5,above=of s3] {};
\node (s11) [diamond,mynode5,below=of s4] {};
\node (s12) [diamond,mynode5,right=of s6] {};

\path[-latex]
  (s-1) edge node[auto,font=\small] {$\lambda$} (s1)
  (s1) edge node[auto,font=\small] {$r_1$} (s2)
  (s2) edge node[auto,font=\small] {$r_2$} (s3)
  (s3) edge node[auto,font=\small] {$r_3$} (s5)    
  (s2) edge node[auto,font=\small,below left] {$r_4$} (s4)
  (s4) edge node[auto,font=\small,below right] {$r_5$} (s5) 
  (s5) edge node[auto,font=\small] {$r_6$} (s6)  
  (s3) edge node[auto,font=\small] {$\nu$} (s10)
  (s4) edge node[auto,font=\small] {$\nu$} (s11)
  (s6) edge node[auto,font=\small] {$\nu$} (s12)
  ;
\end{tikzpicture}
\caption{Schematic representation of the network, consisting of 6 roads $r_i$ between circular nodes, with $d_i \in \N$ cells per road, $i \in \{1,\ldots,6\}$. Arrows with label $\lambda$ or $\nu$ are, respectively, cells at which arrivals and departures occur.}
\label{Fig: Schematic network Sec6}
\end{figure}
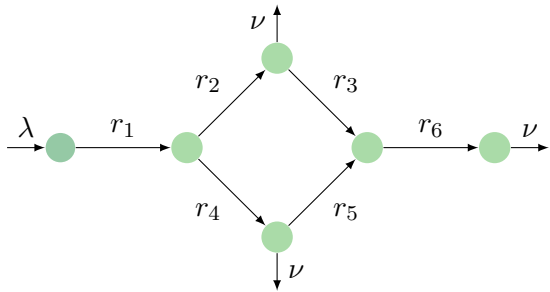

{\autoref{Fig: Schematic network Sec6} visualizes the network we work with in this section. There are six roads, denoted $r_1$ up to $r_6$, where road $k$  consists of $d_k \in \N$ cells. Vehicles arrive at a cell in front of $r_1$ at rate $\lambda$, from which they enter the first cell of $r_1$.
From the last cell of $r_1$, traffic diverges into the first cells of $r_2$ and $r_4$, from which it propagates into these roads. At the last cell of $r_2$, traffic either goes into $r_3$ or to an auxiliary cell from which it leaves the network. The same applies to vehicles leaving the last cell of $r_4$, where traffic goes into $r_5$ or leaves via an auxiliary  cell. From the last cells of $r_3$ and $r_5$, vehicles merge into the first cell of $r_6$. Vehicles then propagate until the last cell of $r_6$, from where they enter an auxiliary cell from where they leave the network.  At any cell where vehicles potentially leave, we let the rate of leaving be  bounded by $\nu$. The padding cells in front of $r_1$ as well as after $r_2$, $r_4$, and $r_6$ have been added to ensure that in no cell merging or diverging simultaneously occurs with arrivals and departures.}


In our setup, we let vehicles diverge according to a fixed routing probability, as in \cite[Section~3.3]{DaganzoNetworks}. In our network traffic diverges at the ends of $r_1$, $r_2$ and $r_4$. Denote by $p_{12}$ the probability of traffic from $r_1$ being routed onto $r_2$, so that the probability of being routed onto $r_4$ is $p_{14}:=1-p_{12}$. The probabilities $p_{23}$  and $p_{45}$ are defined in an analogous manner. 
There is one point in the network where streams merge: at the last cells of $r_3$ and $r_5$, before vehicles are routed onto $r_6$. When the first cell of road~6 cannot accommodate the flows resulting from the last cells of $r_3$ and $r_5$, the parameter $p_{36}$ provides the fraction of available capacity allocated to vehicles stemming from $r_3$ that can enter the first cell of $r_6$  (with  the remainder $p_{56}:=1- p_{36}$ being allocated to vehicles stemming from $r_5$).

Importantly, in principle all parameters can be chosen different from each another, and in addition time-dependent. However, as our experiments are primarily intended to demonstrate our approach, we consider relatively simple instance. In particular, we take time-independent parameters, and use the same $\nu$ at each exit. Moreover, we can work with cell-specific discrete flux-functions, thus facilitating assessing the effect of for instance speed limits in specific cells.

\subsection{Symmetric example}

To demonstrate our approach, we start by considering  a symmetric setting. 
We denote the parameters and of cell~$i$ at road ~$k$ by a subscript~$i$ and superscript~$k$. In particular, the length of cell~$i$ at road~$k$ is denoted $\ell^{(k)}_i$, which in our experiments we take equal to $1$~km, for $i \in \{1,\ldots,d_k\}$, $k\in\{1,\ldots,6\}$. We extend the notation for the vehicle density $\rho_i(\cdot)$ in a similar way to $\rho^{(k)}_i(\cdot)$. 

We choose our parameters  such   that the paths $r_1\to r_2\to r_3\to r_6$ and $r_1\to r_4\to r_5\to r_6$ are equally used: we take $p_{12} = p_{36} = \frac{1}{2}$, and $p_{23} = p_{45} = \frac{3}{4}$.
The other parameters are 
$\lambda = 1800$ veh/h, $\nu = 900$ veh/h. We let each road $k$ consist of $d_k = 5$ cells. 
For the discrete flux-function, we let,  for all cells on all roads,  $v^f = 80$ km/h, $w = 20$ km/h, $\rho^{\textrm{max}} = 108$ veh/km, $q^{\textrm{max}} = 1800$ veh/h.

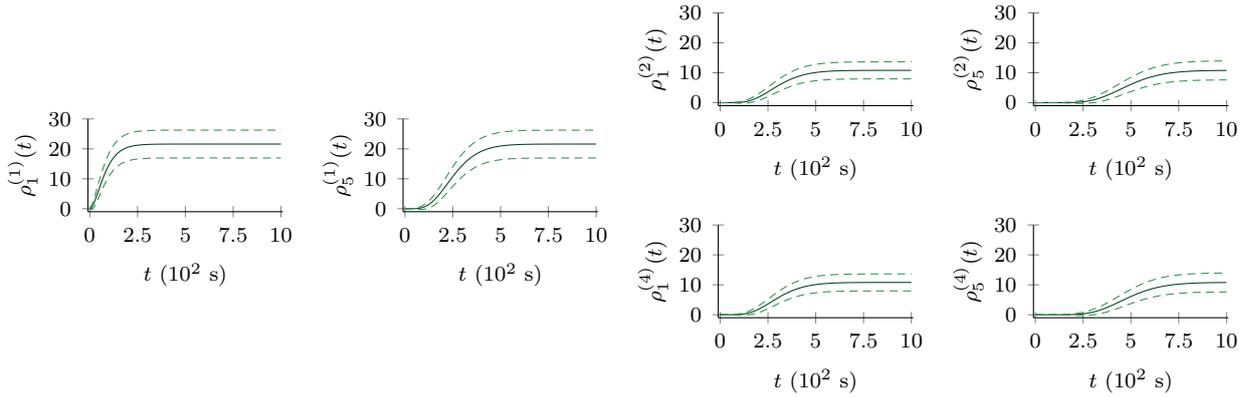
\begin{figure}
\centering
\begin{picture}(472 , 160)(0,0)
\put(0,40){
\tikzsetnextfilename{Fig_sym_cell11}
\begin{tikzpicture}	
\begin{axis}[
    xlabel={$t$ ($10^2$ s)},    ylabel={$\rho^{(1)}_1(t)$},
    xmin=-10, xmax=1010,    ymin=-1, ymax=30,
    xtick={0,250,500,750,1000},    ytick={-0,10,20,30},
    xticklabels={$0$,$2.5$,$5$,$7.5$,$10$},
    xticklabel style={/pgf/number format/fixed, /pgf/number format/precision=2},
	scaled x ticks=false,
    xtick pos=left,    ytick pos=left,    axis lines=left,
    width = 118pt,    height = 80pt,
    y label style={at={(axis description cs:-0.20,.5)}, anchor=south},
    no markers,    x axis line style=-,    y axis line style=-,
    font = \tiny
]
\addplot+[thin, color36] table [x={Time}, y={r1_cell_1}] {df_sec6_network_tau=1000.dat};
\addplot+[densely dashed, color34] table [x={Time}, y={r1_cell_1_upper}] {df_sec6_network_tau=1000.dat};
\addplot+[densely dashed, color34] table [x={Time}, y={r1_cell_1_lower}] {df_sec6_network_tau=1000.dat};
\end{axis}
\end{tikzpicture}
}
\put(118,40){
\tikzsetnextfilename{Fig_sym_cell15}
\begin{tikzpicture}	
\begin{axis}[
    xlabel={$t$ ($10^2$ s)},    ylabel={$\rho^{(1)}_5(t)$},
    xmin=-10, xmax=1010,    ymin=-1, ymax=30,
    xtick={0,250,500,750,1000},    ytick={-0,10,20,30},
    xticklabels={$0$,$2.5$,$5$,$7.5$,$10$},
    xticklabel style={/pgf/number format/fixed, /pgf/number format/precision=2},
	scaled x ticks=false,
    xtick pos=left,    ytick pos=left,    axis lines=left,
    width = 118pt,    height = 80pt,
    y label style={at={(axis description cs:-0.20,.5)}, anchor=south},
    no markers,    x axis line style=-,    y axis line style=-,
    font = \tiny
]
\addplot+[thin, color36] table [x={Time}, y={r1_cell_5}] {df_sec6_network_tau=1000.dat};
\addplot+[densely dashed, color34] table [x={Time}, y={r1_cell_5_upper}] {df_sec6_network_tau=1000.dat};
\addplot+[densely dashed, color34] table [x={Time}, y={r1_cell_5_lower}] {df_sec6_network_tau=1000.dat};
\end{axis}
\end{tikzpicture}
}
\put(236,0){
\tikzsetnextfilename{Fig_sym_cell41}
\begin{tikzpicture}	
\begin{axis}[
    xlabel={$t$ ($10^2$ s)},    ylabel={$\rho^{(4)}_1(t)$},
    xmin=-10, xmax=1010,    ymin=-1, ymax=30,
    xtick={0,250,500,750,1000},    ytick={-0,10,20,30},
    xticklabels={$0$,$2.5$,$5$,$7.5$,$10$},
    xticklabel style={/pgf/number format/fixed, /pgf/number format/precision=2},
	scaled x ticks=false,
    xtick pos=left,    ytick pos=left,    axis lines=left,
    width = 118pt,    height = 80pt,
    y label style={at={(axis description cs:-0.20,.5)}, anchor=south},
    no markers,    x axis line style=-,    y axis line style=-,
    font = \tiny
]
\addplot+[thin, color36] table [x={Time}, y={r4_cell_1}] {df_sec6_network_tau=1000.dat};
\addplot+[densely dashed, color34] table [x={Time}, y={r4_cell_1_upper}] {df_sec6_network_tau=1000.dat};
\addplot+[densely dashed, color34] table [x={Time}, y={r4_cell_1_lower}] {df_sec6_network_tau=1000.dat};
\end{axis}
\end{tikzpicture}
}
\put(354,0){
\tikzsetnextfilename{Fig_sym_cell45}
\begin{tikzpicture}	
\begin{axis}[
    xlabel={$t$ ($10^2$ s)},    ylabel={$\rho^{(4)}_5(t)$},
    xmin=-10, xmax=1010,    ymin=-1, ymax=30,
    xtick={0,250,500,750,1000},    ytick={-0,10,20,30},
    xticklabels={$0$,$2.5$,$5$,$7.5$,$10$},
    xticklabel style={/pgf/number format/fixed, /pgf/number format/precision=2},
	scaled x ticks=false,
    xtick pos=left,    ytick pos=left,    axis lines=left,
    width = 118pt,    height = 80pt,
    y label style={at={(axis description cs:-0.20,.5)}, anchor=south},
    no markers,    x axis line style=-,    y axis line style=-,
    font = \tiny
]
\addplot+[thin, color36] table [x={Time}, y={r4_cell_5}] {df_sec6_network_tau=1000.dat};
\addplot+[densely dashed, color34] table [x={Time}, y={r4_cell_5_upper}] {df_sec6_network_tau=1000.dat};
\addplot+[densely dashed, color34] table [x={Time}, y={r4_cell_5_lower}] {df_sec6_network_tau=1000.dat};
\end{axis}
\end{tikzpicture}
}
\put(236,80){
\tikzsetnextfilename{Fig_sym_cell21}
\begin{tikzpicture}	
\begin{axis}[
    xlabel={$t$ ($10^2$ s)},    ylabel={$\rho^{(2)}_1(t)$},
    xmin=-10, xmax=1010,    ymin=-1, ymax=30,
    xtick={0,250,500,750,1000},    ytick={-0,10,20,30},
    xticklabels={$0$,$2.5$,$5$,$7.5$,$10$},
    xticklabel style={/pgf/number format/fixed, /pgf/number format/precision=2},
	scaled x ticks=false,
    xtick pos=left,    ytick pos=left,    axis lines=left,
    width = 118pt,    height = 80pt,
    y label style={at={(axis description cs:-0.20,.5)}, anchor=south},
    no markers,    x axis line style=-,    y axis line style=-,
    font = \tiny
]
\addplot+[thin, color36] table [x={Time}, y={r2_cell_1}] {df_sec6_network_tau=1000.dat};
\addplot+[densely dashed, color34] table [x={Time}, y={r2_cell_1_upper}] {df_sec6_network_tau=1000.dat};
\addplot+[densely dashed, color34] table [x={Time}, y={r2_cell_1_lower}] {df_sec6_network_tau=1000.dat};
\end{axis}
\end{tikzpicture}
}
\put(354,80){
\tikzsetnextfilename{Fig_sym_cell25}
\begin{tikzpicture}	
\begin{axis}[
    xlabel={$t$ ($10^2$ s)},    ylabel={$\rho^{(2)}_5(t)$},
    xmin=-10, xmax=1010,    ymin=-1, ymax=30,
    xtick={0,250,500,750,1000},    ytick={-0,10,20,30},
    xticklabels={$0$,$2.5$,$5$,$7.5$,$10$},
    xticklabel style={/pgf/number format/fixed, /pgf/number format/precision=2},
	scaled x ticks=false,
    xtick pos=left,    ytick pos=left,    axis lines=left,
    width = 118pt,    height = 80pt,
    y label style={at={(axis description cs:-0.20,.5)}, anchor=south},
    no markers,    x axis line style=-,    y axis line style=-,
    font = \tiny
]
\addplot+[thin, color36] table [x={Time}, y={r2_cell_5}] {df_sec6_network_tau=1000.dat};
\addplot+[densely dashed, color34] table [x={Time}, y={r2_cell_5_upper}] {df_sec6_network_tau=1000.dat};
\addplot+[densely dashed, color34] table [x={Time}, y={r2_cell_5_lower}] {df_sec6_network_tau=1000.dat};
\end{axis}
\end{tikzpicture}
}
\end{picture}
\caption{Mean and variance of vehicle density, in the symmetric setting, for cells~$1$ and $5$ of $r_1$, $r_2$ and $r_4$, for $t \in [0,1000]$.}
\label{Fig: sym network 124}
\end{figure}

\begin{figure}
\centering
\begin{picture}(472 , 160)(0,0)
\put(0,80){
\tikzsetnextfilename{Fig_sym_cell31}
\begin{tikzpicture}	
\begin{axis}[
    xlabel={$t$ ($10^2$ s)},    ylabel={$\rho^{(3)}_1(t)$},
    xmin=-10, xmax=1010,    ymin=-1, ymax=30,
    xtick={0,250,500,750,1000},    ytick={-0,10,20,30},
    xticklabels={$0$,$2.5$,$5$,$7.5$,$10$},
    xticklabel style={/pgf/number format/fixed, /pgf/number format/precision=2},
	scaled x ticks=false,
    xtick pos=left,    ytick pos=left,    axis lines=left,
    width = 118pt,    height = 80pt,
    y label style={at={(axis description cs:-0.20,.5)}, anchor=south},
    no markers,    x axis line style=-,    y axis line style=-,
    font = \tiny
]
\addplot+[thin, color36] table [x={Time}, y={r3_cell_1}] {df_sec6_network_tau=1000.dat};
\addplot+[densely dashed, color34] table [x={Time}, y={r3_cell_1_upper}] {df_sec6_network_tau=1000.dat};
\addplot+[densely dashed, color34] table [x={Time}, y={r3_cell_1_lower}] {df_sec6_network_tau=1000.dat};
\end{axis}
\end{tikzpicture}
}
\put(118,80){
\tikzsetnextfilename{Fig_sym_cell35}
\begin{tikzpicture}	
\begin{axis}[
    xlabel={$t$ ($10^2$ s)},    ylabel={$\rho^{(3)}_5(t)$},
    xmin=-10, xmax=1010,    ymin=-1, ymax=30,
    xtick={0,250,500,750,1000},    ytick={-0,10,20,30},
    xticklabels={$0$,$2.5$,$5$,$7.5$,$10$},
    xticklabel style={/pgf/number format/fixed, /pgf/number format/precision=2},
	scaled x ticks=false,
    xtick pos=left,    ytick pos=left,    axis lines=left,
    width = 118pt,    height = 80pt,
    y label style={at={(axis description cs:-0.20,.5)}, anchor=south},
    no markers,    x axis line style=-,    y axis line style=-,
    font = \tiny
]
\addplot+[thin, color36] table [x={Time}, y={r3_cell_5}] {df_sec6_network_tau=1000.dat};
\addplot+[densely dashed, color34] table [x={Time}, y={r3_cell_5_upper}] {df_sec6_network_tau=1000.dat};
\addplot+[densely dashed, color34] table [x={Time}, y={r3_cell_5_lower}] {df_sec6_network_tau=1000.dat};
\end{axis}
\end{tikzpicture}
}
\put(0,0){
\tikzsetnextfilename{Fig_sym_cell51}
\begin{tikzpicture}	
\begin{axis}[
    xlabel={$t$ ($10^2$ s)},    ylabel={$\rho^{(5)}_1(t)$},
    xmin=-10, xmax=1010,    ymin=-1, ymax=30,
    xtick={0,250,500,750,1000},    ytick={-0,10,20,30},
    xticklabels={$0$,$2.5$,$5$,$7.5$,$10$},
    xticklabel style={/pgf/number format/fixed, /pgf/number format/precision=2},
	scaled x ticks=false,
    xtick pos=left,    ytick pos=left,    axis lines=left,
    width = 118pt,    height = 80pt,
    y label style={at={(axis description cs:-0.20,.5)}, anchor=south},
    no markers,    x axis line style=-,    y axis line style=-,
    font = \tiny
]
\addplot+[thin, color36] table [x={Time}, y={r5_cell_1}] {df_sec6_network_tau=1000.dat};
\addplot+[densely dashed, color34] table [x={Time}, y={r5_cell_1_upper}] {df_sec6_network_tau=1000.dat};
\addplot+[densely dashed, color34] table [x={Time}, y={r5_cell_1_lower}] {df_sec6_network_tau=1000.dat};
\end{axis}
\end{tikzpicture}
}
\put(118,0){
\tikzsetnextfilename{Fig_sym_cell52}
\begin{tikzpicture}	
\begin{axis}[
    xlabel={$t$ ($10^2$ s)},    ylabel={$\rho^{(5)}_5(t)$},
    xmin=-10, xmax=1010,    ymin=-1, ymax=30,
    xtick={0,250,500,750,1000},    ytick={-0,10,20,30},
    xticklabels={$0$,$2.5$,$5$,$7.5$,$10$},
    xticklabel style={/pgf/number format/fixed, /pgf/number format/precision=2},
	scaled x ticks=false,
    xtick pos=left,    ytick pos=left,    axis lines=left,
    width = 118pt,    height = 80pt,
    y label style={at={(axis description cs:-0.20,.5)}, anchor=south},
    no markers,    x axis line style=-,    y axis line style=-,
    font = \tiny
]
\addplot+[thin, color36] table [x={Time}, y={r5_cell_5}] {df_sec6_network_tau=1000.dat};
\addplot+[densely dashed, color34] table [x={Time}, y={r5_cell_5_upper}] {df_sec6_network_tau=1000.dat};
\addplot+[densely dashed, color34] table [x={Time}, y={r5_cell_5_lower}] {df_sec6_network_tau=1000.dat};
\end{axis}
\end{tikzpicture}
}
\put(236,40){
\tikzsetnextfilename{Fig_sym_cell61}
\begin{tikzpicture}	
\begin{axis}[
    xlabel={$t$ ($10^2$ s)},    ylabel={$\rho^{(6)}_1(t)$},
    xmin=-10, xmax=1010,    ymin=-1, ymax=40,
    xtick={0,250,500,750,1000},    ytick={-0,10,20,30},
    xticklabels={$0$,$2.5$,$5$,$7.5$,$10$},
    xticklabel style={/pgf/number format/fixed, /pgf/number format/precision=2},
	scaled x ticks=false,
    xtick pos=left,    ytick pos=left,    axis lines=left,
    width = 118pt,    height = 80pt,
    y label style={at={(axis description cs:-0.20,.5)}, anchor=south},
    no markers,    x axis line style=-,    y axis line style=-,
    font = \tiny
]
\addplot+[thin, color36] table [x={Time}, y={r6_cell_1}] {df_sec6_network_tau=1000.dat};
\addplot+[densely dashed, color34] table [x={Time}, y={r6_cell_1_upper}] {df_sec6_network_tau=1000.dat};
\addplot+[densely dashed, color34] table [x={Time}, y={r6_cell_1_lower}] {df_sec6_network_tau=1000.dat};
\end{axis}
\end{tikzpicture}
}
\put(354,40){
\tikzsetnextfilename{Fig_sym_cell65}
\begin{tikzpicture}	
\begin{axis}[
    xlabel={$t$ ($10^2$ s)},    ylabel={$\rho^{(6)}_5(t)$},
    xmin=-10, xmax=1010,    ymin=-1, ymax=40,
    xtick={0,250,500,750,1000},    ytick={-0,10,20,30},
    xticklabels={$0$,$2.5$,$5$,$7.5$,$10$},
    xticklabel style={/pgf/number format/fixed, /pgf/number format/precision=2},
	scaled x ticks=false,
    xtick pos=left,    ytick pos=left,    axis lines=left,
    width = 118pt,    height = 80pt,
    y label style={at={(axis description cs:-0.20,.5)}, anchor=south},
    no markers,    x axis line style=-,    y axis line style=-,
    font = \tiny
]
\addplot+[thin, color36] table [x={Time}, y={r6_cell_5}] {df_sec6_network_tau=1000.dat};
\addplot+[densely dashed, color34] table [x={Time}, y={r6_cell_5_upper}] {df_sec6_network_tau=1000.dat};
\addplot+[densely dashed, color34] table [x={Time}, y={r6_cell_5_lower}] {df_sec6_network_tau=1000.dat};
\end{axis}
\end{tikzpicture}
}
\end{picture}
\caption{Mean and variance of vehicle density, in the symmetric setting, for cells~$1$ and $5$ of $r_3$, $r_5$ and $r_6$, for $t \in [0,1000]$.}
\label{Fig: sym network 356}
\end{figure}
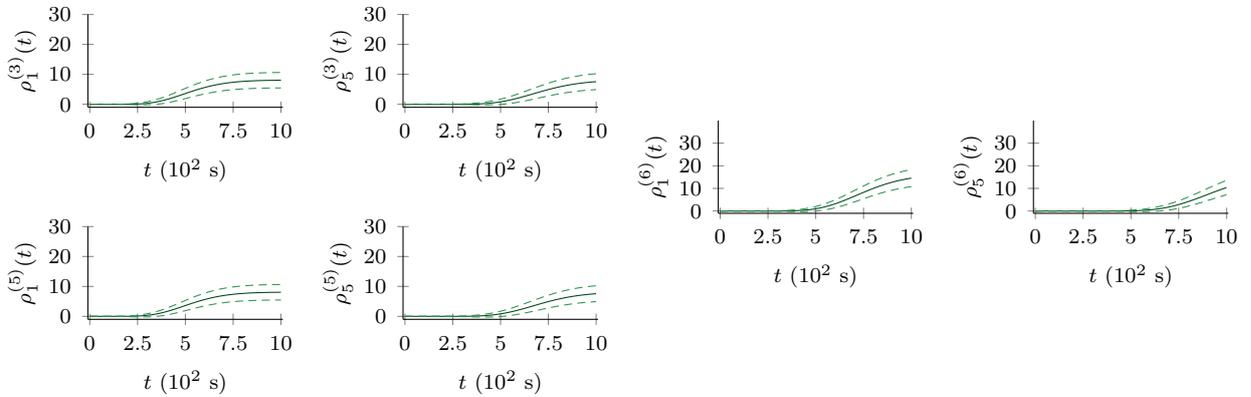

In Figures \ref{Fig: sym network 124} and \ref{Fig: sym network 356} we have plotted the mean of our  Gaussian approximation corresponding to the vehicle densities $\rho^{(k)}_i$, for  the first and last cells of the segments (i.e., $i \in {1,5}$) and all roads (i.e., $k \in \{1,\ldots,6\}$). In addition we show the confidence intervals (whose widths  amount to  two standard deviations) on the time interval $[0,1000]$. In this example we have started from an empty system, but any other starting condition can be dealt with analogously. The mean and standard deviation are found by solving the differential equations that are the network counterpart of \eqref{Eqn: dVar rho}. Observe from the graphs how for each $r_k$, cell~$5$ lags behind cell~$1$ in terms of the build-up of the vehicle density, showing the principle of forward-propagation in this network context. Also note how the mean traffic density $\rho^{(1)}_5(\cdot)$ is correctly divided (according to $p_{12} = \frac{1}{2}$) into two equal parts, with the variance being proportionally split as well. We in addition observe that some  traffic is  departing the system after leaving cell~$5$ of $r_2$ (and cell~$5$ of $r_4$), as can be seen by comparing $\rho^{(2)}_5(\cdot)$ to $\rho^{(3)}_1(\cdot)$.

Figures \ref{Fig: sym2 network 124} and \ref{Fig: sym2 network 356} show the results of the same experiment, but now on the time interval $[0,5400]$. They show  the backpressure that is being exerted throughout the network, starting at the end of $r_6$. Moreover,  the symmetry of $\rho^{(3)}_5(\cdot)$ and $\rho^{(5)}_5(\cdot)$ shows the correct merging behavior:  both roads get equal priority, so that capacity in the first cell or $r_6$ is equally split.

Summarizing, our experiments illustrate that, in our network context, the model is capable of reproducing the proper characteristics. At any point in time, we have the mean and standard deviation of the vehicle densities corresponding to the individual cells, but in addition one can evaluate all covariances (both in the spatial and the temporal sense, that is). As stressed before, having access to these quantities enable the evaluation of a broad range of performance metrics, thus supporting the assessment of the efficacy of various control measures.

\begin{figure}
\centering
\begin{picture}(472 , 160)(0,0)
\put(0,40){
\tikzsetnextfilename{Fig_sym2_cell11}
\begin{tikzpicture}	
\begin{axis}[
    xlabel={$t$ ($10^2$ s)},    ylabel={$\rho^{(1)}_1(t)$},
    xmin=-10, xmax=5410,    ymin=-1, ymax=30,
    xtick={0,1250,2500,3750,5000},    ytick={-0,10,20,30},
    xticklabels={$0$,$12.5$,$25$,$37.5$,$50$},
    xticklabel style={/pgf/number format/fixed, /pgf/number format/precision=2},
	scaled x ticks=false,
    xtick pos=left,    ytick pos=left,    axis lines=left,
    width = 118pt,    height = 80pt,
    y label style={at={(axis description cs:-0.20,.5)}, anchor=south},
    no markers,    x axis line style=-,    y axis line style=-,
    font = \tiny
]
\addplot+[thin, color36] table [x={Time}, y={r1_cell_1}] {df_sec6_network_tau=5400.dat};
\addplot+[densely dashed, color34] table [x={Time}, y={r1_cell_1_upper}] {df_sec6_network_tau=5400.dat};
\addplot+[densely dashed, color34] table [x={Time}, y={r1_cell_1_lower}] {df_sec6_network_tau=5400.dat};
\end{axis}
\end{tikzpicture}
}
\put(118,40){
\tikzsetnextfilename{Fig_sym2_cell15}
\begin{tikzpicture}	
\begin{axis}[
    xlabel={$t$ ($10^2$ s)},    ylabel={$\rho^{(1)}_5(t)$},
    xmin=-10, xmax=5410,    ymin=-1, ymax=30,
    xtick={0,1250,2500,3750,5000},    ytick={-0,10,20,30},
    xticklabels={$0$,$12.5$,$25$,$37.5$,$50$},
    xticklabel style={/pgf/number format/fixed, /pgf/number format/precision=2},
	scaled x ticks=false,
    xtick pos=left,    ytick pos=left,    axis lines=left,
    width = 118pt,    height = 80pt,
    y label style={at={(axis description cs:-0.20,.5)}, anchor=south},
    no markers,    x axis line style=-,    y axis line style=-,
    font = \tiny
]
\addplot+[thin, color36] table [x={Time}, y={r1_cell_5}] {df_sec6_network_tau=5400.dat};
\addplot+[densely dashed, color34] table [x={Time}, y={r1_cell_5_upper}] {df_sec6_network_tau=5400.dat};
\addplot+[densely dashed, color34] table [x={Time}, y={r1_cell_5_lower}] {df_sec6_network_tau=5400.dat};
\end{axis}
\end{tikzpicture}
}
\put(236,0){
\tikzsetnextfilename{Fig_sym2_cell41}
\begin{tikzpicture}	
\begin{axis}[
    xlabel={$t$ ($10^2$ s)},    ylabel={$\rho^{(4)}_1(t)$},
    xmin=-10, xmax=5410,    ymin=-1, ymax=30,
    xtick={0,1250,2500,3750,5000},    ytick={-0,10,20,30},
    xticklabels={$0$,$12.5$,$25$,$37.5$,$50$},
    xticklabel style={/pgf/number format/fixed, /pgf/number format/precision=2},
	scaled x ticks=false,
    xtick pos=left,    ytick pos=left,    axis lines=left,
    width = 118pt,    height = 80pt,
    y label style={at={(axis description cs:-0.20,.5)}, anchor=south},
    no markers,    x axis line style=-,    y axis line style=-,
    font = \tiny
]
\addplot+[thin, color36] table [x={Time}, y={r4_cell_1}] {df_sec6_network_tau=5400.dat};
\addplot+[densely dashed, color34] table [x={Time}, y={r4_cell_1_upper}] {df_sec6_network_tau=5400.dat};
\addplot+[densely dashed, color34] table [x={Time}, y={r4_cell_1_lower}] {df_sec6_network_tau=5400.dat};
\end{axis}
\end{tikzpicture}
}
\put(354,0){
\tikzsetnextfilename{Fig_sym2_cell45}
\begin{tikzpicture}	
\begin{axis}[
    xlabel={$t$ ($10^2$ s)},    ylabel={$\rho^{(4)}_5(t)$},
    xmin=-10, xmax=5410,    ymin=-1, ymax=30,
    xtick={0,1250,2500,3750,5000},    ytick={-0,10,20,30},
    xticklabels={$0$,$12.5$,$25$,$37.5$,$50$},
    xticklabel style={/pgf/number format/fixed, /pgf/number format/precision=2},
	scaled x ticks=false,
    xtick pos=left,    ytick pos=left,    axis lines=left,
    width = 118pt,    height = 80pt,
    y label style={at={(axis description cs:-0.20,.5)}, anchor=south},
    no markers,    x axis line style=-,    y axis line style=-,
    font = \tiny
]
\addplot+[thin, color36] table [x={Time}, y={r4_cell_5}] {df_sec6_network_tau=5400.dat};
\addplot+[densely dashed, color34] table [x={Time}, y={r4_cell_5_upper}] {df_sec6_network_tau=5400.dat};
\addplot+[densely dashed, color34] table [x={Time}, y={r4_cell_5_lower}] {df_sec6_network_tau=5400.dat};
\end{axis}
\end{tikzpicture}
}
\put(236,80){
\tikzsetnextfilename{Fig_sym2_cell21}
\begin{tikzpicture}	
\begin{axis}[
    xlabel={$t$ ($10^2$ s)},    ylabel={$\rho^{(2)}_1(t)$},
    xmin=-10, xmax=5410,    ymin=-1, ymax=30,
    xtick={0,1250,2500,3750,5000},    ytick={-0,10,20,30},
    xticklabels={$0$,$12.5$,$25$,$37.5$,$50$},
    xticklabel style={/pgf/number format/fixed, /pgf/number format/precision=2},
	scaled x ticks=false,
    xtick pos=left,    ytick pos=left,    axis lines=left,
    width = 118pt,    height = 80pt,
    y label style={at={(axis description cs:-0.20,.5)}, anchor=south},
    no markers,    x axis line style=-,    y axis line style=-,
    font = \tiny
]
\addplot+[thin, color36] table [x={Time}, y={r2_cell_1}] {df_sec6_network_tau=5400.dat};
\addplot+[densely dashed, color34] table [x={Time}, y={r2_cell_1_upper}] {df_sec6_network_tau=5400.dat};
\addplot+[densely dashed, color34] table [x={Time}, y={r2_cell_1_lower}] {df_sec6_network_tau=5400.dat};
\end{axis}
\end{tikzpicture}
}
\put(354,80){
\tikzsetnextfilename{Fig_sym2_cell25}
\begin{tikzpicture}	
\begin{axis}[
    xlabel={$t$ ($10^2$ s)},    ylabel={$\rho^{(2)}_5(t)$},
    xmin=-10, xmax=5410,    ymin=-1, ymax=30,
    xtick={0,1250,2500,3750,5000},    ytick={-0,10,20,30},
    xticklabels={$0$,$12.5$,$25$,$37.5$,$50$},
    xticklabel style={/pgf/number format/fixed, /pgf/number format/precision=2},
	scaled x ticks=false,
    xtick pos=left,    ytick pos=left,    axis lines=left,
    width = 118pt,    height = 80pt,
    y label style={at={(axis description cs:-0.20,.5)}, anchor=south},
    no markers,    x axis line style=-,    y axis line style=-,
    font = \tiny
]
\addplot+[thin, color36] table [x={Time}, y={r2_cell_5}] {df_sec6_network_tau=5400.dat};
\addplot+[densely dashed, color34] table [x={Time}, y={r2_cell_5_upper}] {df_sec6_network_tau=5400.dat};
\addplot+[densely dashed, color34] table [x={Time}, y={r2_cell_5_lower}] {df_sec6_network_tau=5400.dat};
\end{axis}
\end{tikzpicture}
}
\end{picture}
\caption{Mean and variance of vehicle density, in the symmetric setting, for cells~$1$ and $5$ of $r_1$, $r_2$ and $r_4$, for $t \in [0,5400]$.}
\label{Fig: sym2 network 124}
\end{figure}
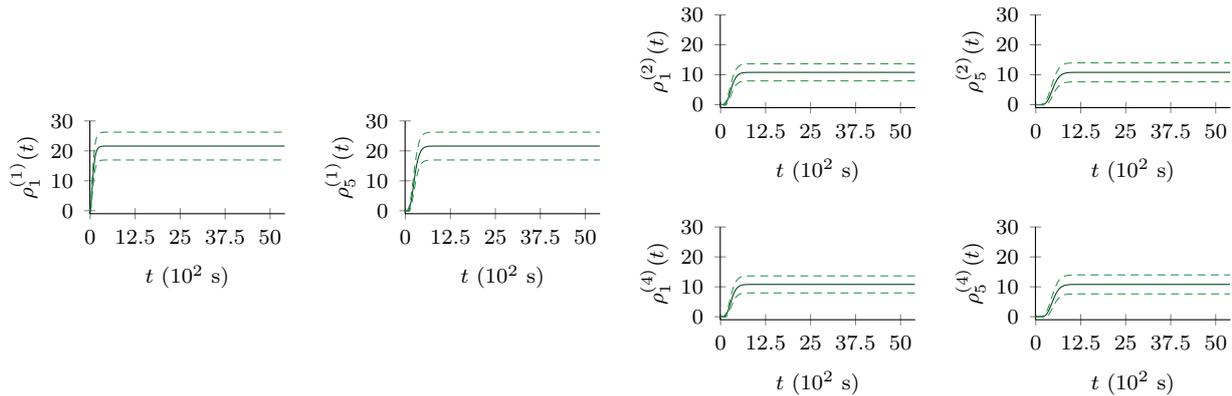

\begin{figure}
\centering
\begin{picture}(472 , 160)(0,0)
\put(0,80){
\tikzsetnextfilename{Fig_sym2_cell31}
\begin{tikzpicture}	
\begin{axis}[
    xlabel={$t$ ($10^2$ s)},    ylabel={$\rho^{(3)}_1(t)$},
    xmin=-10, xmax=5410,    ymin=-1, ymax=30,
    xtick={0,1250,2500,3750,5000},    ytick={-0,10,20,30},
    xticklabels={$0$,$12.5$,$25$,$37.5$,$50$},
    xticklabel style={/pgf/number format/fixed, /pgf/number format/precision=2},
	scaled x ticks=false,
    xtick pos=left,    ytick pos=left,    axis lines=left,
    width = 118pt,    height = 80pt,
    y label style={at={(axis description cs:-0.20,.5)}, anchor=south},
    no markers,    x axis line style=-,    y axis line style=-,
    font = \tiny
]
\addplot+[thin, color36] table [x={Time}, y={r3_cell_1}] {df_sec6_network_tau=5400.dat};
\addplot+[densely dashed, color34] table [x={Time}, y={r3_cell_1_upper}] {df_sec6_network_tau=5400.dat};
\addplot+[densely dashed, color34] table [x={Time}, y={r3_cell_1_lower}] {df_sec6_network_tau=5400.dat};
\end{axis}
\end{tikzpicture}
}
\put(118,80){
\tikzsetnextfilename{Fig_sym2_cell35}
\begin{tikzpicture}	
\begin{axis}[
    xlabel={$t$ ($10^2$ s)},    ylabel={$\rho^{(3)}_5(t)$},
    xmin=-10, xmax=5410,    ymin=-5, ymax=105,
    xtick={0,1250,2500,3750,5000},    ytick={0,25,50,75,100},
    xticklabels={$0$,$12.5$,$25$,$37.5$,$50$},
    xticklabel style={/pgf/number format/fixed, /pgf/number format/precision=2},
	scaled x ticks=false,
    xtick pos=left,    ytick pos=left,    axis lines=left,
    width = 118pt,    height = 80pt,
    y label style={at={(axis description cs:-0.20,.5)}, anchor=south},
    no markers,    x axis line style=-,    y axis line style=-,
    font = \tiny
]
\addplot+[thin, color36] table [x={Time}, y={r3_cell_5}] {df_sec6_network_tau=5400.dat};
\addplot+[densely dashed, color34] table [x={Time}, y={r3_cell_5_upper}] {df_sec6_network_tau=5400.dat};
\addplot+[densely dashed, color34] table [x={Time}, y={r3_cell_5_lower}] {df_sec6_network_tau=5400.dat};
\end{axis}
\end{tikzpicture}
}
\put(0,0){
\tikzsetnextfilename{Fig_sym2_cell51}
\begin{tikzpicture}	
\begin{axis}[
    xlabel={$t$ ($10^2$ s)},    ylabel={$\rho^{(5)}_1(t)$},
    xmin=-10, xmax=5410,    ymin=-1, ymax=30,
    xtick={0,1250,2500,3750,5000},    ytick={-0,10,20,30},
    xticklabels={$0$,$12.5$,$25$,$37.5$,$50$},
    xticklabel style={/pgf/number format/fixed, /pgf/number format/precision=2},
	scaled x ticks=false,
    xtick pos=left,    ytick pos=left,    axis lines=left,
    width = 118pt,    height = 80pt,
    y label style={at={(axis description cs:-0.20,.5)}, anchor=south},
    no markers,    x axis line style=-,    y axis line style=-,
    font = \tiny
]
\addplot+[thin, color36] table [x={Time}, y={r5_cell_1}] {df_sec6_network_tau=5400.dat};
\addplot+[densely dashed, color34] table [x={Time}, y={r5_cell_1_upper}] {df_sec6_network_tau=5400.dat};
\addplot+[densely dashed, color34] table [x={Time}, y={r5_cell_1_lower}] {df_sec6_network_tau=5400.dat};
\end{axis}
\end{tikzpicture}
}
\put(118,0){
\tikzsetnextfilename{Fig_sym2_cell52}
\begin{tikzpicture}	
\begin{axis}[
    xlabel={$t$ ($10^2$ s)},    ylabel={$\rho^{(5)}_5(t)$},
    xmin=-10, xmax=5410,    ymin=-5, ymax=105,
    xtick={0,1250,2500,3750,5000},    ytick={0,25,50,75,100},
    xticklabels={$0$,$12.5$,$25$,$37.5$,$50$},
    xticklabel style={/pgf/number format/fixed, /pgf/number format/precision=2},
	scaled x ticks=false,
    xtick pos=left,    ytick pos=left,    axis lines=left,
    width = 118pt,    height = 80pt,
    y label style={at={(axis description cs:-0.20,.5)}, anchor=south},
    no markers,    x axis line style=-,    y axis line style=-,
    font = \tiny
]
\addplot+[thin, color36] table [x={Time}, y={r5_cell_5}] {df_sec6_network_tau=5400.dat};
\addplot+[densely dashed, color34] table [x={Time}, y={r5_cell_5_upper}] {df_sec6_network_tau=5400.dat};
\addplot+[densely dashed, color34] table [x={Time}, y={r5_cell_5_lower}] {df_sec6_network_tau=5400.dat};
\end{axis}
\end{tikzpicture}
}
\put(236,40){
\tikzsetnextfilename{Fig_sym2_cell61}
\begin{tikzpicture}	
\begin{axis}[
    xlabel={$t$ ($10^2$ s)},    ylabel={$\rho^{(6)}_1(t)$},
    xmin=-10, xmax=5410,    ymin=-5, ymax=105,
    xtick={0,1250,2500,3750,5000},    ytick={0,25,50,75,100},
    xticklabels={$0$,$12.5$,$25$,$37.5$,$50$},
    xticklabel style={/pgf/number format/fixed, /pgf/number format/precision=2},
	scaled x ticks=false,
    xtick pos=left,    ytick pos=left,    axis lines=left,
    width = 118pt,    height = 80pt,
    y label style={at={(axis description cs:-0.20,.5)}, anchor=south},
    no markers,    x axis line style=-,    y axis line style=-,
    font = \tiny
]
\addplot+[thin, color36] table [x={Time}, y={r6_cell_1}] {df_sec6_network_tau=5400.dat};
\addplot+[densely dashed, color34] table [x={Time}, y={r6_cell_1_upper}] {df_sec6_network_tau=5400.dat};
\addplot+[densely dashed, color34] table [x={Time}, y={r6_cell_1_lower}] {df_sec6_network_tau=5400.dat};
\end{axis}
\end{tikzpicture}
}
\put(354,40){
\tikzsetnextfilename{Fig_sym2_cell65}
\begin{tikzpicture}	
\begin{axis}[
    xlabel={$t$ ($10^2$ s)},    ylabel={$\rho^{(6)}_5(t)$},
    xmin=-10, xmax=5410,    ymin=-5, ymax=105,
    xtick={0,1250,2500,3750,5000},    ytick={0,25,50,75,100},
    xticklabels={$0$,$12.5$,$25$,$37.5$,$50$},
    xticklabel style={/pgf/number format/fixed, /pgf/number format/precision=2},
	scaled x ticks=false,
    xtick pos=left,    ytick pos=left,    axis lines=left,
    width = 118pt,    height = 80pt,
    y label style={at={(axis description cs:-0.20,.5)}, anchor=south},
    no markers,    x axis line style=-,    y axis line style=-,
    font = \tiny
]
\addplot+[thin, color36] table [x={Time}, y={r6_cell_5}] {df_sec6_network_tau=5400.dat};
\addplot+[densely dashed, color34] table [x={Time}, y={r6_cell_5_upper}] {df_sec6_network_tau=5400.dat};
\addplot+[densely dashed, color34] table [x={Time}, y={r6_cell_5_lower}] {df_sec6_network_tau=5400.dat};
\end{axis}
\end{tikzpicture}
}
\end{picture}
\caption{Mean and variance of vehicle density, in the symmetric setting, for cells~$1$ and $5$ of $r_3$, $r_5$ and $r_6$, for $t \in [0,5400]$.}
\label{Fig: sym2 network 356}
\end{figure}
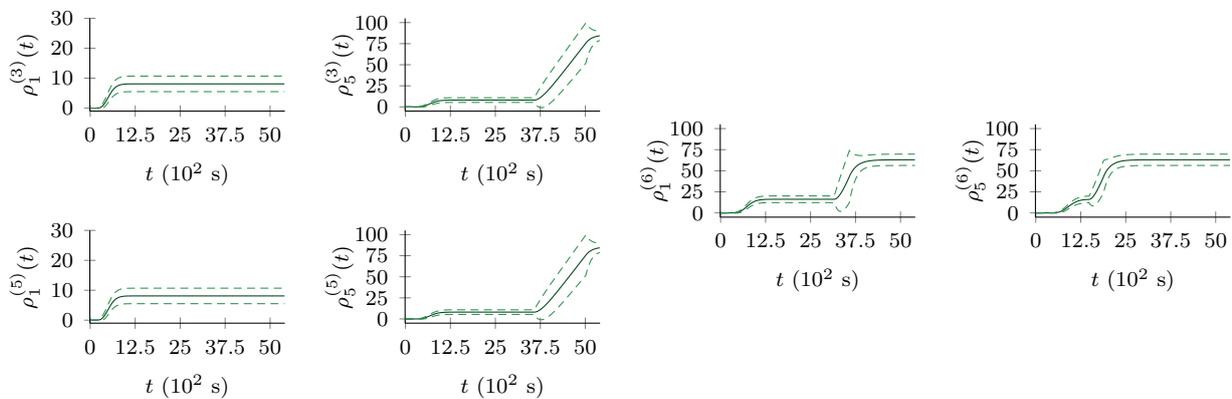

\subsection{Non-symmetric example}

We proceed by an example featuring a non-symmetric network. It shows how one can use the Gaussian approximation to properly quantify the impact of a change in the parameters. We take the network studied in the previous subsection, but we change the maximum velocities of $r_2$, $r_3$, $r_4$, and $r_5$, while varying the routing probability $p_{12}$. More precisely, we increase the maximum velocity of $r_2$ and $r_3$ to $100$ km/h, whereas we lower the maximum velocity of $r_4$ and $r_5$ to $60$ km/h. For these settings, we evaluate the system state for $t \in [0,5400]$, measured in seconds, for $p_{12} \in \{0.25,0.5,0.75\}$, again starting with an empty network at time $0$. Both the  maximum velocity and the routing probability can be considered as control {parameters}; regarding the latter, realize that drivers' routing decisions can be influenced by route-selection software. In a practical setting, one wishes to select the values of the control parameters that optimize a certain objective function, often under some constraints. If one is in a situation in which the objective function can be expressed in terms of the joint vehicle density distribution, one can use our Gaussian methodology to evaluate an approximation of the objective function, which can then be optimized over the control parameter space.

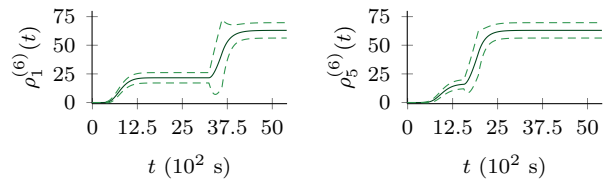
\begin{figure}
\centering
\begin{picture}(236 , 80)(0,0)
\put(0,0){
\tikzsetnextfilename{Fig_Nsym_cell61}
\begin{tikzpicture}	
\begin{axis}[
    xlabel={$t$ ($10^2$ s)},    ylabel={$\rho^{(6)}_1(t)$},
    xmin=-10, xmax=5410,    ymin=-1, ymax=80,
    xtick={0,1250,2500,3750,5000},    ytick={0,25,50,75},
    xticklabels={$0$,$12.5$,$25$,$37.5$,$50$},
    xticklabel style={/pgf/number format/fixed, /pgf/number format/precision=2},
	scaled x ticks=false,
    xtick pos=left,    ytick pos=left,    axis lines=left,
    width = 118pt,    height = 80pt,
    y label style={at={(axis description cs:-0.20,.5)}, anchor=south},
    no markers,    x axis line style=-,    y axis line style=-,
    font = \tiny
]
\addplot+[thin, color36] table [x={Time}, y={r6_cell_1}] {df_sec6_nonsymnetwork2_tau=5400_p_12=0,50.dat};
\addplot+[densely dashed, color34] table [x={Time}, y={r6_cell_1_upper}] {df_sec6_nonsymnetwork2_tau=5400_p_12=0,50.dat};
\addplot+[densely dashed, color34] table [x={Time}, y={r6_cell_1_lower}] {df_sec6_nonsymnetwork2_tau=5400_p_12=0,50.dat};
\end{axis}
\end{tikzpicture}
}
\put(118,0){
\tikzsetnextfilename{Fig_Nsym_cell65}
\begin{tikzpicture}	
\begin{axis}[
    xlabel={$t$ ($10^2$ s)},    ylabel={$\rho^{(6)}_5(t)$},
    xmin=-10, xmax=5410,    ymin=-1, ymax=80,
    xtick={0,1250,2500,3750,5000},    ytick={0,25,50,75},
    xticklabels={$0$,$12.5$,$25$,$37.5$,$50$},
    xticklabel style={/pgf/number format/fixed, /pgf/number format/precision=2},
	scaled x ticks=false,
    xtick pos=left,    ytick pos=left,    axis lines=left,
    width = 118pt,    height = 80pt,
    y label style={at={(axis description cs:-0.20,.5)}, anchor=south},
    no markers,    x axis line style=-,    y axis line style=-,
    font = \tiny
]
\addplot+[thin, color36] table [x={Time}, y={r6_cell_5}] {df_sec6_nonsymnetwork2_tau=5400_p_12=0,50.dat};
\addplot+[densely dashed, color34] table [x={Time}, y={r6_cell_5_upper}] {df_sec6_nonsymnetwork2_tau=5400_p_12=0,50.dat};
\addplot+[densely dashed, color34] table [x={Time}, y={r6_cell_5_lower}] {df_sec6_nonsymnetwork2_tau=5400_p_12=0,50.dat};
\end{axis}
\end{tikzpicture}
}
\end{picture}
\caption{Mean and variance of vehicle density, in the non-symmetric setting with $p_{12} = 0.5$, for cells~$1$ and $5$ of $r_6$, for $t \in [0,5400]$.}
\label{Fig: Nsym network 6}
\end{figure}

\begin{figure}
\centering
\begin{picture}(354 , 160)(0,0)
\put(0,80){
\tikzsetnextfilename{Fig_Nsym_cell35_p12_025}
\begin{tikzpicture}	
\begin{axis}[
    xlabel={$t$ ($10^2$ s)},    ylabel={$\rho^{(3)}_5(t)$},
    xmin=-10, xmax=5410,    ymin=-0.3, ymax=6,
    xtick={0,1250,2500,3750,5000},    ytick={0,2,4,6},
    xticklabels={$0$,$12.5$,$25$,$37.5$,$50$},
    xticklabel style={/pgf/number format/fixed, /pgf/number format/precision=2},
    yticklabel style={text width=1.1em,align=right},
	scaled x ticks=false,
    xtick pos=left,    ytick pos=left,    axis lines=left,
    width = 118pt,    height = 80pt,
    y label style={at={(axis description cs:-0.20,.5)}, anchor=south},
    no markers,    x axis line style=-,    y axis line style=-,
    font = \tiny
]
\addplot+[thin, color36] table [x={Time}, y={r3_cell_5}] {df_sec6_nonsymnetwork2_tau=5400_p_12=0,25.dat};
\addplot+[densely dashed, color34] table [x={Time}, y={r3_cell_5_upper}] {df_sec6_nonsymnetwork2_tau=5400_p_12=0,25.dat};
\addplot+[densely dashed, color34] table [x={Time}, y={r3_cell_5_lower}] {df_sec6_nonsymnetwork2_tau=5400_p_12=0,25.dat};
\end{axis}
\end{tikzpicture}
}
\put(118,80){
\tikzsetnextfilename{Fig_Nsym_cell35_p12_050}
\begin{tikzpicture}	
\begin{axis}[
    xlabel={$t$ ($10^2$ s)},    ylabel={$\rho^{(3)}_5(t)$},
    xmin=-10, xmax=5410,    ymin=-15, ymax=105,
    xtick={0,1250,2500,3750,5000},    ytick={0,30,60,90},
    xticklabels={$0$,$12.5$,$25$,$37.5$,$50$},
    xticklabel style={/pgf/number format/fixed, /pgf/number format/precision=2},
    yticklabel style={text width=1.1em,align=right},
	scaled x ticks=false,
    xtick pos=left,    ytick pos=left,    axis lines=left,
    width = 118pt,    height = 80pt,
    y label style={at={(axis description cs:-0.20,.5)}, anchor=south},
    no markers,    x axis line style=-,    y axis line style=-,
    font = \tiny
]
\addplot+[thin, color36] table [x={Time}, y={r3_cell_5}] {df_sec6_nonsymnetwork2_tau=5400_p_12=0,50.dat};
\addplot+[densely dashed, color34] table [x={Time}, y={r3_cell_5_upper}] {df_sec6_nonsymnetwork2_tau=5400_p_12=0,50.dat};
\addplot+[densely dashed, color34] table [x={Time}, y={r3_cell_5_lower}] {df_sec6_nonsymnetwork2_tau=5400_p_12=0,50.dat};
\end{axis}
\end{tikzpicture}
}
\put(236,80){
\tikzsetnextfilename{Fig_Nsym_cell35_p12_075}
\begin{tikzpicture}	
\begin{axis}[
    xlabel={$t$ ($10^2$ s)},    ylabel={$\rho^{(3)}_5(t)$},
    xmin=-10, xmax=5410,    ymin=-15, ymax=105,
    xtick={0,1250,2500,3750,5000},    ytick={0,30,60,90},
    xticklabels={$0$,$12.5$,$25$,$37.5$,$50$},
    xticklabel style={/pgf/number format/fixed, /pgf/number format/precision=2},
    yticklabel style={text width=1.1em,align=right},    
	scaled x ticks=false,
    xtick pos=left,    ytick pos=left,    axis lines=left,
    width = 118pt,    height = 80pt,
    y label style={at={(axis description cs:-0.20,.5)}, anchor=south},
    no markers,    x axis line style=-,    y axis line style=-,
    font = \tiny
]
\addplot+[thin, color36] table [x={Time}, y={r3_cell_5}] {df_sec6_nonsymnetwork2_tau=5400_p_12=0,75.dat};
\addplot+[densely dashed, color34] table [x={Time}, y={r3_cell_5_upper}] {df_sec6_nonsymnetwork2_tau=5400_p_12=0,75.dat};
\addplot+[densely dashed, color34] table [x={Time}, y={r3_cell_5_lower}] {df_sec6_nonsymnetwork2_tau=5400_p_12=0,75.dat};
\end{axis}
\end{tikzpicture}
}
\put(0,0){
\tikzsetnextfilename{Fig_Nsym_cell55_p12_025}
\begin{tikzpicture}	
\begin{axis}[
    xlabel={$t$ ($10^2$ s)},    ylabel={$\rho^{(5)}_5(t)$},
    xmin=-10, xmax=5410,    ymin=-15, ymax=105,
    xtick={0,1250,2500,3750,5000},    ytick={0,30,60,90},
    xticklabels={$0$,$12.5$,$25$,$37.5$,$50$},
    xticklabel style={/pgf/number format/fixed, /pgf/number format/precision=2},
    yticklabel style={text width=1.1em,align=right},    
	scaled x ticks=false,
    xtick pos=left,    ytick pos=left,    axis lines=left,
    width = 118pt,    height = 80pt,
    y label style={at={(axis description cs:-0.20,.5)}, anchor=south},
    no markers,    x axis line style=-,    y axis line style=-,
    font = \tiny
]
\addplot+[thin, color36] table [x={Time}, y={r5_cell_5}] {df_sec6_nonsymnetwork2_tau=5400_p_12=0,25.dat};
\addplot+[densely dashed, color34] table [x={Time}, y={r5_cell_5_upper}] {df_sec6_nonsymnetwork2_tau=5400_p_12=0,25.dat};
\addplot+[densely dashed, color34] table [x={Time}, y={r5_cell_5_lower}] {df_sec6_nonsymnetwork2_tau=5400_p_12=0,25.dat};
\end{axis}
\end{tikzpicture}
}
\put(118,0){
\tikzsetnextfilename{Fig_Nsym_cell55_p12_050}
\begin{tikzpicture}	
\begin{axis}[
    xlabel={$t$ ($10^2$ s)},    ylabel={$\rho^{(5)}_5(t)$},
    xmin=-10, xmax=5410,    ymin=-15, ymax=105,
    xtick={0,1250,2500,3750,5000},    ytick={0,30,60,90},
    xticklabels={$0$,$12.5$,$25$,$37.5$,$50$},
    xticklabel style={/pgf/number format/fixed, /pgf/number format/precision=2},
    yticklabel style={text width=1.1em,align=right},    
	scaled x ticks=false,
    xtick pos=left,    ytick pos=left,    axis lines=left,
    width = 118pt,    height = 80pt,
    y label style={at={(axis description cs:-0.20,.5)}, anchor=south},
    no markers,    x axis line style=-,    y axis line style=-,
    font = \tiny
]
\addplot+[thin, color36] table [x={Time}, y={r5_cell_5}] {df_sec6_nonsymnetwork2_tau=5400_p_12=0,50.dat};
\addplot+[densely dashed, color34] table [x={Time}, y={r5_cell_5_upper}] {df_sec6_nonsymnetwork2_tau=5400_p_12=0,50.dat};
\addplot+[densely dashed, color34] table [x={Time}, y={r5_cell_5_lower}] {df_sec6_nonsymnetwork2_tau=5400_p_12=0,50.dat};
\end{axis}
\end{tikzpicture}
}
\put(236,0){
\tikzsetnextfilename{Fig_Nsym_cell55_p12_075}
\begin{tikzpicture}	
\begin{axis}[
    xlabel={$t$ ($10^2$ s)},    ylabel={$\rho^{(5)}_5(t)$},
    xmin=-10, xmax=5410,    ymin=-1, ymax=10,
    xtick={0,1250,2500,3750,5000},    ytick={0,3,6,9},
    xticklabels={$0$,$12.5$,$25$,$37.5$,$50$},
    xticklabel style={/pgf/number format/fixed, /pgf/number format/precision=2},
    yticklabel style={text width=1.1em,align=right},    
	scaled x ticks=false,
    xtick pos=left,    ytick pos=left,    axis lines=left,
    width = 118pt,    height = 80pt,
    y label style={at={(axis description cs:-0.20,.5)}, anchor=south},
    no markers,    x axis line style=-,    y axis line style=-,
    font = \tiny
]
\addplot+[thin, color36] table [x={Time}, y={r5_cell_5}] {df_sec6_nonsymnetwork2_tau=5400_p_12=0,75.dat};
\addplot+[densely dashed, color34] table [x={Time}, y={r5_cell_5_upper}] {df_sec6_nonsymnetwork2_tau=5400_p_12=0,75.dat};
\addplot+[densely dashed, color34] table [x={Time}, y={r5_cell_5_lower}] {df_sec6_nonsymnetwork2_tau=5400_p_12=0,75.dat};
\end{axis}
\end{tikzpicture}
}
\end{picture}
\caption{Mean and variance of vehicle density, in the non-symmetric setting with $p_{12} = 0.25$ (left column), $p_{12} = 0.5$ (middle column) and $p_{12} = 0.75$ (right column), for cell~$5$ of $r_3$ (upper row) and $r_5$ (lower row), for $t \in [0,5400]$.}
\label{Fig: Nsym network 35}
\end{figure}
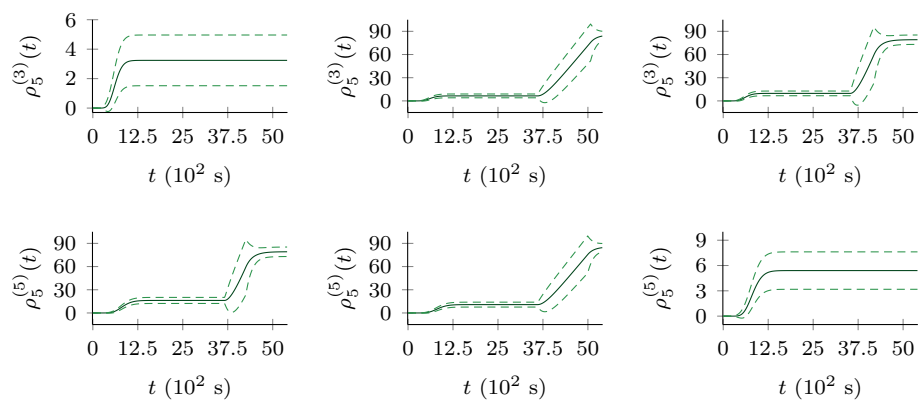

\autoref{Fig: Nsym network 6} shows, for $p_{12} = 0.5$, the mean vehicle density of cells~1 and 5 of $r_6$, on $[0,5400]$, together with a confidence interval (which is again two standard deviations wide). For $p_{12} = 0.25$ and $p_{12} = 0.75$ we obtain almost identical figures, so we have omitted those. \autoref{Fig: Nsym network 35} correspond to  cell~5 of $r_3$ and $r_5$, with $p_{12} \in \{0.25,0.5,0.75\}$.

As expected, when increasing $p_{12}$, the vehicle density on $r_3$ increases and the vehicle density on $r_5$ decreases. By comparing the top left and bottom right figure, one sees that the maximum velocity  has an effect on vehicle densities as well: when  lowering the maximum velocity, the vehicle density increases. This can also be seen, however less clearly, by comparing the top right and bottom left figures, and the figures in the middle column. Another interesting phenomenon is the onset of congestion, due to backpressure from cell~1 in $r_6$. The time of this onset can be controlled for $r_3$ and $r_5$, but also the intensity at which the level of congestion increases, which can be seen from comparing the figures for $r_3$ with $p_{12} = 0.5$ and $p_{12}=0.75$ (and similarly for $r_5$ with $p_{12} = 0.25$ and $p_{12}=0.5$).

Suppose, for instance, that one wishes  to minimize the amount of carbon-dioxide emitted by vehicles (which increases in the vehicles'  velocities), but at the same time  the throughput (which also increases with maximum velocity, but relatively slowly) is to be maximized. For any weight assigned to these two opposite effects, with our approximation one can select the values of $p_{12}$ and the maximum velocities that optimize the resulting objective function. Because we have an approximation of the full distribution of the vehicle densities, our objective function does not necessarily contain mean values only, but can also contain higher moments (and therefore standard deviations) and quantiles. For instance, in the context of the above example, in case one wishes to work with a conservative estimate of the carbon-dioxide emission, one can use its, say, 95th percentile.

\section{Concluding remarks}\label{CR}
In this paper we have demonstrated how the Gaussian process approximation that was developed in  \cite{JL2012,JL2013,MS2019} can be used for the purposes of designing transportation networks and controlling vehicular traffic flows. With a comprehensive analysis of a historical dataset, we have shown that the joint per-cell vehicle density distribution is well approximated by a multivariate Gaussian distribution. In addition, with various  examples, we have illustrated the potential of the Gaussian methodology  for applications in design and control. These examples include (i)~the evaluation of steady-state performance measures, (ii)~route choice based on stochastic travel times, (iii)~various traffic control mechanisms, and (iv)~Gaussian approximations in a network setting, but clearly various other applications can be thought of. The experiments implicitly highlight the apprioriateness of our {stochastic} traffic model over a deterministic counterpart.

This research can be continued in various directions. From an applications point of view, the demonstrated methodologies can be used in various other design-related issues, e.g., evaluating the performance of various network layouts. Moreover, they can be used to develop control algorithms that aim at using real-time information to improve efficiency. In addition, while we can make {the} methodology work for networks of intermediate size, we could explore the computational challenges that arise in relation to the numerical evaluation for large-scale networks.

In terms of more theoretical research themes, there are several promising directions as well. As argued in Section \ref{sec:SPM}, integrals of $\rho(\cdot)$ are relevant in the context of transient performance measures. Noticing now that an integral of a Gaussian process is itself a Gaussian process, one can set up a framework for Gaussian approximations of transient performance measures. In addition, one could focus on developing route-selection algorithms with non-standard cost functions (such as a given percentile of the travel-time distribution, or the mean travel time increased by a given multiple of the corresponding standard deviation). Finally, one could aim at applying our scaling methodology to microscopic or mesoscopic stochastic models, which lend themselves better for describing the specific dynamics of urban environments, e.g., in networks featuring relatively many intersections. 

{\subsubsection*{Acknowledgments}
The authors would like to thank Jesper Slik (Vrije Universiteit Amsterdam) for his help in preparing the historical dataset used in \autoref{sec:DAMV}, and Marko Boon (Eindhoven University of Technology) for his useful feedback.

\bibliographystyle{abbrv} 
\bibliography{segmentbibliography}

\appendix

\newpage

\section{Supplementary material to \autoref{sec:DAMV}}
\label{sec: appendix}

\begin{table}[h]
\centering
\begin{tabular}{|c|c|}
\hline
Index & Name measurement site\\
\hline
\hline
1 & \verb|RWS01_MONIBAS_0021hrl0432ra| \\
2 & \verb|RWS01_MONIBAS_0021hrl0426ra| \\
3 & \verb|RWS01_MONIBAS_0021hrl0420ra| \\
4 & \verb|RWS01_MONIBAS_0021hrl0414ra| \\
5 & \verb|RWS01_MONIBAS_0021hrl0409ra| \\
6 & \verb|RWS01_MONIBAS_0021hrl0403ra| \\
\hline
\end{tabular}
\caption{Measurement site names of the closed-off segment.}
\label{Table: Technical Info}
\end{table}

\begin{figure}[h]
%
\caption{\label{Fig: p-value integrals 1}
Cumulative p-values of $\chi^2$ tests at site~1, as a function of time, for $\tau$-minute intervals, with $\tau=1$ (top left), $\tau =2$ (top right), $\tau=5$ (bottom left), and $\tau = 10$ (bottom right).}
\end{figure}

\begin{figure}[h]
%
\caption{\label{Fig: p-value integrals 3}
Cumulative p-values of $\chi^2$ tests at site~3, as a function of time, for $\tau$-minute intervals, with $\tau=1$ (top left), $\tau =2$ (top right), $\tau=5$ (bottom left), and $\tau = 10$ (bottom right).}
\end{figure}

\begin{figure}
%
\caption{\label{Fig: distr flows 1}Distribution of flows at site 1, at times 05:20 {\sc am} (left) and 07:20 {\sc am}, measured in intervals of $\tau$ minutes, for $\tau = 1,2,5,10,20$ (read from above to below).}
\end{figure}

\begin{figure}
%
\caption{\label{Fig: distr flows 3}Distribution of flows at site 3, at times 05:20 {\sc am} (left) and 07:20 {\sc am}, measured in intervals of $\tau$ minutes, for $\tau = 1,2,5,10,20$ (read from above to below).}
\end{figure}

\begin{figure}
%
\caption{Cumulative p-values of $\chi^2$-test ($P_{(\alpha,\beta)}(t)$), for linear combinations of measurements at site~1 and site~2, as a function of time $t$, for $\tau = 1$.}
\label{Fig: p-value integrals joints, sites 1 and 2, t = 1}
\end{figure}

\begin{figure}
%
\caption{Cumulative p-values of $\chi^2$-test ($P_{(\alpha,\beta)}(t)$), for linear combinations of measurements at site~1 and site~2, as a function of time $t$, for $\tau = 2$.}
\label{Fig: p-value integrals joints, sites 1 and 2, t = 2}
\end{figure}

\begin{figure}
%
\caption{Cumulative p-values of $\chi^2$-test ($P_{(\alpha,\beta)}(t)$), for linear combinations of measurements at site~2 and site~3, as a function of time $t$, for $\tau = 1$.}
\label{Fig: p-value integrals joints, sites 1 and 6, t = 1}
\end{figure}

\begin{figure}
%
\caption{Cumulative p-values of $\chi^2$-test ($P_{(\alpha,\beta)}(t)$), for linear combinations of measurements at site~1 and site~6, as a function of time $t$, for $\tau = 2$.}
\label{Fig: p-value integrals joints, sites 1 and 6, t = 2}
\end{figure}

\end{document}